\documentclass[preprint,12pt,3p]{elsarticle}

\usepackage{amssymb}
\usepackage{amsmath,amsthm}
\usepackage{mathrsfs}
\usepackage{hyperref}
\usepackage{cleveref}
\usepackage[all,cmtip]{xy}
\usepackage{epsf, graphicx}
\usepackage{latexsym,amsfonts,amsbsy,amssymb}
\usepackage{geometry}
\usepackage{titletoc}
\usepackage{color,xcolor}
\usepackage{rotating}
\usepackage{algorithm}
\usepackage{algorithmic}
\usepackage{amscd}
\makeatletter

\@addtoreset{equation}{section} \makeatother
\newtheorem{theorem}{Theorem}[section]
\newtheorem{lemma}[theorem]{Lemma}
\newtheorem{notation}[theorem]{Notation}
\newtheorem{conjecture}[theorem]{Conjecture}
\newtheorem{corollary}[theorem]{Corollary}%[section]
\newtheorem{remark}[theorem]{Remark}%[section]
\newtheorem{void}[theorem]{}
\newtheorem{definition}[theorem]{Definition}%[section]
\newtheorem{proposition}[theorem]{Proposition}%[section]
\def\Ind{{\rm Ind}}
\def\Res{{\rm Res}}

\def\Aut{{\rm Aut}}

\def\Hom{{\rm Hom}}

\def\O{\mathcal{O}}

\def\F{\mathcal{F}}

\def\FF{\mathbb{F}}

\def\Y{\mathcal{Y}}

\def\T{\mathcal{T}}
\def\Z{\mathbb{Z}}

\def\ws{\widetilde{S}}
\def\wa{\widetilde{A}}

\def\0{\bar{0}}
\def\1{\bar{1}}

\def\tr{{\rm tr}}

\makeatletter
\def\ps@pprintTitle{%
\let\@oddhead\@empty
\let\@evenhead\@empty
\def\@oddfoot{\reset@font\hfil\thepage\hfil}
\let\@evenfoot\@oddfoot
}
\makeatother

\begin{document}

\begin{frontmatter}

\title{The refined Brou\'{e} conjecture for RoCK blocks of double covers of symmetric and alternating groups}

\author[label1]{Yucong Du}%\fnref{label3}}  %\corref{cor1}  correspondence auth.
\ead{duyc@cqu.edu.cn}

\author[label2]{Xin Huang}%\fnref{label3}}  %\corref{cor1}  correspondence auth.
\ead{xinhuang@mails.ccnu.edu.cn}

\address[label1]{College of Mathematics and Statistics, Center for Discrete Mathematics, Chongqing University, Chongqing 401331, China}
\address[label2]{School of Mathematics and Statistics, Central China Normal University, Wuhan 430079, China}

\begin{abstract}
Recently, Kleshchev and Livesey \cite{KL22,KL25} proved the existence of RoCK $p$-blocks for double covers of symmetric and alternating groups over large enough coefficient rings. They proved that these RoCK blocks of double covers are Morita equivalent to standard ``local" blocks via bimodules with endopermutation source. Based on this, Kleshchev and Livesey proved that RoCK blocks are splendidly Rickard equivalent to their Brauer correspondents. The analogous result for blocks of symmetric groups, a theorem of Chuang and Kessar \cite{CK}, was an important step in Chuang and Rouquier ultimately proving Brou\'{e}'s abelian defect group conjecture for symmetric groups. In this paper we show that the Morita and splendid Rickard equivalences constructed by Kleshchev and Livesey descend to the ring $\Z_p$ of $p$-adic integers, hence prove Kessar and Linckelmann's refinement of Brou\'{e}'s abelian defect group conjecture for these RoCK blocks.

\end{abstract}


\end{frontmatter}

\begingroup
\renewcommand{\baselinestretch}{0.9}\selectfont
\tableofcontents
\endgroup

% \linenumbers

\section{Introduction}\label{s1}

In \cite{Kessar_Linckelmann}, Kessar and Linckelmann proposed a refined version of Brou\'{e}'s abelian defect group conjecture by considering non splitting coefficient rings:

	\begin{conjecture}[The refined Brou\'{e} conjecture {\cite{Kessar_Linckelmann}}]\label{KLconj}
		For an arbitrary $p$-modular system $(K,\O,k)$ and a block idempotent $b$ of a finite group $G$ over $\O$ with an abelian defect group, there is a splendid Rickard equivalence between $\O Gb$ and its Brauer correspondent.
	\end{conjecture}

Comparing to the refined version, the usual Brou\'{e} conjecture is with the assumption that the complete discrete valuation ring is large enough for $G$. Hence Conjecture \ref{KLconj} contains more local-global information and has more applications; see \cite[Corollary 1.9]{Kessar_Linckelmann}, \cite[Theorem 1.4]{Boltje} and \cite[Theorem 1]{H24b}. Conjecture \ref{KLconj} has been proved for blocks with cyclic defect groups in \cite{Kessar_Linckelmann}, blocks with a Klein four defect group in \cite{H23}, blocks of alternating groups in \cite{H24}, blocks of ${\rm SL}_n(q)$ and ${\rm GL}_n(q)$ in defining characteristic (see \cite{HLZ23}), and unipotent blocks of ${\rm GL}_n(q)$ in non-defining characteristics (see \cite{HLZ24}). Conjecture \ref{KLconj} has also been verified for certain blocks of $p$-nilpotent and $p$-solvable groups in \cite{BY}, \cite{M1} and \cite{M2}. Since any field is a splitting field for symmetric groups, Conjecture \ref{KLconj} - which in the symmetric group case is equivalent to the usual Brou\'{e} conjecture - has been proved by Chuang and Rouquier \cite[Theorem 7.6]{CR}.

For a non-negative integer $n$, denote by $S_n$ the symmetric group on $n$ letters and by $\widetilde{S}^\eta_n$ ($\eta\in \{\pm\}$) one of the two double covers $\ws_n^+$ and $\ws_n^-$ of $S_n$; see \ref{void:double covers}. Denote by $A_n$ the corresponding alternating group and by $\widetilde{A}_n$ its double cover. In this paper we provide additional evidence for Kessar and Linckelmann's refinement of Brou\'{e}'s abelian defect group conjecture by considering blocks of $\widetilde{S}^\eta_n$ and $\widetilde{A}_n$. We first prove the following theorem:

\begin{theorem}[Theorem {\ref{theo:proof of p=2}}]\label{theo:p=2}
	If $p=2$, then Conjecture \ref{KLconj} holds for $\ws^\eta_n$ and $\wa_n$.
\end{theorem}

From now on we assume that $p$ is an odd prime. Let $(K,\O, k)$ be a $p$-modular system. That is, $k$ is a field of characteristic $p$, and $\O$ is a complete discrete valuation ring with quotient field $K$ of characteristic $0$ and residue field $k$.	Denote by $z$ the central element in $\ws^\eta_n$ of order $2$. Since $\O \ws^\eta_n(1+z)/2$ is isomorphic to $\O S_n$ and $\O \widetilde{A}_n(1+z)/2$ is isomorphic to $\O A_n$,  Conjecture \ref{KLconj} is true for blocks of $\O \ws^\eta_n(1+z)/2$ and of $\O \widetilde{A}_n(1+z)/2$. So we restrict our attention in this paper to blocks of $\O \ws^\eta_n(1-z)/2$ and of $\O \widetilde{A}_n(1-z)/2$, which are called the {\it spin blocks of the symmetric and alternating groups}, respectively. For block with trivial defect group, Conjecture \ref{KLconj} is trivially true, hence we only need to consider Conjecture \ref{KLconj} for block with nontrivial defect groups. For the rest of this introduction we assume that when dealing with a finite group $G$, $K$ (and consequently $\O$) contains a primitive $|G|$-th root of unity; following \cite{KL22,KL25} we also assume that $-1$ and $2$ have square roots (say $\sqrt{-1}$ and $\sqrt{2}$) in $K$ (and hence $\O$). Let $\Z_p$ denote the ring of $p$-adic integers, $\mathbb{Q}_p$ the field of $p$-adic numbers and $\FF_p$ the field with $p$ elements. 

 Let $(R,R_p)\in\{(\O,\Z_p),(k,\FF_p)\}$ and let $e_z:=(1-z)/2\in R\ws^\eta_n$. From \cite{KL22,KL25} we see that it is important to consider $R\ws^\eta_n$ as a superalgebra via
$$(R \ws^\eta_n)_{\bar{0}}:=R \widetilde{A}_n~~~{\rm and}~~~(R \ws^\eta_n)_{\bar{1}}:=R (\ws^\eta_n\backslash\widetilde{A}_n).$$
 Recall from \cite{KL25} that the {\it spin superblock algebras} (i.e. the indecomposable superideals) of $R \ws^\eta_ne_z$ are labeled by pairs $(\rho,d)$, where $\rho$ is a $\bar{p}$-core ($p$-bar core) and $d\in \mathbb{N}$ (called the {\it weight} of the block) such that $n=|\rho|+dp$. Denote such a superblock algebra by $B^{\eta,\rho,d}_R$ and denote by $e_{\rho,d}$ the central idempotent of $R \ws^\eta_n$ such that $B^{\eta,\rho,d}_R=R\ws^\eta_ne_{\rho,d}$. We prove (in Proposition \ref{prop:descent of block idempotents} below) that the idempotent $e_{\rho,d}$ is contained in the subring $R_p\wa^\eta_n$. Hence we can define $B^{\eta,\rho,d}_{R_p}:=R_p \ws^\eta_n e_{\rho,d}$. Note that the case $d=0$ corresponds to the defect zero situations, so we assume that $d>0$ in the rest of this section. In this case $B^{\eta,\rho,d}_R$ is in fact a spin block algebra (and not just a superblock algebra) of $R\ws^\eta_n$.
 
We denote by $D^{\rho,d}$ a defect group of $B^{\eta,\rho,d}_R$. The defect group $D^{\rho,d}$ of $B^{\eta,\rho,d}_R$ is abelian if and only if $d<p$. Moreover, if $d<p$, then $D^{\rho,d}\cong (C_p)^{\times d}$. By \cite[Theorem 4.5]{KL25} (or \cite[Propsition 3.16]{K96}), as an element of $R\wa_n$, $e_{\rho,d}$ is a single block unless $d=0$ and $\rho$ is even. So $B^{\eta,\rho,d}_R$ can be regarded as an $R$-superalgebra with even part $(B^{\eta,\rho,d}_R)_{\0}:=R\wa_ne_{\rho,d}$. We denote $R_p\wa_ne_{\rho,d}$ by $(B^{\eta,\rho,d}_{R_p})_{\0}$. It is known that $D^{\rho,d}$ is also a defect group of $(B^{\eta,\rho,d}_R)_{\0}$; see \cite[Proposition 3.16]{K96} or Theorem \ref{thm:An_blocks} below. Let $b^{\eta,\rho,d}_R$ be the Brauer correspondent of $B^{\eta,\rho,d}_R$ in $RN_{\ws^\eta_n}(D^{\rho,d})$, and let $(b^{\eta,\rho,d}_R)_{\0}$ be the Brauer correspondent of $(B^{\eta,\rho,d}_R)_{\0}$ in $RN_{\wa_n}(D^{\rho,d})$. Similarly, we let $b^{\eta,\rho,d}_{R_p}$ be the Brauer correspondent of $B^{\eta,\rho,d}_{R_p}$ in $R_pN_{\ws^\eta_n}(D^{\rho,d})$, and let $(b^{\eta,\rho,d}_{R_p})_{\0}$ be the Brauer correspondent of $(B^{\eta,\rho,d}_{R_p})_{\0}$ in $R_pN_{\wa_n}(D^{\rho,d})$. Note that $(b^{\eta,\rho,d}_R)_{\0}$ (resp. $(b^{\eta,\rho,d}_{R_p})_{\0}$) is indeed the even part of the superalgebra $b^{\eta,\rho,d}_R$ (resp. $b^{\eta,\rho,d}_{R_p}$).
If $d>0$, we have $b^{\eta,\rho,d}_R=RN_{\ws^\eta_n}(D^{\rho,d})e_{\rho,0}$, $(b^{\eta,\rho,d}_R)_{\0}=RN_{\wa_n}(D^{\rho,d})e_{\rho,0}$, $b^{\eta,\rho,d}_{R_p}=R_pN_{\ws^\eta_n}(D^{\rho,d})e_{\rho,0}$ and $(b^{\eta,\rho,d}_{R_p})_{\0}=R_pN_{\wa_n}(D^{\rho,d})e_{\rho,0}$; see Theorem \ref{theo:defect group and Brauer correspondent} (ii) below. 

Let $\T_{n,R}$ be the {\it twisted group superalgebra} defined by Kleshchev and Livesey \cite{KL22,KL25}; see \ref{void:twisted group superalgebra} below. For a superalgebra $A$, the twisted wreath superproduct $A\wr_s \T_{n,R}$ is defined in \cite[\S3.7]{KL25} (see also \ref{void:twisted group superalgebra} below). In order to consider the refined Brou\'{e} conjecture, for $\eta\in \{\pm\}$, we define (in \ref{void:epsilon-twisted group superalgebra}) the {\it $\eta$-twisted group superalgebra} $\T^\eta_{n,R_p}$ and the {\it $\eta$-twisted wreath superproduct} $A\wr_s \T^\eta_{n,R_p}$. %Since $\sqrt{-1}$ is assumed to be in $R$, we prove in Proposition \ref{prop:descent of TnR} that $\T_{n,R}\cong \T^\eta_{n,R}\cong R\otimes_{R_p} \T^\eta_{n,R_p}$ and $A\wr_s\T_{n,R}\cong A\wr_s\T^\eta_{n,R}$ as $R$-superalgebras. Note that if $\sqrt{-1}\notin R_p$, then it seems that $\T^+_{n,R_p}$ and $\T^-_{n,R_p}$ are not isomorphic as $R_p$-algebras.

In \cite{KL22}, Kleshchev and Livesey proved the following structure theorem of $b^{\eta,\rho,d}_{k}$: 

\begin{theorem}[{\cite[Theorem A]{KL22}}]\label{theo:intro Kle Liv Brauer corr over k}
Let $0<d<p$ and $\rho$ be any $\bar{p}$-core. Then
\begin{enumerate}[{\rm (i)}]
	\item $b^{\eta,\rho,d}_k\cong B^{\eta,\rho,0}_k\otimes_k (b^{\eta,\varnothing,1}_k\wr_s \T_{d,k})$ as $k$-superalgebras.
	\item $(b^{\eta,\rho,d}_k)_{\0}\cong (B^{\eta,\rho,0}_k\otimes_k (b^{\eta,\varnothing,1}_k\wr_s \T_{d,k}))_{\0}$ as $k$-algebras.
\end{enumerate}
\end{theorem}

Based on the above theorem we prove a similar structure theorem of $b^{\eta,\rho,d}_{\FF_p}$: 

\begin{theorem}[Theorem \ref{theo: Brauer corr over Fp}]\label{theo:intro Brauer corr over Fp}
	Let $0<d<p$ and $\rho$ be any $\bar{p}$-core. Then
	\begin{enumerate}[{\rm (i)}]
		\item $b^{\eta,\rho,d}_{\FF_p}\cong B^{\eta,\rho,0}_{\FF_p}\otimes_{\FF_p} (b^{\eta,\varnothing,1}_{\FF_p}\wr_s \T^\eta_{d,\FF_p})$ as $\FF_p$-superalgebras.
		\item $(b^{\eta,\rho,d}_{\FF_p})_{\0}\cong (B^{\eta,\rho,0}_{\FF_p}\otimes_{\FF_p} (b^{\eta,\varnothing,1}_{\FF_p}\wr_s \T^\eta_{d,\FF_p}))_{\0}$ as $\FF_p$-algebras.
	\end{enumerate}
\end{theorem}

In \cite{KL22}, Kleshchev and Livesey proved the following analogue of ``the Rickard and Mǎrcuş step for symmetric groups" for double covers of symmetric and alternating groups:
\begin{theorem}[{\cite[Theorem C]{KL22}}]\label{theo:intro Brauer corr splendid over k}
Let $0<d<p$ and $\rho$ be any $\bar{p}$-core. Then 
\begin{enumerate}[{\rm (i)}]
	\item $B^{\eta,\rho,0}_k\otimes_k (b^{\eta,\varnothing,1}_k \wr_s\T_{d,k})$ is derived equivalent to $B^{\eta,\rho,0}_k\otimes_k (B^{\eta,\varnothing,1}_k\wr_s \T_{d,k})$.
	\item $(B^{\eta,\rho,0}_k\otimes_k (b^{\eta,\varnothing,1}_k \wr_s\T_{d,k}))_{\0}$ is derived equivalent to $(B^{\eta,\rho,0}_k\otimes_k (B^{\eta,\varnothing,1}_k\wr_s \T_{d,k}))_{\0}.$
\end{enumerate}
\end{theorem}

We prove a similar result over $\FF_p$:

\begin{theorem}[Theorem \ref{theorem:Brauer corr over Fp splendid}]\label{theorem:introduction Brauer corr over Fp splendid}
	Let $0<d<p$ and $\rho$ be any $\bar{p}$-core. Then 
	\begin{enumerate}[{\rm (i)}]
		\item $B^{\eta,\rho,0}_{\FF_p}\otimes_{\FF_p} (b^{\eta,\varnothing,1}_{\FF_p} \wr_s\T^\eta_{d,\FF_p})$ is splendidly Rickard equivalent to $B^{\eta,\rho,0}_{\FF_p}\otimes_{\FF_p} (B^{\eta,\varnothing,1}_{\FF_p}\wr_s \T^\eta_{d,\FF_p})$.
		\item 	$(B^{\eta,\rho,0}_{\FF_p}\otimes_{\FF_p} (b^{\eta,\varnothing,1}_{\FF_p} \wr_s\T^\eta_{d,\FF_p}))_{\0}$ is splendidly Rickard equivalent to $(B^{\eta,\rho,0}_{\FF_p}\otimes_{\FF_p} (B^{\eta,\varnothing,1}_{\FF_p}\wr_s \T^\eta_{d,\FF_p}))_{\0}$.
	\end{enumerate}
\end{theorem}

See Definitions \ref{defi:splendid over group algebras} and \ref{defi: splendid Rickard equivalences for interior P algebras} below for splendid Rickard equivalences between general interior $P$-algebras. %Our strategy for proving Theorem \ref{theorem:introduction Brauer corr over Fp splendid} is slightly different from Kleshchev and Livesey's proof of Theorem  \ref{theo:intro Brauer corr splendid over k}. 
In \cite{KL22}, Kleshchev and Livesey  defined the notion of a {\it $d$-Rouquier $\bar{p}$-core} for $d\in \mathbb{N}$ (see \ref{void:Rouquier core} below) and proved the following result:

\begin{theorem}[{\cite[Theorem D]{KL22}, \cite[Theorem A and Lemma 9.2]{KL25}}]\label{theo:KLRoCK}
	Let $0<d<p$ and $\rho$ be a $d$-Rouquier $\bar{p}$-core. Then $B^{\eta,\rho,d}_R$ is Morita superequivalent to $B^{\eta,\rho,0}_R\otimes_R (B^{\eta,\varnothing,1}_R\wr_s \T_{d,R})$ via a bimodule $X^R$ with an endopermutation source $U^R$.
\end{theorem}

For the definition of a Morita superequivalence, see \cite[\S3.3]{KL25} or \ref{void: Morita superequivalences} below.  Following \cite{KL25} we call the block $e_{\rho,d}$ of $R\ws^\eta_n$ in Theorem \ref{theo:KLRoCK} a {\it RoCK block} of $\ws^\eta_n$. Since $e_{\rho,d}$ is also a block idempotent of $R\wa_n$, following \cite{KL25} we also called the the block $e_{\rho,d}$ of $R\wa_n$ a {\it RoCK block} of $\wa_n$. Using Theorem \ref{theo:KLRoCK}, Kleshchev and Livesey proved Brou\'{e}'s abelian defect group conjecture for RoCK blocks of $R\ws^\eta_n$ and $R\wa_n$:

\begin{theorem}[{\cite[Theorem B]{KL25}}]\label{theo:KLBroueconjecture}
Let $0<d<p$ and $\rho$ be a $d$-Rouquier $\bar{p}$-core. Then $B^{\eta,\rho,d}_R$ and $(B^{\eta,\rho,d}_R)_{\0}$ are splendidly Rickard equivalent to $b^{\eta,\rho,d}_R$ and $(b^{\eta,\rho,d}_R)_{\0}$ respectively.
\end{theorem}

 In this paper we prove that Theorems \ref{theo:KLRoCK} and \ref{theo:KLBroueconjecture} descend to $R_p$:
 \begin{theorem}[Theorem \ref{theorem: Proof of Theorem main1}]\label{theo:main1}
 	Let $0<d<p$ and $\rho$ be a $d$-Rouquier $\bar{p}$-core. 
 	 Then $B^{\eta,\rho,d}_{R_p}$ is Morita superequivalent to $B^{\eta,\rho,0}_{R_p}\otimes_{R_p}(B^{\eta,\varnothing,1}_{R_p}\wr_s \T^\eta_{d,R_p})$ via a bimodule $X^{R_p}$ with an endopermutation source $U^{R_p}$, such that $R\otimes_{R_p} X^{R_p}$ and $R\otimes_{R_p} U^{R_p}$ are isomorphic to the $X^R$ and $U^R$ in Theorem \ref{theo:KLRoCK}, respectively.	 
 \end{theorem}
 
Corollary \ref{corollary:Morita equ for wa_n} (a corollary of Theorem \ref{theo:main1}) answers a question of Kleshchev and Livesey proposed in \cite[Remark 8.15]{KL25} in the highest generality of their question.
 
\begin{theorem}[Theorem \ref{theo: proof of main2}]\label{theo:main2}
	Let $0<d<p$ and $\rho$ be a $d$-Rouquier $\bar{p}$-core. Then $B^{\eta,\rho,d}_{R_p}$ and $(B^{\eta,\rho,d}_{R_p})_{\0}$ are splendidly Rickard equivalent to $b^{\eta,\rho,d}_{R_p}$ and $(b^{\eta,\rho,d}_{R_p})_{\0}$ respectively. In particular, Conjecture \ref{KLconj} holds for RoCK blocks of $\ws^\eta_n$ and $\wa_n$.
\end{theorem}

Besides Theorems \ref{theo:intro Brauer corr over Fp}, \ref{theorem:introduction Brauer corr over Fp splendid} and \ref{theo:main1}, our proof of Theorem \ref{theo:main2} uses extended versions of known results on Morita equivalences with endopermutation source as presented in Section \ref{s: Morita equivalences with endopermutation source over arbitrary fields}.

Using \cite[Theorem 1.4]{Boltje} and \cite[Theorem 1]{H24b}, we obtain a corollary of Theorem \ref{theo:main2}:

\begin{corollary}
Refinements of the Alperin--McKay conjecture by Isaacs--Navarro, Navarro, Turull and Boltje {\rm (see \cite[(1.5), (1.6), (1.7), (1.8)]{Boltje})} and the blockwise Navarro/Galois Alperin weight conjecture {\rm (see \cite[Conjecture 2]{H24b})} hold for RoCK blocks of $\ws^\eta_n$ and $\wa_n$ with abelian defect groups.
\end{corollary}

The blockwise Navarro Alperin weight conjecture has been proved for all blocks of $\ws^\eta_n$ and $\wa_n$ in \cite{DH} using different method.

For a partition $\lambda$ we denote by $h(\lambda)$ the number of nonzero parts of $\lambda$.
Now by Theorem \ref{theo:main2} and Theorem \ref{theorem:splendid Rickard equ pres Brauer corr} below,
%\ref{theo:intro Brauer corr over Fp}, \ref{theorem:introduction Brauer corr over Fp splendid} and \ref{theo:main1}, 
to prove Conjecture \ref{KLconj} for all spin blocks of symmetric and alternating groups, it suffices to prove the following refined version of the Kessar--Schaps conjecture \cite{KS}:

\begin{conjecture}
Arbitrary spin blocks $B^{\eta,\rho,d}_{R_p}$ and $B^{\eta,\rho',d'}_{R_p}$ of symmetric groups (or spin blocks $(B^{\eta,\rho,d}_{R_p})_{\0}$ and $(B^{\eta,\rho',d'}_{R_p})_{\0}$ of alternating groups) are splendidly Rickard equivalent if and only if $d=d'$ and $|\rho|-h(\rho)\equiv |\rho'|-h(\rho')~({\rm mod}~2)$.
\end{conjecture}

If one replaces ``splendidly Rickard equivalent" with ``derived equivalent" and replaces $R_p$ with $k$, the above conjecture is the usual Kessar--Schaps conjecture. Recently, proofs of the usual Kessar--Schaps conjecture were given in \cite{ELV} and \cite{BK}.

This paper is organised as follows: In Section \ref{section:Preliminaries} we review some relevant results on superalgebras and supermodules which are necessary for later use, and we prove the refined Brou\'{e} conjecture for double covers in characteristic $2$. In Section \ref{s3:descent of twisted wreth superproducts} we prove descent results for twisted group superalgebras and twisted wreath superproducts. In Section \ref{section:Descent and supermodules}, we prove general descent results for supermodules and superequivalences. In Section \ref{section:On Rickard equivalences and splendid complexes}, we formulate convenient reference for some useful results on graded Rickard equivalences and splendid complexes. In Section \ref{section:On spin blocks of ws_n and wa_n} we review some block theory of $\ws^\eta_n$ and $\wa_n$. Based on \cite{KL22}, we prove in Section \ref{section:On Brauer corr} the structure results for Brauer correspondents of spin blocks over $\FF_p$. In Section \ref{s: Morita equivalences with endopermutation source over arbitrary fields}, we extend some known results on Morita equivalences with endopermutation source over splitting fields to arbitrary fields. In Section \ref{section:The Morita superequivalence bimodule}, we analyse the Morita superequivalence bimodule for a RoCK block constructed by Kleshchev and Livesey \cite{KL25}, and then realise this bisupermodule over $\Z_p$; we prove the refined Brou\'{e} conjecture for RoCK spin blocks in the end of Section \ref{section:The Morita superequivalence bimodule}.

\section{Preliminaries}\label{section:Preliminaries}

\begin{void}\label{void:generalities}
{\rm \textbf{Generalities.} Throughout the rest of this paper $p$ is a prime, $\O$ is a complete discrete valuation ring with a residue field $k$ of odd characteristic $p$ and a quotient field $K$ of characteristic $0$; unless otherwise specified, we assume that $p$ is odd. 
Such a triple $(K,\O,k)$ is called a $p$-modular system. (Note that we didn't assume $(K,\O,k)$ to be large enough.) An {\it extension} of the $p$-modular system $(K,\O,k)$ is a $p$-modular system $(K',\O',k')$ such that $\O$ is a subring of $\O'$ with $J(\O)\subseteq J(\O')$. For an algebra $A$, we denote by $A^{\rm op}$ the opposite algebra of $A$, by $A^\times$ the set of units in $A$, by $A$-mod the category of finitely generated $A$-modules, and by $K^b(A)$ the homotopy category of bound complexes over $A$-mod (see e.g. \cite[page 117]{Lin18a}). Unless otherwise stated, modules in this paper are left modules. By a {\it block} of $\O G$ (resp. $kG$) we usually mean an indecomposable direct factor $B$ of $\O G$ (resp. $kG$); some times we also mean the corresponding primitive central idempotent. 
	
Let $R\in\{\O,k\}$ and $A$ a finitely generated $R$-algebra. A {\it projective cover} $(P,\varphi)$ of a finitely generated $A$-module $M$ is a projective $A$-module $P$ and a surjective $A$-module homomorphism $\varphi:P \to M$ such that no proper direct summand of $P$ surjects onto $M$. The assumption on $A$ ensures that projective covers of $M$ always exist and are unique up to isomorphism. The $A$-module $\Omega_A(M):={\rm ker}(\varphi)$ is called the {\it Heller translate} of $M$. We inductively define powers of the Heller translate $\Omega_A^n$, for a positive integer $n$. Finally, $\Omega_A^0(M)$ is defined as the direct sum of all the non-projective indecomposable summands of $M$.

Let $R\subseteq R'$ be two commutative rings. Let $A$ be an $R$-algebra and let $A':=R'\otimes_R A$. Let $\Gamma$ be the group of automorphisms of $R'$ which restricts to the identity map on $R$.  For an $A'$-module
$U'$ and a $\sigma\in \Gamma$, denote by ${}^\sigma U'$ the $A'$-module which is equal to $U'$ as a module over the subalgebra $1\otimes A$ of $A'$, such that $x\otimes a$ acts on $U'$ as $\sigma^{-1}(x) \otimes a$ for all $a \in A$ and $x\in R'$.
The $A'$-module $U'$ is {\it $\Gamma$-stable} if ${}^\sigma U'\cong U'$ for all $\sigma\in \Gamma$. We say that $U'$ is {\it defined over} $R$ or $U'$ {\it descends to} $R$, if there is an $A$-module $U$ such that $U'\cong R'\otimes_R U$.
In this special case $U'$ is $\Gamma$-stable, because for any $\sigma\in\Gamma$, the map sending $x\otimes u$ to $\sigma^{-1}(x)\otimes u$ induces an isomorphism $R'\otimes_R U\cong {}^\sigma(R'\otimes_R U)$, where $u\in U$ and $x\in R'$.
}
\end{void}

\begin{void}
{\rm \textbf{Superalgebras and supermodules.} We review some definitions and notation from \cite[\S3]{KL25}. Let $R$ be a commutative ring. An {\it $R$-superspace} is a finitely generated $R$-module $V$ with decomposition $V=V_{\0}\oplus V_{\1}$. For $\varepsilon\in\Z/2\Z$, we refer to the elements of $V_\varepsilon$ as the homogeneous elements of {\it parity} $\varepsilon$, and write $|v|=\varepsilon$ for $v\in V_\varepsilon$. The subspace $V_{\0}$ (resp. $V_{\1}$) is called the {\it even} (resp. {\it odd}~) part of $V$. An {\it $R$-subsuperspace} $W$ of $V$ is an $R$-submodule of $V$ such that $W = (W \cap V_{\0}) \oplus (W \cap V_{\1})$. Of course we can treat $W$ as an $R$-superspace in its own right. Let $V$ and $W$ be $R$-superspaces. The direct sum $V \oplus W$ is considered as an $R$-superspace in an obvious way. The tensor product $V\otimes_R W$ is considered as an $R$-superspace via $|v\otimes w|:=|v|+|w|$.  A {\it homomorphism} $f:V\to W$ of $R$-superspaces is an $R$-linear map satisfying $f(V_\varepsilon)\subseteq W_\varepsilon$, for $\varepsilon \in \Z/2\Z$.  
	
An {\it$R$-superalgebra} is an $R$-superspace $A$ that is also an $R$-algebra such that $A_\varepsilon A_\delta\subseteq A_{\varepsilon+\delta}$ for any $\varepsilon,\delta\in \Z/2\Z$. A {\it homomorphism} $f:A\to B$ of $R$-superalgebras is an algebra homomorphism that is also a homomorphism of $R$-superspaces. Let $A$ be an $R$-superalgebra. If $A_{\1}\cap A^\times\neq \varnothing$, then $A$ is called a {\it superalgebra with superunit}, and any $u\in A_{\1}\cap A^\times$ is called a {\it superunit}. Any two superunits of $A$ differ from each other by an element of $A_{\bar{0}}^\times$. For any superunit $u$ of $A$, we have $A_{\bar{1}}=uA_{\bar{0}}=A_{\bar{0}}u$.
	
An {\it $A$-supermodule} $M$ is an $A$-module which is also an $R$-superspace such that $A_\varepsilon M_\delta\subseteq M_{\varepsilon+\delta}$ for any $\varepsilon,\delta\in\Z/2\Z$. If $B$ is another $R$-superalgebra, an {\it $A$-$B$-bisupermodule module} $M$ is an $A$-$B$-bimodule which is also an $R$-superspace such that $A_\varepsilon M_\delta\subseteq M_{\varepsilon+\delta}$ and $ M_\delta B_\varepsilon\subseteq M_{\delta+\varepsilon}$ for any $\varepsilon,\delta\in\Z/2\Z$. A {\it homomorphism} $f:M\to N$ of $A$-supermodules is an $A$-module homomorphism that is also a homomorphism of $R$-superspaces; a {\it homomorphism} of bisupermodules is defined similarly.
	
We defined the $R$-superalgebra homomorphism
\begin{equation}\label{equation:involution map}
\sigma=\sigma_A:A\to A,~a\mapsto (-1)^{|a|}a.
\end{equation}
If $2\in R^\times$, the superstructure (i.e. the $\Z/2\Z$-grading on $A$) is completely determined by $\sigma_A$.
The tensor product $A\otimes_R B$ of $R$-superalgebras $A$ and $B$ is considered to be an $R$-superalgebra via
$$
(a\otimes b)(a'\otimes b') = (-1)^{|b| |a'|}aa'\otimes bb'.
$$

Following \cite{KL25}, we define the superalgebra $A^{\rm sop}$ (denoted $A^{\rm op}$ in \cite{KL22}) to be equal to $A$ as an $R$-superspace but with multiplication given by $a\cdot b=(-1)^{|a||b|}ba$, for all $a,b\in A$. %Note that the notation $A^{\rm sop}$ and $A^{\rm op}$ in this paper are different. %We define the superalgebra $A^{\rm op}$ to be equal to $A$ as an $R$-superspace but with multiplication given by $a\cdot b=ba$, for all $a,b\in A$. 
We can view an $(A,B)$-bisupermodule $M$ as an $(A \otimes_R B^{\rm sop})$-supermodule via
\begin{align*}
	(a\otimes b)\cdot m=(-1)^{|b||m|}amb,
\end{align*}
for all $a\in A$, $b\in B$ and $m\in M$. In this way, we identify the notions of an $(A,B)$-bisupermodule and an $(A \otimes_R B^{\rm sop})$-supermodule.

If $V$ is an $R$-superspace we denote by $|V|$ the $R$-space which is $V$ with grading forgotten. In particular, if  
$A$ is an $R$-superalgebra and $M$ is an $A$-supermodule, we have an $R$-algebra $|A|$ and an 
$|A|$-module $|M|$. Similarly, 
if $B$ is another superalgebra and $M$ is an $(A,B)$-bisupermodule, we have an $(|A|,|B|)$-bisupermodule $|M|$. An $A$-supermodule $M$ is called {\it module-indecomposable} if $|M|$ is indecomposable as an $|A|$-module.
An $(A,B)$-bisupermodule $M$ is called {\it module-indecomposable} if $|M|$ indecomposable as an $(|A|,|B|)$-bimodule.

Let $A,B,C$ and $D$ be $R$-superalgebras, $M$ an $(A,C)$-bisupermodule and $N$ a $(B,D)$-bisupermodule. We define the $(A\otimes_R B,C\otimes_R D)$-bisupermodule $M\boxtimes_R N$ to be the $R$-superspace $M\otimes_R N$ with the action 
$$
(a\otimes b)(m\otimes n)=(-1)^{|b||m|}(am\otimes bn),\qquad
(m\otimes n)(c\otimes d)=(-1)^{|c||n|}(mc\otimes nd),
$$
for all $a \in A$, $b \in B$, $c \in C$, $d \in D$, $m \in M$ and $n \in N$. In particular, given an $A$-supermodule $M$ and a $B$-supermodule $N$, we have the $(A\otimes_R B)$-supermodule $M\boxtimes_R N$ with the action 
$
(a\otimes b)(m\otimes n)=(-1)^{|b||m|}(am\otimes bn).$ If $X$ is a complex of $(A,C)$-bisupermodules with differential $\delta=(\delta_n:X_n\to X_{n-1})_{n\in \Z}$ and $Y$ is a complex of $(B,D)$-bisupermodules with differential $\epsilon=(\epsilon_n:Y_n\to Y_{n-1})_{n\in \Z}$, then we can define the complex of $(A\otimes_R B,C\otimes_R D)$-bisupermodules $X\boxtimes_R Y$ with $(X\boxtimes_R Y)_n=\oplus_{i+j=n}(X_i\boxtimes_R Y_j)$ and with differential given by taking the direct sum of the maps $(-1)^i{\rm id}_{X_i}\otimes \epsilon_j:X_i\boxtimes_R Y_j\to X_i\boxtimes_R Y_{j-1}$ and $\delta_i\otimes {\rm id}_{Y_j}:X_i\boxtimes_R Y_j\to X_{i-1}\boxtimes_R Y_j$.

}
\end{void}

\begin{remark}\label{remark: usual algebras as superalgebras}
{\rm Note that any finitely generated algebra (module, bimodule) $A$ can be regarded as a superalgebra (supermodule, bisupermodule) with $A$ as the even part and $0$ as the odd part. So the above box tensor products (as well as most other results in this paper) for superalgebras, supermodules and bisupermodules apply to usual algebras, modules and bimodules. When applying the box tensor products to usual algebras and modules, we write $\otimes$ instead of $\boxtimes$.}
\end{remark}

\begin{void}\label{void:supersymmetric algebras}
{\rm \textbf{Supersymmetric superalgebras.} Let $R$ be a commutative ring and $A$ an $R$-superalgebra. A symmetrizing form on $A$ is a linear form $\tr\in \Hom_R(A,R)$ such that the bilinear form $(a,b):=\tr(ab)$ is a perfect symmetric pairing, where ``symmetric" means that $(a,b)=(b,a)$ for any $a,b\in A$. Following \cite[page 97]{KL25}, a superalgebra $A$ is called {\it supersymmetric} if it is finitely generated projective as an $R$-module and there is a symmetrizing form $\tr$ such that $\tr(A_{\1})=\{0\}$. In this case we call $\tr$ a {\it supersymmetrizing form}. Note that by considering a (non-super) algebra as a superalgebra considered in parity $\0$, we recover the definition of a {\it symmetric algebra} from \cite[Definition 2.11.1]{Lin18a}. 
}	
\end{void}

\begin{void}\label{void: Morita superequivalences}
{\rm \textbf{Morita superequivalences.} Let $R$ be a commutative ring and let $A$, $B$ be $R$-superalgebras. Following \cite{KL25}, a {\it Morita superequivalence} between $A$ and $B$ is a Morita equivalence between $A$ and $B$ induced by an $(A,B)$-bisupermodule $M$ and a $(B,A)$-bisupermodule $N$, i.e. $M\otimes_B N\cong A$ as $(A,A)$-bisupermodules and $N\otimes_A M\cong B$ as $(B,B)$-bisupermodules. Here $M\otimes_B N$ is an $(A,A)$-bisupermodule with
$$(M\otimes_B N)_\varepsilon=\operatorname{span}_R\{m\otimes n\mid m\in M_\alpha,\ n\in N_\beta,\ \alpha+\beta=\varepsilon\},$$
	where $\varepsilon,\alpha,\beta\in\mathbb Z/2\mathbb Z$.
}	
\end{void}

\begin{void}
	{\rm \textbf{Dual supermodules.} Let $R$ be a commutative ring and let $A$ be an $R$-superalgebra and $M$ a finitely generated $A$-supermodule. Following \cite{KL25} we define the {\it dual} of $M$ to be $M^*:=\Hom_R(M, R)$ with the structure of an $R$-superspace through
		\begin{align*}
			{\rm Hom}_R(M, R)_{\varepsilon} := \{f \in {\rm Hom}_R(M, R)\mid f(M_{\varepsilon+\1})=0\}\qquad(\varepsilon\in\Z/2\Z).
		\end{align*}
		We treat $M^*$ as a right $A$-supermodule via $(f\cdot a)(m):=f(am)$, 
		for all $a \in A$, $f \in M^*$ and $m \in M$, hence also as an $A^{\rm sop}$-module via $(a\cdot f)(m):=(-1)^{|a||f|}f(am)$. For more properties of dual supermodules and relevant results on Morita superequivalences, we refer to \cite[\S3.7]{KL25}.
	}
\end{void}

\begin{void}
{\rm \textbf{Crossed superproducts.} Let $R$ be a commutative ring and $G$ a finite group. A {\it $G$-graded crossed $R$-superproduct} is an $R$-superalgebra $A$ with a decomposition $\bigoplus_{g\in G} A_g$ into subsuperspaces such that $A_gA_h \subseteq A_{gh}$, for all $g,h\in G$, and such that, for all $g\in G$, we have $A_g\cap A^{\times} \neq \varnothing$.
	If $A$ is a $G$-graded crossed $R$-superproduct, then so is $A^{\rm sop}$ by defining $(A^{\rm sop})_g=A_{g^{-1}}$, for all $g\in G$.
	Note that $A_{1_G}$ is always a subsuperalgebra of $A$, and $(A_{1_G})^{\rm sop}=(A^{\rm sop})_{1_G}$.
	
%A $G$-graded crossed superproduct $A$ is called {\it supersymmetric} if it has a supersymmetrizing form ${\rm tr}: A \to \O$ that turns $A$ into a supersymmetric algebra and ${\rm tr}(A_g)=0$, for all $g \in G\setminus\{1_G\}$. In particular, ${\rm tr}$ gives $A_{1_G}$ the structure of a supersymmetric algebra.
	
 If $A$ and $B$ are $G$-graded crossed $R$-superproducts, following \cite[equation (3.5)]{KL25} we define
	\begin{equation}\label{equation:diag}
		(A,B)_{G}:= \sum_{g\in G}A_g\otimes_R B_{g^{-1}}
		=
		\sum_{g\in G}A_g\otimes_R (B^{\rm sop})_g  \subseteq A \otimes_R B^{\rm sop}.
	\end{equation}
	The definition of a $G$-graded crossed superproduct ensures that 
	$$A_{1_G}\otimes_R B_{1_G}^{\rm sop}\subseteq (A,B)_{G}\subseteq A \otimes_R B^{\rm sop}$$ 
	are $R$-subsuperalgebras.

}
\end{void}

\begin{void}\label{void:twisted group superalgebra}
	{\rm Let $R\in \{\O,k\}$ and let $n$ be a positive integer. In \cite{KL22,KL25}, Kleshchev and Livesey defined the {\it twisted group superalgebra} $\T_{n,R}$ over $R$ of the symmetric group $S_n$ to be the $R$-superalgebra given by odd generators $\tau_1,\cdots,\tau_{n-1}$ subject to the relations 
		\begin{align*}
			\tau_r^2=1,\quad \tau_r\tau_s = -\tau_s\tau_r\text{ if }|r-s|>1,\quad (\tau_r\tau_{r+1})^3=1.
		\end{align*}
		For each $w\in S_n$, choosing a reduced expression $w=s_{r_1}\cdots s_{r_k}$ in terms of simple transpositions, we define $\tau_w:=\tau_{r_1}\cdots \tau_{r_k}$ (which depends up to a sign on the choice of a reduced expression). It is well-known that   $\{\tau_w\mid w\in S_n\}$ is an $R$-basis of $\T_{n,R}$ (one way to see that is to use the first isomorphism in (\ref{equation:Sn_Tn_isom}) below). 
		
			Given an $R$-superalgebra $A$, Kleshchev and Livesey \cite[Section 5.1.a]{KL22} defined the {\it twisted wreath superproduct} $A\wr_s \T_{n,R}$  as the free product $A^{\otimes n}\star \T_{n,R}$ of superalgebras subject to the relations
		\begin{align}\label{algn:A_swr_Tn^R wr}
			\tau_r\,(a_1\otimes\dots\otimes a_n) = (-1)^{\sum_{u\neq r,r+1}|a_u|}\left({}^{s_r}(a_1\otimes\dots\otimes a_n)\,\tau_r\right)
		\end{align}
		for all $a_1,\cdots,a_n\in A$ and $1\leq r \leq n-1$, where $s_r$ is the transposition $(r,r+1)\in S_n$.
	
	}
\end{void}

\begin{proposition}\label{proposition:AwrT-graded}
Let $R\in \{\O,k\}$, $n$ a positive integer and $A$ an $R$-superalgebra. Then 
	\begin{align}\label{algn:AwrT_graded}
	A\wr_s \T_{n,R} = \bigoplus_{w \in S_n} A^{\otimes n}\tau_w,
	\end{align}
where $A^{\otimes n}$ denotes $\underbrace{A\otimes_R\cdots\otimes_R A}_n$. Hence $A\wr_s \T_{n,R}$ is an $S_n$-graded crossed $R$-superproduct. 
\end{proposition}
\begin{proof} 
	Let $(K',\O',k')$ be a $p$-modular system which is an extension of $(K,\O,k)$ such that $\sqrt{-2}$ (a root of the polynomial $x^2+2\in \O[x]$) is contained in $\O'$. Let $R'\in\{\O',k'\}$ and let $A'=R'\otimes_R A$. Then by \cite[equation (3.8)]{KL25}, $A'\wr_s \T_{n,R'}$ is an $S_n$-graded crossed $R'$-superproduct via
\begin{align}\label{algn:AwrT_graded-R'}
	A'\wr_s \T_{n,R'} = \bigoplus_{w \in S_n} (A')^{\otimes n}\tau_w.
\end{align}
Let $B:=A\wr_s \T_{n,R}$, and for any $w\in S_n$, let $B_w$ be the $R$-subspace $A^{\otimes n}\tau_w$ of $B$. Since $\tau_w\in B_w\cap B^\times$, $B_w\cap B^\times\neq \varnothing$. By the relations (\ref{algn:A_swr_Tn^R wr}), we see that $B_{w_1}B_{w_2}\subseteq B_{w_1w_2}$ for any $w_1,w_2\in S_n$. Hence by the definition of a free product, we obtain that $B=\sum_{w\in S_n} B_w$. Since the right side of (\ref{algn:AwrT_graded-R'}) is a direct sum, the sum $\sum_{w\in S_n} B_w$ is a direct sum as well. Hence $B=\bigoplus_{w\in S_n} B_w$, and this completes the proof.
\end{proof}

%\begin{remark}
%{\rm To explain the equality (\ref{algn:AwrT_graded-R'}), \cite[equation (3.8)]{KL25} uses the element $\sqrt{-2}$; see \cite[Proposition 3.22]{KL25}. But after the explanation, %$\sqrt{-2}$ can be throw away - because $\sqrt{-2}$ would not involved in this equality. 
%}
%\end{remark}

\begin{void}
{\rm Keep the notation of Proposition \ref{proposition:AwrT-graded}. Given the $S_n$-grading in (\ref{algn:AwrT_graded}), recall the definition of $(A \wr_s \T_{n,R},A \wr_s \T_{n,R})_{S_n}$ from (\ref{equation:diag}).
Let $M$ be an $(A,A)$-bisupermodule. Following \cite[equation (3.9)]{KL25}, we denote by $M^{\boxtimes n}_{S_n}$ the $(A \wr_s \T_{n,R},A \wr_s \T_{n,R})_{S_n}$-supermodule given by extending the $(A^{\otimes n},A^{\otimes n})$-bisupermodule $M^{\boxtimes n}:=\underbrace{M\boxtimes_R \cdots\boxtimes_R M}_n$, via
	\begin{align}\label{algn:ext_M^n}
		(\tau_r\otimes \tau_r^{-1})\cdot(m_1\otimes\cdots\otimes m_n)
		:=(-1)^{|m_r|+|m_{r+1}|}\left({}^{s_r}\,(m_1\otimes \cdots \otimes m_n)\right),
	\end{align}
	for all $m_1,\cdots,m_n\in M$ and $1\leq r\leq n-1$.
	
Following \cite[equation (3.10)]{KL25} we now define
\begin{align}\label{algn:sup_wreath_bimod}
		M \wr_s \T_{n,R}:= {\rm Ind}_{(A \wr_s \T_{n,R}, A \wr_s \T_{n,R})_{S_n}}^{(A \wr_s \T_{n,R}) \otimes_R (A \wr_s \T_{n,R})^{\rm sop}}(M^{\boxtimes n}_{S_n}).
	\end{align}
%	Note that
%	$$
%	(A \wr_s \T_{n,R}) \otimes_R (A \wr_s \T_{n,R})^{\rm sop} = \bigoplus_{w \in S_d} (\tau_w \otimes 1)(A \wr_s \T_{n,R},A \wr_s \T_{n,R})_{S_n}.
%	$$
%	Therefore, using the bimodule notation, we can write
%	\begin{align}\label{algn:M_wr_Td_sum}
%		M \wr_s \T_{n,R} = \bigoplus_{w \in S_n}\tau_w M^{\boxtimes n}.
%	\end{align}
%	In particular,
%	\begin{align}\label{M_wr_Td_left}
%		M \wr_s \T_{n,R} \cong (A \wr_s \T_{n,R}) \otimes_{A^{\otimes n}} M^{\boxtimes n},
%	\end{align}
%	as $(A \wr_s \T_{n,R},A^{\otimes n})$-bisupermodules.
		
	}
\end{void}

\begin{void}
{\rm  Let $R\in\{\O,k\}$, $n$ a positive integer and $A$ an $R$-superalgebra with a superunit $u$. If $C_2$ is a cyclic group of order $2$ with generator $g$, we can give $A$ the structure of a $C_2$-graded crossed $R$-superproduct  with graded components 
	$A_{g^{\varepsilon}}=A_{\0}u^{\varepsilon}=A_\varepsilon$
	for $\varepsilon\in\Z/2\Z$. 
	For any $n\in\Z_{>0}$, following \cite[5.1.c]{KL22}, we consider the wreath product 
	$$C_2\wr S_n=\{g_1^{\varepsilon_1}\cdots g_n^{\varepsilon_n}w\mid \varepsilon_1,\dots,\varepsilon_n\in\Z/2\Z,\, w\in S_n\}$$
	where $g_r$ is the generator of the $r$-th $C_2$ factor in the base group $C_2^{\times n}$. Then we can give $A\wr_s \T_{n,R}$ the structure of a crossed product of $A_{\0}^{\otimes n}$ and $C_2 \wr S_n$ (see \ref{void:grad algebras}) with 
	graded components 
	$$
	(A\wr_s \T_{n,R})_{g_1^{\varepsilon_1}\cdots g_n^{\varepsilon_n}w}=
	A_{\0}^{\otimes n}u_1^{\varepsilon_1}\dots u_n^{\varepsilon_n} \tau_w
	$$
	where each $u_r$ is $1^{\otimes (r-1)}\otimes u\otimes 1^{\otimes (n-r)}$.
	
	Let $A$ now be a superalgebra with a superunit $u_A$ and 
	$B$ be a superalgebra with a superunit $u_B$. 
	Recalling (\ref{equation:diag}), we have the subsuperalgebras  
	\begin{align*}
		&(A,B)_{C_2}=(A_{\0}\otimes_R B_{\0}^{\rm sop})\oplus (A_{\bar{1}}\otimes_R B_{\bar{1}}^{\rm sop})\subseteq A\otimes_R B^{\rm sop},\\
		&(A\wr_s \T_{n,R},B\wr_s \T_{n,R})_{C_2\wr S_n}\subseteq (A\wr_s \T_{n,R})\otimes_R (B\wr_s \T_{n,R})^{\rm sop}.
	\end{align*}
Setting $U_r:=(u_A)_r \otimes (u_B^{-1})_r$, for $1\leq r\leq n$, we have that $A_{\0}^{\otimes n}\otimes_R (B_{\0}^{\otimes n})^{\rm sop}$ together with the $U_r$'s generate a subsuperalgebra of $(A\wr_s \T_{n,R},B\wr_s \T_{n,R})_{C_2\wr_s S_n}$ isomorphic to $((A,B)_{C_2})^{\otimes n}$. We identify the two superalgebras. 
Then $(A\wr_s \T_{n,R},B\wr_s \T_{n,R})_{C_2\wr_s S_n}$ is generated by $((A,B)_{C_2})^{\otimes n}$ and $\{T_r:=\tau_r \otimes \tau_r^{-1} \mid 1\leq r\leq n-1\}$ subject to the relations listed in \cite[(5.1.13)]{KL22}.

}	
\end{void}

\begin{proposition}[a slight generalisation of {\cite[Proposition 5.1.14]{KL22}}]\label{prop:ext_der_wreath}
	Let $R\in\{\O,k\}$, let $n$ be a positive integer, and let $A$ and $B$ be $R$-superalgebras with superunit. If $X$ is a complex of $A_{\0}$-$B_{\0}$-bimodules that extends to a complex of $(A,B)_{C_2}$-modules, then $X^{\otimes n}$ extends to a complex of $(A\wr_s \T_{n,R},B\wr_s \T_{n,R})_{C_2\wr S_n}$-modules.
\end{proposition}

\begin{proof} 
	In \cite{KL22}, the ground field is assumed to be algebraically closed. But all the coefficients used in the proof of \cite[Proposition 5.1.14]{KL22} are in the prime field, so the ``algebraically closed" assumption is not necessary for \cite[Proposition 5.1.14]{KL22}. Hence the statement holds for $R=k$. The proof for the case $R=\O$ is identical.  \end{proof}

\begin{void}\label{void:double covers}
{\rm \textbf{The groups $\ws^\eta_n$ and $\wa_n$.} The double covers of $S_n$ (where $n$ is a non-negative integer) are given by 
	$$\ws_n^+:=\langle z,t_1,\cdots,t_{n-1}~|~z^2=1,~t_iz=zt_i,~t_i^2=1,~t_it_j=zt_jt_i~{\rm if}~|i-j|>1,~(t_it_{i+1})^3=1 \rangle,$$
$$\ws_n^-:=\langle z,t_1,\cdots,t_{n-1}~|~z^2=1,~t_iz=zt_i,~t_i^2=z,~t_it_j=zt_jt_i~{\rm if}~|i-j|>1,~(t_it_{i+1})^3=1 \rangle.$$ 
We will use $\ws^\eta_n$ ($\eta\in\{\pm\}$) to simultaneously represent both double covers. Denote by $\pi_n$ the surjective group homomorphism $\ws^\eta_n\to S_n$, $z\mapsto 1$, $t_i\mapsto s_i$, where $s_i$ denotes the transposition $(i,i+1)\in S_n$. Then the double cover of $A_n$ is given by $\wa_n:=\pi_n^{-1}(A_n)$.

}	
\end{void}

\begin{theorem}\label{theo:proof of p=2}
If $p=2$, then Conjecture \ref{KLconj} holds for $\ws^\eta_n$ and $\wa_n$.
\end{theorem}

\begin{proof}
Let $(K',\O',k')$ be a $2$-modular system which is an extension of $(K,\O,k)$ such that $K'$ (and hence $\O'$) contains a primitive $(2n!)$-th root of unity.

$\bullet$ Assume that $b$ is a block of $\O \ws^\eta_n$ with an abelian defect group $P$. Let $b'$ be a block of $\O'\ws_n^\eta$ such that $bb'\neq 0$. Then $P$ is also a defect group of $b'$; see e.g. \cite[Lemma 4.3 (iii)]{H26}. Since $z\in Z(\ws^\eta_n)$ and $\langle z\rangle$ is a $2$-group, $\langle z\rangle$ is contained in $P$; see e.g. \cite[Theorem 6.2.6 (i)]{Lin18b}. Moreover, by \cite[Chap. 5, Theorems 8.10, 8.11]{NT}, the image $\bar{b}'$ of $b'$ under the canonical homomorphism $\O' \ws^\eta_n\to \O' S_n$ is a block of $\O' S_n$ with a defect group $\bar{P}=P/\langle z\rangle$. Assume that $\bar{b}'$ has weight $w$. Then $\bar{P}$ is isomorphic to a Sylow $2$-subgroup of $S_{2w}$. Since $\bar{P}$ is abelian, we have $w\leq 1$. It follows that $|\bar{P}|\leq2$, and hence $P$ is a cyclic group or a Klein four group. Thus by \cite[Theorem 1.12]{Kessar_Linckelmann} and \cite[Theorem 3]{H23}, Conjecture \ref{KLconj} holds for $\O\ws^\eta_nb$.

$\bullet$ Now assume that $b$ is a block of $\O \wa_n$ with an abelian defect group $P$. Let $b'$ be a block of $\O'\wa_n$ such that $bb'\neq 0$. As in the first case, the image $\bar{b}'$ of $b'$ under the canonical homomorphism $\O' \wa_n\to \O' A_n$ is a block of $\O' A_n$ with a defect group $\bar{P}=P/\langle z\rangle$. Since $\bar{P}$ is abelian, $\bar{P}$ is a trivial group or a Klein four group; see \cite[Lemma (7A)]{FH}. If $\bar{P}=1$, then $P=\langle z\rangle$. In this case Conjecture \ref{KLconj} holds for $\O \wa_nb$ by \cite[Theorem 1.12]{Kessar_Linckelmann}. Now assume that $\bar{P}$ is a Klein four group. Without loss of generality we may assume that $\bar{P}=\{(1),(12)(34),(13)(24),(14)(23)\}$. Let $x=t_1t_3$, a preimage of $(12)(34)$ in $P$. Since $t_3t_1=zt_1t_3$, $x^2=t_1t_3t_1t_3=t_1(zt_1t_3)t_3=z$. This implies that the two preimages of $(12)(34)$ in $P$ are of order $4$. By the same reason, the preimages of $(13)(24)$ or $(14)(23)$ in $P$ are of order $4$ as well. So $z$ is the unique involution in $P$, and this implies that $P$ is isomorphic to the quaternion group $Q_8$ of order $8$. This contradicts the assumption that $P$ is abelian.
\end{proof}

\section{Descent of twisted wreath superproducts}\label{s3:descent of twisted wreth superproducts}

\begin{void}\label{void:epsilon-twisted group superalgebra}
	{\rm Let $R\in \{\O,k\}$ and let $n$ be a positive integer. For $\eta\in\{\pm\}$, we define the {\it $\eta$-twisted group superalgebra} $\T^\eta_{n,R}$ of the symmetric group $S_n$ over $R$ to be the $R$-superalgebra given by odd generators $\xi_1,\cdots,\xi_{n-1}$ subject to the relations 
		\begin{align*}
			\xi_r^2=\eta(-1)^{\frac{(p-1)}{2}},\quad \xi_r\xi_s = -\xi_s\xi_r\text{ if }|r-s|>1,\quad (\xi_r\xi_{r+1})^3=1.
		\end{align*}
		For each $w\in S_n$, choosing a reduced expression $w=s_{r_1}\cdots s_{r_k}$ in terms of simple transpositions, we define $\xi_w:=\xi_{r_1}\cdots \xi_{r_k}$ (which depends up to a sign on the choice of a reduced expression).

	}
\end{void}

\begin{proposition}[Descent of $\T_{n,R}$]\label{prop:descent of TnR}
Let $(R,R_p)\in \{(\O,\Z_p),(k,\FF_p)\}$. If $\sqrt{-1}\in R$, then there is an isomorphism of $R$-superalgebras $\T^\eta_{n,R}\cong \T_{n,R}$ sending $\eta^i(-1)^{\frac{(p-1)i}{2}}\sqrt{\eta(-1)^{\frac{(p-1)}{2}}}\xi_i$ to $\tau_i$. Hence in this case, we have an isomorphism of $R$-superalgebras $ R\otimes_{R_p} \T^\eta_{n,R_p}\cong \T_{n,R}$ sending $\eta^i(-1)^{\frac{(p-1)i}{2}}\sqrt{\eta(-1)^{\frac{(p-1)}{2}}}\otimes\xi_i$ to $\tau_i$. In particular,  $\{\xi_w\mid w\in S_n\}$ is an $R_p$-basis of $\T_{n,R_p}^\eta$.
\end{proposition}

\begin{proof}
Since the sets $\{\eta^i(-1)^{\frac{(p-1)i}{2}}\sqrt{\eta(-1)^{\frac{(p-1)}{2}}}\xi_i\mid 1\leq i\leq n-1\}$ and $\{\tau_i\mid 1\leq i\leq n-1\}$ generate $\T^\eta_{n,R}$ and $\T_{n,R}$ respectively, and satisfy the same relations, we have the desired isomorphisms. Since $\{\tau_w\mid w\in S_n\}$ is an $R$-basis of $\T_{n,R}$, by the isomorphism $\T^\eta_{n,R}\cong \T_{n,R}$, $\{\xi_w\mid w\in S_n\}$ is an $R$-basis of $\T^\eta_{n,R}$. Hence $\{\xi_w\mid w\in S_n\}$ is also an $R_p$-basis of $\T_{n,R_p}^\eta$. 
\end{proof}

\begin{definition}
{\rm	Given an $R$-superalgebra $A$, we define the {\it $\eta$-twisted wreath superproduct} $A\wr_s \T^\eta_{n,R}$  as the free product $A^{\otimes n}\star \T^\eta_{n,R}$ of superalgebras subject to the relations
\begin{align*}%\label{algn:A_swr_Tn^R}
	\xi_r\,(a_1\otimes\dots\otimes a_n) = (-1)^{\sum_{u\neq r,r+1}|a_u|}\left({}^{s_r}(a_1\otimes\dots\otimes a_n)\,\xi_r\right)
\end{align*}
for all $a_1,\cdots,a_n\in A$ and $1\leq r \leq n-1$, where $s_r$ is the transposition $(r,r+1)\in S_n$.
}
\end{definition}

\begin{proposition}\label{proposition:epsilon-AwrT-graded}
Let $(K',\O',k')$ be a $p$-modular system which is an extension of $(K,\O,k)$ such that $\sqrt{-1}\in \O'$. Let $(R',R)\in \{(\O',\O),(k',k)\}$, $n$ a positive integer, and $A'$ an $R'$-superalgebra. Then the following hold:
\begin{enumerate}[{\rm (i)}]
	\item  The isomorphism $\T^\eta_{n,R'}\cong \T_{n,R'}$ in Proposition \ref{prop:descent of TnR} extends to an isomorphism of $R'$-superalgebras $A'\wr_s \T^\eta_{n,R'}\cong A'\wr_s\T_{n,R'}$.
	\item  If $A'\cong R'\otimes_R A$ for some $R$-superalgebra $A$, then $ R'\otimes_{R} (A\wr_s\T^\eta_{n,R})\cong A'\wr_s\T_{n,R'}$. 
	\item We have a decomposition \begin{align}\label{algn:epsilon-AwrT_graded}
		A\wr_s \T^\eta_{n,R} = \bigoplus_{w \in S_n} A^{\otimes n}\xi_w.
	\end{align}
 Hence $A\wr_s \T^\eta_{n,R}$ is an $S_n$-graded crossed $R$-superproduct. 
\end{enumerate}

\end{proposition}

\begin{proof}
 The $\eta^i(-1)^{\frac{(p-1)i}{2}}\sqrt{\eta(-1)^{\frac{(p-1)}{2}}}\xi_i$-conjugation action and the $\tau_i$-conjugation action on an element $(a_1\otimes\dots\otimes a_n)\in A'^{\otimes n}$ are the same, hence statement (i) holds. By statement (i) we have 
$$ R'\otimes_{R} (A\wr_s\T^\eta_{n,R})\cong (R'\otimes_{R}A)\wr_s (R'\otimes_{R} \T^\eta_{n,R})\cong A'\wr_s \T^\eta_{n,R'}\cong A'\wr_s \T_{n,R'}$$
as $R'$-superalgebras, whence (ii). By Proposition \ref{proposition:AwrT-graded}, we have
$$A'\wr_s \T_{n,R'} = \bigoplus_{w \in S_n} A'^{\otimes n}\tau_w.$$ 
Since $\eta^i(-1)^{\frac{(p-1)i}{2}}\sqrt{\eta(-1)^{\frac{(p-1)}{2}}}$ is invertible in $R'$, by statement (ii), we have
$$	A'\wr_s \T^\eta_{n,R'} = \bigoplus_{w \in S_n} A'^{\otimes n}\xi_w.$$
This equality restricts to the desired equality (\ref{algn:epsilon-AwrT_graded}), proving (iii).
\end{proof}

\begin{void}
	{\rm Keep the notation of Proposition \ref{proposition:epsilon-AwrT-graded}. Given the $S_n$-grading in (\ref{algn:epsilon-AwrT_graded}), recall the definition of $(A \wr_s \T^\eta_{n,R},A \wr_s \T^\eta_{n,R})_{S_n}$ from (\ref{equation:diag}).
		Let $M$ be an $(A,A)$-bisupermodule. Analogous to (\ref{algn:ext_M^n}), we denote by $M^{\boxtimes n}_{S_n}$ the $(A \wr_s \T^\eta_{n,R},A \wr_s \T^\eta_{n,R})_{S_n}$-supermodule given by extending the $(A^{\otimes n},A^{\otimes n})$-bisupermodule $M^{\boxtimes n}:=\underbrace{M\boxtimes_{R} \cdots\boxtimes_{R} M}_n$, via
		\begin{align}\label{algn:epsilon-ext_M^n}
			(\xi_r\otimes \xi_r^{-1})\cdot(m_1\otimes\cdots\otimes m_n)
			:=(-1)^{|m_r|+|m_{r+1}|}\left({}^{s_r}\,(m_1\otimes \cdots \otimes m_n)\right),
		\end{align}
		for all $m_1,\cdots,m_n\in M$ and $1\leq r\leq n-1$.
		
		Analogous to (\ref{algn:sup_wreath_bimod}), we now define
		\begin{align}\label{algn:epsilon-sup_wreath_bimod}
			M \wr_s \T^\eta_{n,R}:= {\rm Ind}_{(A \wr_s \T^\eta_{n,R}, A \wr_s \T^\eta_{n,R})_{S_n}}^{(A \wr_s \T^\eta_{n,R}) \otimes_{R} (A \wr_s \T^\eta_{n,R})^{\rm sop}}(M^{\boxtimes n}_{S_n}).
		\end{align}

	}
\end{void}

\begin{proposition}\label{proposition:epsilon-MwrT}
	Let $(K',\O',k')$ be a $p$-modular system which is an extension of $(K,\O,k)$ such that $\sqrt{-1}\in \O'$. Let $(R',R)\in \{(\O',\O),(k',k)\}$, $n$ a positive integer, $A'$ an $R'$-superalgebra and $M'$ an $(A',A')$-bisupermodule. Through the isomorphism $A'\wr_s \T^\eta_{n,R'}\cong A'\wr_s\T_{n,R'}$, we regard the $(A'\wr_s \T_{n,R'},A'\wr_s \T_{n,R'})$-bisupermodule $M'\wr_s \T_{n,R'}$ as an $(A'\wr_s \T^\eta_{n,R'},A'\wr_s \T^\eta_{n,R'})$-bisupermodule. Then the following hold:
	\begin{enumerate}[{\rm (i)}]
		\item  We have an isomorphism of $(A'\wr_s \T^\eta_{n,R'},A'\wr_s \T^\eta_{n,R'})$-bisupermodules $M'\wr_s \T^\eta_{n,R'}\cong M'\wr_s\T_{n,R'}$.
		\item  If $A'\cong R'\otimes_R A$ for some $R$-superalgebra $A$ and $M'\cong R'\otimes_R M$ for some $(A,A)$-bisupermodule $M$, then $ R'\otimes_{R} (M\wr_s\T^\eta_{n,R})\cong M'\wr_s\T_{n,R'}$ as $(A'\wr_s \T^\eta_{n,R'},A'\wr_s \T^\eta_{n,R'})$-bisupermodules. 
	\end{enumerate}

\end{proposition}

\begin{proof}
This is a routine excise. 
  \end{proof}               

\begin{void}
	{\rm  Let $R\in\{\O,k\}$, $n$ a positive integer and $A$ an $R$-superalgebra with a superunit $u$. If $C_2$ is a cyclic group of order $2$ with generator $g$, we can give $A$ the structure of a $C_2$-graded crossed $R$-superproduct  with graded components 
		$A_{g^{\varepsilon}}=A_{\0}u^{\varepsilon}=A_\varepsilon$
		for $\varepsilon\in\Z/2\Z$. 
		For any $n\in\Z_{>0}$, again we consider the wreath product 
		$$C_2\wr S_n=\{g_1^{\varepsilon_1}\cdots g_n^{\varepsilon_n}w\mid \varepsilon_1,\dots,\varepsilon_n\in\Z/2\Z,\, w\in S_n\}$$
		where $g_r$ is the generator of the $r$-th $C_2$ factor in the base group $C_2^{\times n}$. Then we can give $A\wr_s \T^\eta_{n,R}$ the structure of a crossed product of $A_{\0}^{\otimes n}$ and $C_2\wr S_n$ (see \ref{void:grad algebras}) with 
		graded components 
		$$
		(A\wr_s \T^\eta_{n,R})_{g_1^{\varepsilon_1}\cdots g_n^{\varepsilon_n}w}=
		A_{\0}^{\otimes n}u_1^{\varepsilon_1}\dots u_n^{\varepsilon_n} \xi_w
		$$
		where each $u_r$ is $1^{\otimes (r-1)}\otimes u\otimes 1^{\otimes (n-r)}$.
		
		Let $A$ now be a superalgebra with a superunit $u_A$ and 
		$B$ be a superalgebra with a superunit $u_B$. 
		Recalling (\ref{equation:diag}), we have the subsuperalgebras  
		\begin{align*}
			&(A,B)_{C_2}=(A_{\0}\otimes_R B_{\0}^{\rm sop})\oplus (A_{\bar{1}}\otimes_R B_{\bar{1}}^{\rm sop})\subseteq A\otimes_R B^{\rm sop},\\
			&(A\wr_s \T^\eta_{n,R},B\wr_s \T^\eta_{n,R})_{C_2\wr S_n}\subseteq (A\wr_s \T^\eta_{n,R})\otimes_R (B\wr_s \T^\eta_{n,R})^{\rm sop}.
		\end{align*}
		Setting $U_r:=(u_A)_r \otimes (u_B^{-1})_r$, for $1\leq r\leq n$, we have that $A_{\0}^{\otimes n}\otimes_R (B_{\0}^{\otimes n})^{\rm sop}$ together with the $U_r$'s generate a subsuperalgebra of $(A\wr_s \T^\eta_{n,R},B\wr_s \T^\eta_{n,R})_{C_2\wr_s S_n}$ isomorphic to $((A,B)_{C_2})^{\otimes n}$. We identify the two superalgebras. 
		Then $(A\wr_s \T^\eta_{n,R},B\wr_s \T^\eta_{n,R})_{C_2\wr_s S_n}$ is generated by $((A,B)_{C_2})^{\otimes n}$ and $\{\Xi_r:=\xi_r \otimes \xi_r^{-1}\mid 1\leq r\leq n-1\}$ subject to the following relations:
\begin{align}\label{algn:complex_rel}
	\begin{split}
		&\Xi_i (\underline{a} \otimes \underline{b}) = ({}^{s_i}\underline{a} \otimes {}^{s_i}\underline{b}) \Xi_i, \text{ for all }1\leq i\leq n-1,\\
		&\Xi_i^2 =\eta(-1)^{\frac{(p-1)}{2}}\otimes \eta(-1)^{\frac{(p-1)}{2}}=1, \text{ for all }1\leq i\leq n-1,\\
		&\Xi_i\Xi_j = \Xi_j\Xi_i, \text{ for all }1\leq i,j\leq n-1\text{ with } |i-j|>1,\\
		&\Xi_i\Xi_{i+1}\Xi_i = \Xi_{i+1}\Xi_i\Xi_{i+1}, \text{ for all }1\leq i\leq n-2,\\
		&\Xi_i U_j = U_j \Xi_i, \text{ for all }1\leq i \leq n-1, 1\leq j\leq n\text{ and }j\neq i,i+1,\\
		&\Xi_i U_i = U_{i+1} \Xi_i, \text{ for all }1\leq i \leq n-1,\\
		&\Xi_i U_{i+1} = U_i \Xi_i, \text{ for all }1\leq i \leq n-1,
	\end{split}
\end{align}
where $\underline{a} = a_1\otimes \dots \otimes a_n \in A_{\0}^{\otimes n}$ and $\underline{b} = b_1 \otimes \dots \otimes b_n \in (B_{\0}^{\otimes n})^{\rm op}$.

	}	
\end{void}

\begin{proposition}[a slight modification of {\cite[Proposition 5.1.14]{KL22}}]\label{prop:epsilon-ext_der_wreath}
	Let $R\in\{\O,k\}$, let $n$ be a positive integer, and let $A$ and $B$ be $R$-superalgebras with superunit. If $X$ is a complex of $A_{\0}$-$B_{\0}$-bimodules that extends to a complex of $(A,B)_{C_2}$-modules, then $X^{\otimes n}$ extends to a complex of $(A\wr_s \T^\eta_{n,R},B\wr_s \T^\eta_{n,R})_{C_2\wr S_n}$-modules.
\end{proposition}

\begin{proof}
Note that the set $\{\Xi_r\mid1\leq r\leq n-1\}$ in (\ref{algn:complex_rel}) and the set $\{T_r\mid1\leq r\leq n-1\}$ in \cite[(5.1.13)]{KL22} satisfy the same relations. So the proof of the present proposition is analogous to the proof of \cite[Proposition 5.1.14]{KL22} - we only need to replace ``$T_r$" there with $\Xi_r$ in (\ref{algn:complex_rel}). In \cite{KL22}, the ground field is assumed to be algebraically closed. But all the coefficients used in the proof of \cite[Proposition 5.1.14]{KL22} are in the prime field, so the ``algebraically closed" assumption is not necessary for \cite[Proposition 5.1.14]{KL22}. Hence the statement holds for $R=k$. The proof for the case $R=\O$ is identical.  
\end{proof}

\section{Descent and supermodules}\label{section:Descent and supermodules}

Let $(K',\O',k')$ be a $p$-modular system which is an extension of $(K,\O,k)$. 

\begin{proposition}\label{prop:extension}
Let $R\subseteq R'$ be two commutative rings. Let $A$ be an $R$-superalgebra, $M$ an $A$-supermodule. Then $A':=R'\otimes_R A$ is a superalgebra with $A'_{\0}:=R'\otimes_R A_{\0}$ and $A'_{\1}:=R'\otimes_R A_{\1}$, and $M':=R'\otimes_R M$ is an $A'$-supermodule with $M'_{\0}:=R'\otimes_R M_{\0}$ and $M'_{\1}:=R'\otimes_R M_{\1}$. A similar result also holds for bisupermodules instead of supermodules. %Moreover, similar results hold with $\O'$ and $\O$ replaced by $k'$ and $k$ respectively.
\end{proposition}

\begin{proof} 
This can be checked straightforward by the definitions of superalgebras, supermodules and bisupermodules. 
\end{proof}

\begin{void}\label{void:descents of supermodules}
{\rm Let $R\subseteq R'$ be two commutative rings. Let $A$ be an $R$-superalgebra and let $A':=R'\otimes_R A$. Let $U'$ be an $R'$-supermodule. We say that $U'$ is {\it defined over} $R$ or $U'$ {\it descends to} $R$, if there is an $A$-supermodule $U$ such that $U'\cong R'\otimes_R U$ as $A'$-supermodules. In this case, by the definition of an isomorphism of $R'$-supermodules, we have $U'_{\0}\cong R'\otimes_R U_{\0}$ and $U'_{\1}\cong R'\otimes_R U_{\1}$ as $A'_{\0}$-modules (recall from Proposition \ref{prop:extension} that $A'_{\0}=R'\otimes_R A_{\0}$). The converse holds if $A$ has a super unit:
}
\end{void}

\begin{proposition}\label{proposition: descent of zero-component}
Let $R\subseteq R'$ be two commutative rings. Let $A$ be an $R$-superalgebra with a superunit, and let $A':=R'\otimes_R A$. Let $U'$ be an $A'$-supermodule. If the $A'_{\0}$-module $U'_{\0}$ descends to $R$, then $U'$ descends to $R$.
\end{proposition}
\begin{proof}
By assumption, there is an $A_{\0}$-module $U_0$ such that $U'_{\0}\cong R'\otimes_R U_{\0}$ as $A'_{\0}$-module. Let $U:=A\otimes_{A_{\0}} U_{\0}$. It is easy to check that $U$ is an $A$-supermodule. By \cite[Lemma 3.7 (i)]{KL25}, we have $U'\cong  A'\otimes_{A'_{\0}} U'_{\0}$. Hence we have 
$$U'\cong  (R'\otimes_R A)\otimes_{A'_{\0}} (R'\otimes_{R} U_{\0})\cong R'\otimes_R (A\otimes_{A_{\0}} U_{\0})   =   R'\otimes_R U$$ as $A'$-supermodules.   
\end{proof}

\begin{remark}
{\rm In Proposition \ref{proposition: descent of zero-component}, if $A$ does not have a superunit, we don't know whether the statement is true or false, because in this case, $A_{\bar{1}}$ can not be expressed by $A_{\bar{0}}$. So we don't know the full information of $A_{\bar{1}}$ via $A_{\bar{0}}$.}
\end{remark}

\begin{proposition}[Descent of dual supermodules]\label{proposition:Descent of dual supermodules}
Let $R\subseteq R'$ be two commutative rings. Let $A$ be an $R$-superalgebra and let $A':=R'\otimes_R A$. Let $U$ be an $A$-supermodule which is projective as $R$-module. Then ${\rm Hom}_{R'}(R'\otimes_R U,R')\cong R'\otimes_R {\rm Hom}_R(U,R)$ as right $A'$-supermodules.
\end{proposition}

\begin{proof}
By \cite[Corollary 1.12.12]{Lin18a}, we have an $R'$-module isomorphism
$$R'\otimes_R {\rm Hom}_R(U,R)\cong {\rm Hom}_{R'}(R'\otimes_R U,R') $$
sending $1\otimes f$ to the map $1\otimes u\mapsto 1\otimes f(u)$, where $f\in {\rm Hom}_R(U,R)$ and $u\in U$.
One checks that this $R'$-module isomorphism is also an isomorphism of $A'$-supermodules.  
\end{proof}

\begin{proposition}[Descent of Morita superequivalences]\label{proposition:Descent of Morita superequivalences}
Let $(K',\O',k')$ be a $p$-modular system which is an extension of $(K,\O,k)$. Let $A$ and $B$ be supersymmetric $\O$-superalgebras, and let $A':=\O'\otimes_\O A$, $B':=\O'\otimes_\O B$. Let $M$ be an $(A,B)$-bisupermodule, and let $M':=\O'\otimes_\O M$. Then $M'$ and $M'^*$ induce a Morita superequivalence between $A'$ and $B'$ if and only if $M$ and $M^*$ induce a Morita superequivalence between $A$ and $B$. (Note that in this paper, all superspaces are defined to be finitely generated.) A similar statement holds with $\O'$ and $\O$ replaced by $k'$ and $k$ respectively.
\end{proposition}	

%\begin{proof}
%One direction of the implications is trivial. Now suppose that $M'$ and $M'^*$ induce a Morita superequivalence between $A'$ and $B'$. Now assume that we have an $(A',A')$-bisupermodule isomorphism
%\begin{equation}\label{equation:descent of morita superequ 1}
%$f':M'\otimes_{B'} M'^*\cong A'$
%\end{equation} 
%and a $(B',B')$-bisupermodule isomorphism $g': M'^*\otimes_{A'} M'\cong B'$. Then by \cite[Proposition 4.5]{Kessar_Linckelmann}, we have an $(A,A)$-bimodule isomorphism
%$f:M\otimes_{B} M^*\cong A$ 
%and a $(B,B)$-bimodule isomorphism $g: M^*\otimes_{A} M\cong B$. Moreover, by the proof of \cite[Proposition 4.5]{Kessar_Linckelmann}, we have $f'={\rm id}_{\O'}\otimes f$ as $(A',A')$-bimodule isomorphisms, and $g'={\rm id}_{\O'}\otimes g$ as $(B',B')$-bimodule isomorphisms. By Propositions \ref{prop:extension} and \ref{proposition:Descent of dual supermodules}, we have 
%$$A'_{\0}= \O'\otimes_\O A_{\0},~~A'_{\1}= \O'\otimes_\O A_{\1},$$ 
%$$(M'\otimes_{B'} M'^*)_{\0}= \O'\otimes_\O (M\otimes_{B} M'^*)_{\0}$$ and $$(M'\otimes_{B'} M'^*)_{\1}= \O'\otimes_\O (M\otimes_{B} M^*)_{\1}.$$ Since $f'$ is an $(A',A')$-bisupermodule isomorphism, we have $f'((M'\otimes_{B'} M'^*)_{\0})\subseteq A'_{\0}$ and $f'((M'\otimes_{B'} M'^*)_{\1})\subseteq A'_{\1}$.   This forces $f((M\otimes_{B} M^*)_{\0})\subseteq A_{\0}$ and $f((M\otimes_{B} M^*)_{\1})\subseteq A_{\1}$. Hence $f$ is an $(A,A)$-bisupermodule isomorphism. A similar argument shows that $g$ is a $(B,B)$-bisupermodule isomorphism. 
%\end{proof}

\begin{proof}
One direction of the implications is trivial. Conversely, suppose that $M'$ and $M'^*$ induce a Morita superequivalence between $A'$ and $B'$. Then $M'$ is projective as a left $A'$-module and as a right $B'$-module. By \cite[Lemma 4.4 (a)]{Kessar_Linckelmann}, $M$ is projective as a left $A$-module and as a right $B$-module. In particular, $M$ is projective as an $\O$-module. Note that $\Hom_A(M,A)$ and $\operatorname{Hom}_{B^{\mathrm{op}}}(M,B)$ are $(B,A)$-bisupermodules with $$\operatorname{Hom}_A(M,A)_\varepsilon=\left\{
	h\in\operatorname{Hom}_A(M,A)
	\;\middle|\;
	h(M_\delta)\subseteq A_{\delta+\varepsilon}
	\text{ for all }\delta\in\mathbb Z/2\mathbb Z
	\right\}$$
and $$\operatorname{Hom}_{B^{\mathrm{op}}}(M,B)_\varepsilon
=
\left\{
h\in \operatorname{Hom}_{B^{\mathrm{op}}}(M,B)\;\middle|\;
h(M_\delta)\subseteq B_{\delta+\varepsilon}
\text{ for all }\delta\in\mathbb Z/2\mathbb Z
\right\},$$	
where $\varepsilon\in \Z/2\Z$.
	
Fix supersymmetrizing forms $t_A:A\to\O$ and $t_B:B\to\O$. These forms induce isomorphisms of $(B,A)$-bisupermodules
	$$\begin{aligned}
		\theta_A:\operatorname{Hom}_A(M,A)&\xrightarrow{\cong}M^*,
		& h&\mapsto t_A\circ h,\\
		\theta_B:\operatorname{Hom}_{B^{\rm op}}(M,B)&\xrightarrow{\cong}M^*,
		& h&\mapsto t_B\circ h.
	\end{aligned}$$
Here, one checks that the maps $\theta_A$ and $\theta_B$ preserve the grading by using $t_A(A_{\bar1})=0$ and $t_B(B_{\bar1})=0$ (see \ref{void:supersymmetric algebras}). Let $A$-smod (resp. $B$-smod) denote the category of finitely generated $A$- (resp. $B$-) supermodules. Using the isomorphisms $\theta_A$ and $\theta_B$, the tensor functors
$$M\otimes_B-:B\text{-smod}\longrightarrow A\text{-smod}~~~{\rm and}~~~M^*\otimes_A-:A\text{-smod}\longrightarrow B\text{-smod}$$
are biadjoint; see \cite[Theorem 2.12.7]{Lin18a}. Denote by $$f:M\otimes_B M^*\longrightarrow A~~~{\rm and}~~~g:M^*\otimes_A M\longrightarrow B$$ the counits of these adjunctions. Explicitly, we have
	$$f(m\otimes\varphi)=\theta_A^{-1}(\varphi)(m)~~~{\rm and}~~~g(\varphi\otimes m)=\theta_B^{-1}(\varphi)(m),$$
for any $m\in M$ and $\varphi\in M^*$. One easily checks that $f$ and $g$ are bimodule homomorphisms. For $m\in M$ and $\varphi\in M^*$, the maps $\theta_A^{-1}(\varphi)$ and $\theta_B^{-1}(\varphi)$ have parity $|\varphi|$. Hence
$$f(m\otimes\varphi)\in A_{|m|+|\varphi|}~~~{\rm and}~~~g(\varphi\otimes m)\in B_{|\varphi|+|m|}.$$
Thus $f$ and $g$ are homomorphisms of bisupermodules.
	
Define supersymmetrizing forms $t_{A'}:A'\to \O'$ and $t_{B'}:B'\to \O'$ by setting
$$t_{A'}(r'\otimes a)=r't_A(a)~~~{\rm and}~~~t_{B'}(r'\otimes b)=r't_B(b)$$
for any $r'\in\O'$, $a\in A$ and $b\in B$.
Let$$f':M'\otimes_{B'}M'^*\longrightarrow A'~~~{\rm and}~~~g':M'^*\otimes_{A'}M'\longrightarrow B'$$ be the corresponding counits for the biadjoint tensor functors $M'\otimes_{B'}-$ and $M'^*\otimes_{A'}-$, constructed by using $t_{A'}$ and $t_{B'}$. The same argument shows that $f'$ and $g'$ are homomorphisms of bisupermodules.
	
By Proposition~\ref{proposition:Descent of dual supermodules}, we have a natural isomorphism $M'^*\cong\O'\otimes_\O M^*$ of $(B',A')$-bisupermodules. Under this isomorphism and the explicit definitions of $f$, $g$, $f'$ and $g'$, we have
$$f'=\operatorname{id}_{\O'}\otimes f~~~{\rm and}~~~g'=\operatorname{id}_{\O'}\otimes g.$$
Since $M'$ and $M'^*$ induce quasi-inverse equivalences, their counits $f'$ and $g'$ are isomorphisms. As $\O'$ is faithfully flat over $\O$ (see e.g. \cite[Lemma 4.1]{Kessar_Linckelmann}), it follows that $f$ and $g$ are isomorphisms. They are therefore isomorphisms of bisupermodules, proving that $M$ and $M^*$ induce a Morita superequivalence between $A$ and $B$. The proof over $k\subseteq k'$ is identical.
\end{proof}

\begin{remark}
{\rm Proposition \ref{proposition:Descent of Morita superequivalences} can easily be extended to the case of crossed products (see \ref{void:grad algebras}) instead of superalgebras, and to graded Morita equivalences (see \cite[Definition 3.3]{MarCA}) or graded Rickard equivalences (see \ref{void:grad algebras}) instead of Morita superequivalences. The proof is almost the same.
	
}
\end{remark}

\begin{void} 
{\rm \textbf{Projective supermodules and Heller translates.} Let $R\in\{\O,k\}$ and let $A$ be a superalgebra with a superunit. Following \cite[\S3.4]{KL25}, we call an $A$-supermodule $P$ a {\it projective supermodule} if $P$ is isomorphic to a direct summand of $A^{\oplus n}$ for some $n$. A {\it projective cover} of an $A$-supermodule is defined exactly like for modules. Our conditions on the superalgebra $A$ ensure that a projective cover of a finitely generated $A$-supermodule $M$ exists and is unique up to isomorphism. So we can define {\it Heller translate} $\Omega_A(M)$ of $M$ just like for modules. 
}
\end{void}

\begin{proposition}[Descent of Heller translates]\label{proposition:Descent of Heller translates}
Let $(K',\O',k')$ be a $p$-modular system which is an extension of $(K,\O,k)$. Let $A$ be a finite generated $\O$-algebra (resp. $\O$-superalgebra), $A':=\O'\otimes_\O A$, and $M'$ a finitely generated $A'$-module (resp. $A'$-supermodule). If $M'\cong \O'\otimes_\O M$ for some $A$-module (resp. $A$-supermodule) $M$, then $\Omega_{A'}(M')\cong \O'\otimes_\O \Omega_A(M)$ as $A'$-modules (resp. $A'$-supermodules). A similar statement holds with $\O'$ and $\O$ replaced by $k'$ and $k$ respectively.
\end{proposition}

\begin{proof}  We identify $M'$ and $\O'\otimes_\O M$. Let $\varphi:P\to M$ be a projective cover of $M$.   Then ${\rm id}_{\O'}\otimes \varphi:\O'\otimes_\O P\to M'$ is a projective cover of $M'$. Since $\O'$ is a flat $\O$-module (see \cite[Lemma 4.1]{Kessar_Linckelmann}), ${\rm ker}({\rm id}_{\O'}\otimes \varphi)\cong \O'\otimes_\O {\rm ker}(\varphi)$ as $\O'$-module, and hence also as $A'$-modules (resp. $A'$-supermodules).                \end{proof}

We refer to \cite[Theorem 5.2.1]{Lin18a} for the definition of Green correspondents. The following is a descent criterion for Green correspondents.

\begin{proposition}[Descent of Green correspondents]\label{prop:descent Green cor}
Let $(K',\O',k')$ be a $p$-modular system which is an extension of $(K,\O,k)$. Let $G$ be a finite group, $P$ a subgroup of $G$ and $H$ a subgroup of $G$ such that $N_G(P)\subseteq H$. Let $U'$ be an $\O'$-free indecomposable $\O'G$-module having $P$ as a vertex, and let $V'$ be an $\O'H$-module which is the Green correspondent of $U'$ with respect to $P$. Assume that $\O'$ is finitely generated as $\O$-module. Then $U'$ descends to $\O$ if and only if $V'$ descends to $\O$. A similar statement holds for $k'$ and $k$ instead of $\O'$ and $\O$ respectively.
\end{proposition}

\begin{proof}
Assume that $U'$ descends to $\O$, that is, there exists an $\O G$-module $U$ such that $U'\cong \O'\otimes_\O U$. By definition, $V'$ is up to isomorphism the unique direct summand of ${\rm Res}_H^G(U')$ with vertex $P$. By \cite[Lemma 5.1]{Kessar_Linckelmann}, we deduce that there is an indecomposable direct summand $V$ of ${\rm Res}_H^G(U)$ such that $V'\cong \O'\otimes_\O V$.

Conversely, assume that $V'$ descends to $\O$, that is, there exists an $\O H$-module $V$ such that $V'\cong \O'\otimes_\O V$. By definition, $U'$ is up to isomorphism the unique direct summand of ${\rm Ind}_H^G(V')$ with vertex $P$. By \cite[Lemma 5.1]{Kessar_Linckelmann}, again we deduce that there is an indecomposable direct summand $U$ of ${\rm Ind}_H^G(V)$ such that $U'\cong \O'\otimes_\O U$.   
\end{proof}

\begin{void}\label{void:super group algebra}
	{\rm \textbf{Super group algebras.} Throughout the rest of this section $G$ will denote a finite group with an index $2$ subgroup $G_{\0}$. Let $R$ be a commutative ring. The group algebra $RG$ has the structure of an $R$-superalgebra with superunit via
		$$(RG)_{\0}:=R G_{\0}~~~{\rm and}~~~(RG)_{\1}:=\Big\{\sum_{g\in G\backslash G_{\0}}\alpha_g g \mid \alpha_g\in R \Big\}.$$
	Moreover, if $b$ is a central idempotent of $RG$ which is contained in $RG_{\0}$, then $RGb$ has the structure of an $R$-superalgebra with superunit via
	$$(RGb)_{\0}:=R G_{\0}b~~~{\rm and}~~~(RGb)_{\1}:=(RG)_{\1}b.$$
		For example $G$ and $G_{\0}$ can be chosen to be $\ws^\eta_n$ and $\wa_n$ ($n\geq 2$) respectively. 
	For any subgroup $H$ of $G$ we will use $H_{\0}$ to denote $H\cap G_{\0}$ and $R H$ will inherit its superalgebra structure from that of $R G$.

Let $R\in \{\O, k\}$. Following \cite{KL25}, for a module-indecomposable $RG$-supermodule $M$, a {\it vertex} of $M$ means a vertex of $|M|$. The same applies to relative projectivity with respect to a subgroup and other standard notions of finite group representation theory.

	}
\end{void}

In \cite[Theorem 3.39]{KL25}, Kleshchev and Livesey introduced a super version of the usual Green correspondence:

\begin{theorem}[{\cite[Theorem 3.39]{KL25}}]\label{theo:defi-of-supergreen}
Let $P$ be a $p$-subgroup of $G$ and $H$ a subgroup of $G$ such that $N_G(P)\subseteq H$ and $H\nsubseteq G_{\0}$. Then there exists a correspondence between the module-indecomposable $RG$-supermodules with a vertex $P$ and the module-indecomposable $RH$-supermodules with a vertex $P$. Suppose that $U$ is a module-indecomposable $RG$-supermodule with a vertex $P$ corresponding to a module-indecomposable $RH$-supermodule $V$ with a vertex $P$. Then
$$\Ind_H^G(V)\cong U\oplus M_G~~~{\rm and}~~~\Res_H^G(U)\cong V\oplus M_H$$
for some $RG$-supermodule $M_G$ and some $RH$-supermodule $M_H$. Moreover, as an $R G$-module, $M_G$ is a direct sum of modules each with a vertex contained in $P\cap {}^gP$ for some $g\in G\backslash H$ and, as an $RH$-module, $M_H$ is a direct sum of modules each with a vertex contained in $H\cap {}^gP$ for some $g\in G\backslash H$. 
\end{theorem}

 The following is a super version of Proposition \ref{prop:descent Green cor}.

\begin{proposition}[Descent of super Green correspondents]\label{prop:descent super Green cor}
Let $(K',\O',k')$ be a $p$-modular system which is an extension of $(K,\O,k)$. Let $G$ be a finite group, $P$ a subgroup of $G$, and $H$ a subgroup of $G$ such that $N_G(P)\subseteq H$ and $H\nsubseteq G_{\0}$.  Let $U'$ be an $\O'$-free module-indecomposable $\O'G$-supermodule having $P$ as a vertex, and let $V'$ be an $\O'H$-module which is the super Green correspondent of $U'$. Assume that $\O'$ is finitely generated as $\O$-module. Then $U'$ descends to $\O$ if an only if $V'$ descends to $\O$. A similar statement holds for $k'$ and $k$ instead of $\O'$ and $\O$ respectively.
	\end{proposition}
	
\begin{proof} By the proof of \cite[Theorem 3.39]{KL25}, the $\O'H_{\0}$-module $V'_{\0}$ is the ordinary Green correspondent of the $\O'G_{\0}$-module $U'_{\0}$. By the descent criterion for the ordinary Green correspondents (Proposition \ref{prop:descent Green cor}), $U'_{\0}$ descends to $\O$ if and only if $V'_{\0}$ descends to $\O$. By \ref{void:descents of supermodules} and Proposition \ref{proposition: descent of zero-component}, $U'$ descends to $\O$ if and only if $U'_{\0}$ descends to $\O$, and $V'$ descends to $\O$ if and only if $V'_{\0}$ descends to $\O$. This completes the proof.    \end{proof}

\section{On Rickard equivalences and splendid complexes}\label{section:On Rickard equivalences and splendid complexes}

In this section $p$ is not necessarily being odd.
\begin{void}\label{void:grad algebras}
{\rm 	Let $R$ be a commutative ring and $G$ a finite group.  A finitely generated $R$-algebra $A$ is called {\it $G$-graded} (resp. {\it fully $G$-graded}) if, as an $R$-module, $A$ is a direct sum $A=\oplus_{g\in G} A_g$ satisfying $A_gA_h\subseteq A_{gh}$ (resp. $A_gA_h= A_{gh}$) for all $g,h\in G$; if $A^\times \cap A_g\neq \varnothing$, then $A$ is called a {\it crossed product} (of $A_1$ and $G$).  Let $A$, $B$ and $C$ be fully $G$-graded $R$-algebras. A finitely generated $A$-$B$-bimodule $M$ is called {\it $G$-graded} if, as an $R$-module, $M$ is a direct sum $M=\oplus_{x\in G} M_x$ satisfying $A_gM_xB_h=M_{gxh}$ for all $g,h,x\in G$. Hence with this terminology, if $G=\Z/2\Z$, then $A$ and $B$ are exactly $R$-superalgebras and $M$ is exactly an $(A,B)$-bisupermodule. Following \cite{MarCA}, we denote by $A$-Gr-$B$ the categroy of $G$-graded finitely generated $A$-$B$-bimodules and grade preserving $A$-$B$-bimodule homomorphisms. Let $K^b(\mbox{$A$-Gr-$B$})$ denote the homotopy category of bound complexes over $\mbox{$A$-Gr-$B$}$; see e.g. \cite[page 117]{Lin18a} for homotopy categories.
Let $M$ be a $G$-graded $A$-$B$-bimodule and $N$ a $G$-graded $B$-$C$-bimodule. Then $M\otimes_B N$ is a $G$-graded $A$-$C$-bimodule with $(M\otimes_B N)_x=\sum_{gh=x}M_g\otimes_B N_h$.
Hence for an object $X$ of $K^b(\mbox{$A$-Gr-$B$})$ and an object $Y$ of $K^b(\mbox{$B$-Gr-$C$})$, the complex $X\otimes_B Y$ is an object of $K^b(\mbox{$A$-Gr-$C$})$.

Let $A$ and $B$ be symmetric $R$-algebras. If $X$ is a bounded complex of $A$-$B$-bimodules that are finitely generated projective as left $A$-modules and as right $B$-modules, such that $X\otimes_B X^*$ and $A$ are homotopy equivalent as complexes of $A$-$A$-bimomdules, and that $X^*\otimes_A X$ and $B$ are homotopy equivalent as complexes of $B$-$B$-bimodules. Then we say that $A$ and $B$ are {\it Rickard equivalent}, and $X$ is called a {\it Rickard complex} of $A$-$B$-bimodules. Here $A$ (resp. $B$) denotes the complex with only one nonzero term $A$ (resp. $B$) in degree 0, and we adopt \cite[Definition 1.17.9]{Lin18a} for the definition of tensor product of complexes. In the literature, many authors investigate the Brou\'{e} conjecture in terms of the notion of derived equivalences, which looks weaker than Rickard equivalences. However, according to \cite[Theorem 2.21.3]{Lin18a}, for symmetric algebras, derived equivalences imply Rickard equivalences as well.

Let $A$, $B$ and $C$ be fully $G$-graded symmetric $R$-algebras. If there is a complex (resp. module) $X$ in $K^b(\mbox{$A$-Gr-$B$})$ (resp. $A$-Gr-$B$) where each term of $X$ is finitely generated projective as left $A$-module and as right $B$-module, such that $X\otimes_B X^*$ and $A$ are isomorphic in  $K^b(\mbox{$A$-Gr-$A$})$ (resp. $A$-Gr-$A$) and $X^*\otimes_A X$ and $B$ are isomorphic in $K^b(\mbox{$B$-Gr-$B$})$ (resp. (resp. $B$-Gr-$B$)). Then $A$ and $B$ are {\it graded Rickard (resp. Morita) equivalent}. Note that if $G=\Z/2\Z$, then a graded Morita equivalence is exactly a Morita superequivalence. 

}
\end{void}

Using \cite[Lemmas 3.6 and 3.19]{KL25}, it is straightforward to check the following:

\begin{proposition}\label{prop:boxtimes Rickard equivalences}
Let $R$ be a commutative ring, and let $A_1,A_2,B_1,B_2$ be superalgebras (in other words, $\Z/2\Z$-graded algebras). If 
$A_i$ and $B_i$ are $\Z/2\Z$-graded Rickard equivalent via the complexes $X_i$ of $(A_i,B_i)$-bisupermodules (in other words, $\Z/2\Z$-graded $A_i$-$B_i$-bimodules), for $i=1,2$, then $A_1\otimes_R A_2$ and $B_1\otimes_R B_2$ are $\Z/2\Z$-graded Rickard equivalent via the complex $X_1\boxtimes_R X_2$ of $(A_1\otimes_R A_2,B_1\otimes_R B_2)$-bisupermodules.	
\end{proposition}

Let $R$ be a commutative ring, $A$ and $B$ be fully $G$-graded $R$-algebras. An $(A,B)$-bimodule can be regarded as a module of the $R$-algebra $A\otimes_R B^{\rm op}$. We regard $B^{\rm op}$ as a $G$-graded algebra with components $B^{\rm op}_g=B_{g^{-1}}$, and $A\otimes_R B^{\rm op}$ as a $(G\times G)$-graded $R$-algebra with $(A\otimes_R B)_{(g,h)}=A_g\otimes_R B^{\rm op}_h$. Let $\Delta(A,B):=\oplus_{g\in G}A_g\otimes_R B^{\rm op}_g$, which is a fully $G$-graded $R$-algebra and which is also an $R$-subalgebra of $A\otimes_R B^{\rm op}$.

%\begin{lemma}[{\cite[Lemma 2.6]{MarCA}}]\label{lemma: Marcus CA Lemma 2.6}
%Let $R$ be a commutative ring, $G$ a finite group, $A$ and $B$ fully $G$-graded $R$-algebras.	There are mutually inverse functors induce a equivalence between the categories $\Delta(A,B)$-{\rm mod} and $A$-{\rm Gr}-$B$. In particular, the categories $K^b(\mbox{$\Delta(A,B)$-{\rm mod}})$ and $K^b(\mbox{$A$-{\rm Gr}-$B$})$ are equivalent.
%\end{lemma}

%An explicit pair of mutually inverse functors which induces the equivalence in Lemma \ref{lemma: Marcus CA Lemma 2.6} has been given in \cite[Lemma 2.6]{MarCA}. 

The following theorem is a homotopy category version of \cite[Theorem 4.7]{MarCA}, which is, as \cite[Theorem 4.7]{MarCA}, a straightforward consequence of \cite[4.2, 4.3]{MarCA}.

\begin{theorem}\label{theorem:Marcus CA Theorem 4.7}
Let $R$ be a commutative ring, $G$ a finite group, $A$ and $B$ fully $G$-graded symmetric $R$-algebras. The following are equivalent:
\begin{enumerate}[{\rm (i)}]
	\item There is an object $X\in K^b(\mbox{$A$-{\rm Gr}-$B$})$ inducing a graded Rickard equivalent between $A$ and $B$.
	\item There is an object $X_1\in K^b(\Delta(A,B))$, and isomorphisms $\alpha_1:X_1\otimes_{B_1} X_1^*\to A_1$ in $K^b(\Delta(A,A))$ and $\beta_1: X_1^*\otimes_{A_1} X_1\to B_1$ in $K^b(\Delta(B,B))$.
\end{enumerate}
Moreover, if the above conditions hold, then $X_1$ is the $1$-component of $X$,  and 
$$X\cong {\rm Ind}_{\Delta(A,B)}^{A\otimes_R B^{\rm op}}(X_1)\cong A\otimes_{A_1} X_1\cong X_1\otimes_{B_1} B.$$
\end{theorem}
 
If we only know that a Rickard complex for $A_1$ and $B_1$ extends to $\Delta(A,B)$, then we need $|G|$ to be invertible in $R$:

\begin{theorem}[A homotopy category version of {\cite[Theorem 4.8]{MarCA}}]\label{theorem:Marcus CA Theorem 4.8}
	Let $R$ be a commutative ring, $G$ a finite group, $A$ and $B$ fully $G$-graded symmetric $R$-algebras. Assume that $|G|$ is invertible in $R$. Then the following are equivalent:
	\begin{enumerate}[{\rm (i)}]
		\item There is an object $X\in K^b(\mbox{$A$-{\rm Gr}-$B$})$ inducing a graded Rickard equivalent between $A$ and $B$.
		\item There is an object $X_1\in K^b(\Delta(A,B))$ which, regarded as an object of $K^b(A_1\otimes_R B_1^{\rm op})$, induces a Rickard equivalence between $A_1$ and $B_1$.
	\end{enumerate}
	Moreover, if above conditions hold, then $X_1$ is the $1$-component of $X$,  and 
	$$X\cong {\rm Ind}_{\Delta(A,B)}^{A\otimes_R B^{\rm op}}(X_1)\cong A\otimes_{A_1} X_1\cong X_1\otimes_{B_1} B.$$ %are correspond to each other via the category equivalence in Lemma \ref{lemma: Marcus CA Lemma 2.6}.
\end{theorem}

\begin{proof} The proof is analogous to the proof of \cite[Theorem 4.8]{MarCA}.   \end{proof}

\begin{void}\label{void:Boltje splendid}
{\rm Let $G$ and $H$ be finite groups. An $(\O G,\O H)$-bimodule $M$ can be regarded as an $\O (G\times H)$-module (and vice versa) via $(g,h)m=gmh^{-1}$, where $g\in G$, $h\in H$ and $m\in M$. If $M$ is indecomposable as an $(\O G,\O H)$-bimodule, then $M$ is indecomposable as an $\O (G\times H)$-module, hence has a vertex (in $G\times H$) and a source.
	
Splendid equivalences was defined by Rickard in \cite{Rickard96}. Since Rickard's definition, it has been customary to use the adjective ``splendid" to any of the Rickard/derived or Morita equivalences.  In \cite[Remark 1.2]{Boltje}, Boltje suggested a more general version of splendid Rickard equivalences, which is more flexible and preserves all the good properties of the original version of Rickard \cite{Rickard96} and also later definitions (see for instance \cite[Definition 9.7.5]{Lin18b}). See \cite{BM} for the proofs of ``good properties" of this version of splendid Rickard equivalences. We adopt Boltje's wider splendidness. If $\varphi:P\xrightarrow{\cong} Q$ is an isomorphism between subgroups $P\leq G$ and $Q\leq H$, we set
$$\Delta\varphi:=\{(u,\varphi(u))\mid u\in P\},$$
and call it a {\it twisted diagonal subgroup} of $G\times H$. Let $R\in \{\O,k\}$. Whenever useful, we regard an $R\Delta \varphi$-module as an $R P$-module and vice versa via the isomorphism $P\cong \Delta \varphi$ sending $u\in P$ to $(u,\varphi(u))$. We say that an $RG$-$RH$-bimodule $M$ is {\it splendid} if every indecomposable direct summand of $M$, regarded as an $R(G\times H)$-module, is a $p$-permutation module with twisted diagonal vertex. We say that a complex $X=(X_i)_{i\in \Z}$ of $RG$-$RH$-bimodule is {\it splendid} if $X_i$ is splendid for any $i\in \Z$. Using the Mackey formula, one checks that a splendid bimodule is projective as left module and as right module. With this definition of ``splendid", splendid modules and splendid complexes are closed under induction and restriction:
}
\end{void}

\begin{proposition}\label{prop:splendid modules ind and res}
Let $R\in\{\O,k\}$. Let $G$ and $H$ be finite groups, $G_1$ a subgroup of $G$ and $H_1$ a subgroup of $H$. 
\begin{enumerate}[{\rm (i)}]
	\item If $M$ is a splendid (complex of) $RG_1$-$RH_1$-bimodule(s), then ${\rm Ind}_{G_1\times H_1}^{G\times H}(M)$ is a splendid (complex of) $RG$-$RH$-bimodule(s).
	\item  If $M$ is a splendid (complex of) $RG$-$RH$-bimodule(s), then ${\rm Res}_{G_1\times H_1}^{G\times H}(M)$ is a splendid (complex of) $RG_1$-$RH_1$-bimodule(s).
	\item Let $Y$ be a subgroup of $G\times H$ containing $G_1\times H_1$, and $M$ a/an (complex of) $RY$-module(s). Assume that $p\nmid|Y: G_1\times H_1|$ and ${\rm Res}_{G_1\times H_1}^Y(M)$ is splendid. Then ${\rm Ind}_Y^{G\times H}(M)$ is splendid.
\end{enumerate}
\end{proposition}

\begin{proof} It suffices to consider bimodules. If $M$ is a $RG_1$-$RH_1$-bimodule, then a vertex of ${\rm Ind}_{G_1\times H_1}^{G\times H}(M)$ is contained in a vertex of $M$, and this proves (i). One easily checks statement (ii) by using the Mackey formula. It remains to prove statement (iii). Since $p\nmid|Y: G_1\times H_1|$, $M$ is relatively $G_1\times H_1$-projective, hence $M$ is isomorphic to a direct summand of $${\rm Ind}_{G_1\times H_1}^Y{\rm Res}_{G_1\times H_1}^Y(M);$$ see e.g. \cite[Theorem 2.6.2]{Lin18a}. It follows that ${\rm Ind}_{Y}^{G\times H}(M)$ is isomorphic to a direct summand of 
$${\rm Ind}_{G_1\times H_1}^{G\times H}{\rm Res}_{G_1\times H_1}^Y(M),$$
 which is splendid by (i).  \end{proof}

  With the definition of ``splendid" in \ref{void:Boltje splendid}, tensor products and duals of splendid bimodules/complexes are splendid:

\begin{proposition}\label{prop:tensor products and duality}
	Let $R\in\{\O,k\}$. Let $G$, $H$ and $L$ finite groups. If $M$ is a splendid $RG$-$RH$-bimodule and $N$ is a splendid $RH$-$RL$-bimodule, then $M\otimes_{RH} N$ is a splendid $RG$-$RL$-bimodule, and $M^*:={\rm Hom}_R(M,R)$ is a splendid $RH$-$RG$-bimodule.
\end{proposition} 

\begin{proof} This first statement follows from \cite[Lemma 7.2 (b)]{BP}. Let $\varphi:P\xrightarrow{\cong} Q$ be an isomorphism between subgroups $P\leq G$ and $Q\leq H$. The second statement follows from the isomorphism $({\rm Ind}_{\Delta \varphi}^{G\times H}(R))^*\cong {\rm Ind}_{\Delta \varphi^{-1}}^{H\times G}(R)$. \end{proof}

 Propositions \ref{prop:splendid modules ind and res} (ii) and \ref{prop:tensor products and duality} explain why the definition of ``splendid" in \ref{void:Boltje splendid} is more flexible than the versions in \cite{Rickard96} and \cite[Definition 9.7.5]{Lin18b}.

\begin{definition}\label{defi:splendid over group algebras}
{\rm Let $R\in \{\O,k\}$. Let $G$ and $H$ be finite groups, $b\in Z(RG)$ and $c\in Z(RH)$ be nonzero idempotents. If there is a bounded splendid complex $X$ of finitely generated $RGb$-$RHc$-bimodules inducing a Rickard equivalence between $RGb$ and $RHc$. Then we say that $RGb$ and $RHc$ are {\it splendidly Rickard equivalent}, and $X$ is called a {\it splendid Rickard complex}.  }
\end{definition}

\begin{theorem}[a slight generalisation of {\cite[Theorem 5.2]{Rickard96}}]\label{theorem: lifting theorem of splendid Rickard complexes}
 Let $G$ and $H$ be finite groups, $b\in Z(\O G)$ and $c\in Z(\O H)$ be nonzero idempotents. Let $\bar{b}$ and $\bar{c}$ be the image of $b$ and $c$ in $kG$ and $kH$ respectively. If $\bar{X}$ is a splendid Rickard complex of $kG\bar{b}$-$kH\bar{c}$-bimodules. Then there is a splendid Rickard complex $X$ of $\O Gb$-$\O Hc$-bimodules such that $\bar{X}\cong k\otimes_\O X$.  Moreover, $X$ is unique up to isomorphism.
\end{theorem}

\begin{proof} With Rickard's original definition of a splendid complex (see \cite[page 336]{Rickard96}), this theorem has been proved in \cite[Theorem 5.2]{Rickard96}. Note that \cite[Theorem 5.2]{Rickard96} is stated under the assumption that the fraction field $K$ of $\O$ is a splitting field for $G$ and $H$, but the proof does not use this assumption. The only ingredient of the proof of \cite[Theorem 5.2]{Rickard96} is the completeness of $\O$ and the fact that $p$-permutation modules lift uniquely from $k$ to $\O$. Also, the lifting property of $p$-permutation modules does not make any requirement on size the field $k$; see \cite[Theorem 5.10.2 (iv)]{Lin18a}. Hence \cite[Theorem 5.2]{Rickard96} still holds for Definition \ref{defi:splendid over group algebras} instead of the splendid complex in \cite[page 336]{Rickard96}.   \end{proof}

\begin{void}\label{void:p-Permutation equivalences}
{\rm \textbf{$p$-Permutation equivalences.} Let $R\in \{\O,k\}$. Let $G$, $H$, $I$ be finite groups, and $b\in RG$, $c\in RH$, $d\in RI$ be central idempotents. 
	
(i) We denote by $T(RGb)$ (resp.  $T^\Delta(RGb,RHc)$) the free abelian group with standard basis given by the elements $[M]$, where $M$ runs over a set of representatives of the isomorphism classes of indecomposable $p$-permutation $RGb$-modules (resp. $p$-permutation $RGb$-$RHc$-bimodules with twisted diagonal vertices (as subgroups of $G\times H$)). Then taking $R$-duals of bimodules induces a $\Z$-linear isomorphism 
	$$-^*:T^\Delta(RGb,RHc)\cong T^\Delta(RHc,RGb),$$
and tensor products of bimodules induce a $\Z$-bilinear map
$$-\cdot_H-: T^\Delta(RGb,RHc)\times T^\Delta(RHc,RId)\to T^\Delta(RGb,RId).$$
According to Boltje and  Perepelitsky \cite[Definition 9.8]{BP}, a {\it $p$-permutation equivalence} between $RGb$ and $RHc$ is an element $\gamma\in T^\Delta(RGb,RHc)$ with the property that $\gamma\cdot_H \gamma^*=[RGb]$ and $\gamma^*\cdot_G \gamma=[RHc]$. If $X=(X_i)_{i\in \Z}$ is a splendid Rickard complex of $RGb$-$RHc$-bimodules, then by \cite[Theorem 1.6 (a)]{BP}, the element 
$\gamma:=\sum_{i\in \Z}(-1)^i[X_i]\in T^\Delta (RGb,RHc)$ is a $p$-permutation equivalence between $RGb$ and $RHc$. Note that the proof of \cite[Theorem 1.6]{BP} holds for arbitrary $p$-modular systems. For more properties of $p$-permutation equivalences, we refer the reader to \cite{BP}.

(ii) 	Let $(K',\O',k')$ be a $p$-modular system which is an extension of $(K,\O,k)$. Let $(R',R)\in \{(\O',\O),(k',k)\}$. Let $\gamma\in T^\Delta(RGb, RHc)$ be a $p$-permutation equivalence. Then $\gamma$ can be written as $[U]-[V]$, where $U$ and $V$ are direct sums of indecomposable $RGb$-$RHc$-bimodules with twisted diagonal vertices. Let $U':=R'\otimes_R U$, $V'=R'\otimes_R V$. Then by \cite[Lemma 5.1]{Kessar_Linckelmann}, $[U']-[V']$ is an element in $T^\Delta(R'Gb, R'Hc)$.  We define $R'\otimes_R \gamma\in T^\Delta(R'Gb, R'Hc)$ to be $[U']-[V']$, and this element is independent of the choice of the pair $(U,V)$. Set $\gamma'=R'\otimes_R \gamma$. Clearly we have $\gamma'\cdot_H \gamma'^*=[R'Gb]$ and $\gamma'^*\cdot_G \gamma'=[R'Hc]$, and hence $\gamma'$ is a $p$-permutation equivalence between $R'Gb$ and $R'Hc$.
}	
\end{void}

 The following theorem should be a useful tool for the (refined) Brou\'{e} conjecture. 

\begin{theorem}\label{theorem:splendid Rickard equ pres Brauer corr}
	Let $(K',\O',k')$ be a $p$-modular system which is an extension of $(K,\O,k)$. Assume that $\O'$ is finitely generated as $\O$-module. Let $(R',R)\in \{(\O',\O),(k',k)\}$. Let $G$ and $H$ be finite groups. Assume that $K'$ contains a primitive ${\rm exp}(G\times H)$-th root of unity. Let $b$ be a block of $R'G$ and $c$ a block of $R' H$. Assume that $b\in R G$ and $c\in R H$. Assume that $RGb$ and $RHc$ are $p$-permutation equivalent (in particular, splendidly Rickard equivalent). Then the Brauer correspondent of $RGb$ and the Brauer correspondent of $RHc$ are splendidly Morita equivalent.
\end{theorem}

To prove Theorem \ref{theorem:splendid Rickard equ pres Brauer corr}, we need the following lemma:

\begin{lemma}\label{lemma:scalar extension and idecomp summands}
Let $G$ be a finite group and $k'$ a finite extension of $k$. Let $U$ and $V$ be non-isomorphic finitely generated indecomposable $kG$-modules. Then any indecomposable direct summand of the $k'G$-module $k'\otimes_k V$ is not isomorphic a direct summand of $k'\otimes_k U$.

\end{lemma}

\begin{proof}
Suppose that there exists an indecomposable $k'G$-module $X$ which is isomorphic to a direct summand of both $k'\otimes_k U$ and $k'\otimes_k V$. Let $d=[k':k]$. Then we have
$$\Res^{k'G}_{kG}(k'\otimes_k U)\cong U^{\oplus d}~~~{\rm and}~~~\Res^{k'G}_{kG}(k'\otimes_k V)\cong V^{\oplus d}.$$
Hence $\Res^{k'G}_{kG}(X)$ is a direct summand of both $U^{\oplus d}$ and $V^{\oplus d}$. Since $U$ is indecomposable, the Krull--Schmidt theorem implies that every indecomposable direct summand of $\Res^{k'G}_{kG}(X)$ is isomorphic to $U$. Thus
$\Res^{k'G}_{kG}(X)\cong U^{\oplus a}$ for some integer $a\geq 1$. Similarly, $\Res^{k'G}_{kG}(X)\cong V^{\oplus b}$
for some integer $b\geq 1$. Therefore
$U^{\oplus a}\cong V^{\oplus b}$.
Applying the Krull--Schmidt theorem once more yields $U\cong V$, contradicting the assumption that $U$ and $V$ are non-isomorphic.
\end{proof}

\begin{proof}[Proof of Theorem \ref{theorem:splendid Rickard equ pres Brauer corr}]  By lifting properties of $p$-permutation modules (see \cite[Theorem 5.10.2 (iv)]{Lin18a}), we may assume that $(R',R)=(k',k)$. Let $\gamma\in T^\Delta(kGb, kHc)$ be a $p$-permutation equivalence. Then $\gamma':=k'\otimes_k \gamma$ is a $p$-permutation equivalence; see \ref{void:p-Permutation equivalences} (ii). By \cite[Theorem 14.1]{BP}, there is a unique (up to isomorphism) indecomposable $k'Gb$-$k'Hc$-bimodule $M'$ appearing in $\gamma'$, who has a twisted diagonal vertex $\Delta\varphi$, where $\varphi:P\xrightarrow{\cong} Q$ is an isomorphism between subgroups $P\leq G$ and $Q\leq H$, such that $P$ is a defect group of $k'Gb$ and $Q$ is a defect group of $k'Hc$. Again by \cite[Theorem 14.1]{BP}, any other indecomposable $k'Gb$-$k'Hc$-bimodule appearing in $\gamma'$ has vertex with order strictly smaller than $|P|$. Then by Lemma \ref{lemma:scalar extension and idecomp summands} and \cite[Lemma 5.1]{Kessar_Linckelmann}, we see that there is a unique (up to isomorphism) indecomposable $kGb$-$kHc$-bimodule $M$ appearing in $\gamma$, who has a twisted diagonal vertex $\Delta\varphi$, such that $M'\cong k'\otimes_k M$. Let $b'$ be the Brauer correspondent of $b$ in $k'N_G(P)$ and $c'$ the Brauer correspondent of $c$ in $k'N_H(Q)$. Since $b\in kG$ and $c\in kH$, we have $b'\in kN_G(P)$ and $c'\in kN_H(Q)$. Let $L'$ be a $k'N_{G\times H}(\Delta \varphi)$-module which is the Green correspondent of $M'$. Then by \cite[Theorem 14.5]{BP}, the $p$-permutation $k'N_G(P)b'$-$k'N_H(Q)c'$-bimodule 
$$\Ind_{N_{G\times H}(\Delta\varphi)}^{N_G(P)\times N_H(Q)}(L')$$ induces a Morita equivalence between $k'N_G(P)b'$ and $k'N_H(Q)c'$. By Proposition \ref{prop:descent Green cor}, there is a $kN_{G\times H}(\Delta \varphi)$-module $L$ such that $L'\cong k'\otimes_k L$. By \cite[Lemma 5.2]{Kessar_Linckelmann}, $L$ is a $p$-permutation $kN_{G\times H}(\Delta \varphi)$-module. Now by \cite[Proposition 4.5 (c)]{Kessar_Linckelmann}, the $p$-permutation $kN_G(P)b'$-$kN_H(Q)c'$-bimodule $${\rm Ind}_{N_{G\times H}(\Delta \varphi)}^{N_G(P)\times N_H(Q)}(L)$$
induces a Morita equivalence between $kN_G(P)b'$ and $kN_H(Q)c'$.
\end{proof}

\begin{void}
{\rm Let $R\in\{\O,k\}$, and let $G$ and $H$ be finite groups. By definition, an $RG$-$RH$ module $M$ is splendid if it is a direct sum of summands to of the $R(G\times H)$-bimodules ${\rm Ind}_{\Delta \varphi}^{G\times H}(R)$, with $\Delta \varphi$ running over the twisted diagonal subgroups of $G\times H$. For any $\varphi:P\cong Q$, note that we have an isomorphism of $RG$-$RH$-bimodules
	$${\rm Ind}_{\Delta \varphi}^{G\times H}(R)\cong RG\otimes_{R P} {}_\varphi RH\cong RG_{\varphi^{-1}}\otimes_{RQ} RH $$
	sending $(x,y)\otimes 1$ to  $x\otimes y^{-1}$ for all $x\in G$ and $y\in H$. Here ${}_\varphi RH$ is the $RP$-$RH$-bimodule which is equal to $RH$ as a right $RH$-module endowed with the left structural homomorphism $P\xrightarrow{\varphi} Q\to {\rm End}_R(RH)$; the notation $RG_{\varphi^{-1}}$ is similarly defined. 
	
	More generally, as in \cite[Definition 1.3]{HL}, we can define splendid modules, complexes, and Rickard equivalences for symmetric interior $P$-algebras:
}	

\end{void}

\begin{definition}\label{defi: splendid Rickard equivalences for interior P algebras}
{\rm Let $R\in\{\O,k\}$ and let $P$ be a finite $p$-group. Let $A$ and $B$ be symmetric interior $P$-algebras over $R$. An $A$-$B$-bimodule $M$ is called {\it splendid} if it is a direct sum of summands of the $A$-$B$-bimodule $A\otimes_{RQ} {}_\varphi B\cong A_{\varphi^{-1}}\otimes_{RT} B$, with $Q$, $T$ running over subgroups of $P$ and $\varphi:Q\to T$ running over the set of isomorphisms between $Q$ and $T$.  We say that a complex $X=(X_i)_{i\in \Z}$ of $A$-$B$-bimodule is {\it splendid} if $X_i$ is splendid for any $i\in \Z$. If there is a bounded splendid complex $X$ of finitely generated $A$-$B$-bimodules inducing a Rickard equivalence between $A$ and $B$. Then we say that $A$ and $B$ are {\it splendidly Rickard equivalent}, and $X$ is called a {\it splendid Rickard complex}. 
} 
\end{definition} 

If one specializes Definition \ref{defi: splendid Rickard equivalences for interior P algebras} to the case of blocks of group algebras, the previous definitions are recovered.

\begin{proposition}\label{prop:boxtimes splendid}
	Let $R\in\{\O,k\}$. Let $P$ be a finite $p$-group, and let $A_1,A_2,B_1,B_2$ be interior $P$-algebras over $R$, which are also $R$-superalgebras with the images of $P$ contained in their even parts $(A_1)_{\0},(A_2)_{\0},(B_1)_{\0},(B_2)_{\0}$. If $M_i$ is a splendid (complex of) $A_i$-$B_i$-bimodule(s) which is also a/an (complex of) $A_i$-$B_i$-bisupermodule, for $i=1,2$, then $M_1\boxtimes_R M_2$ is a splendid (complex of) $(A_1\otimes_R A_2, B_1\otimes_R B_2)$-bimodule(s).	
\end{proposition}

\begin{proof} It suffices to consider bimodules. Let $Q_1,Q_2,T_1,T_2$ be subgroups of $P$ and let $\varphi_1:Q_1\cong T_1$, $\varphi_2:Q_2\to T_2$ be group isomorphisms. It suffices to show that $$(A_1\otimes_{RQ_1} {}_{\varphi_1} B_1)\boxtimes_R (A_2\otimes_{RQ_2} {}_{\varphi_2} B_2)$$
is a splendid $(A_1\otimes_R A_2, B_1\otimes_R B_2)$-bimodule. By \cite[Lemma 3.6]{KL25}, we have an isomorphism of $(A_1\otimes_R A_2, B_1\otimes_R B_2)$-bimodules
$$(A_1\otimes_{RQ_1} {}_{\varphi_1} B_1)\boxtimes_R (A_2\otimes_{RQ_2} {}_{\varphi_2} B_2)\cong (A_1\otimes_R A_2)\otimes_{R(Q_1\times Q_2)} {}_\varphi(B_1\otimes_R B_2)$$
sending $(a_1\otimes b_1)\otimes (a_2\otimes b_2)$ to $(-1)^{|b_1||a_2|}(a_1\otimes a_2)\otimes(b_1\otimes b_2)$ for all $a_1\in A_1$, $a_2\in A_2$, $b_1\in B_1$ and $b_2\in B_2$, where $\varphi$ is the group isomorphism 
$\varphi_1\times \varphi_2:Q_1\times Q_2\cong T_1\times T_2$, $(u_1,u_2)\mapsto (\varphi_1(u_1),\varphi_2(u_2))$. This completes the proof. \end{proof}

%By \cite[Theorem 1.1]{BP}, the isomorphism $\varphi:P\cong Q$ is an isomorphism between the fusion systems of $k'Gb$ and $k'Hc$ determined by the maximal Brauer pairs $(P,e)$ and $(Q,f)$, respectively.   $\hfill\square$

\section{On spin blocks of $\ws^\eta_n$ and $\wa_n$}\label{section:On spin blocks of ws_n and wa_n}

\begin{notation}
{\rm Let $n$ be a positive integer. We use $[n]$ to denote the set $\{1,\cdots,n\}$. For conventional convenience we set $[0]:=\varnothing$. We denote by $\mathscr{P}_0(n)$ the set of strict partitions of $n$, i.e. partitions of $n$ without repeat parts. Set $\mathscr{P}_0(0):=\{{\rm empty~partion}\}$.
For $\lambda\in\mathscr{P}_0(n)$, ${\rm Ab}_\lambda$ will denote its $\bar{p}$-abacus ($p$-bar-abacus) display; see \cite[\S2.3.b]{KL22} for more details on this. For any non-empty partition $\lambda$, set $h(\lambda):={\rm max}\{k~|~\lambda_k>0\}$ to be the {\it length} of $\lambda$. We adopt the convention that the $0$-th position on the abacus display is always unoccupied, i.e. we always use $h(\lambda)$ beads in ${\rm Ab}_\lambda$.
	 Let $\mathscr{P}_0:=\bigsqcup_{n\in \mathbb{N}} \mathscr{P}_0(n)$.  Following \cite{Morris}, we can associate to every $\lambda\in \mathscr{P}_0$ its $\bar{p}$-core ${\rm core}(\lambda)\in \mathscr{P}_0$. We can also define a non-negative integer called the {\it$\bar{p}$-weight} of $\lambda$: ${\rm wt}(\lambda):=(|\lambda|-|{\rm core}(\lambda)|)/p\in \mathbb{N}$.

Assume in this section that $(K,\O,k)$ is large enough for $\ws^\eta_n$; this means that $K$ (and hence $\O$) contains a primitive $(2n!)$-th root of unity. Following \cite{KL22,KL25} we also assume that $-1$ and $2$ have square roots in $K$ (and hence $\O$); this condition is automatic if $n\geq 4$. Let $(R,R_p)\in \{(\O,\Z_p),(k,\FF_p)\}$. If $n>1$, we consider $R\ws^\eta_n$ as a superalgebra corresponding to $(\ws^\eta_n)_{\0}=\wa_n$; see \ref{void:super group algebra}.  We set 
$$e_z:=(1-z)/2\in R\ws^\eta_n.$$
We will use the adjective ``spin" when referring to objects  associated to $R\ws^\eta_n e_z$ rather than the whole of $R\ws^\eta_n$. For example the spin characters of $\ws^\eta_n$ will refer to the 
$K$-characters $\chi$ of $\ws^\eta_n$ satisfying $\chi(z)=-\chi(1)$, while spin blocks refers to the blocks of $R\ws^\eta_ne_z$. Note that $R\ws^\eta_ne_z$ inherits the structure of superalgebra from $R\ws^\eta_n$ since the idempotent $e_z$ is in $(R\ws^\eta_n)_{\0}$.

Let $G \leq \ws^\eta_n$ with $z \in G$. Following the notation from \ref{void:super group algebra}, we set 
\begin{equation}\label{equation:super-group-algebra-in-doublecover}
	G_{\0} := G \cap \wa_n\quad \text{and} \quad G_{\1} := G \setminus \wa_n
\end{equation}
and treat $R Ge_z$ as a superalgebra in the appropriate way. (Note that unlike in \ref{void:super group algebra} we are yet not assuming that $G_{\0}$ is a proper subgroup of $G$.)

As stated in \cite[equation (4.3)]{KL25}, we have the following isomorphisms of $R$-superalgebras (see \ref{void:twisted group superalgebra} for the notation $\T_{n,R}$):
\begin{equation}\label{equation:Sn_Tn_isom}
	\begin{split}
		R \ws_n^+ e_z & \cong \T_{n,R},~t_i e_z  \mapsto \tau_i,\ z e_z  \mapsto -1
		\\
		R \ws_n^- e_z & \cong \T_{n,R},~t_i e_z  \mapsto (-1)^i\sqrt{-1}\tau_i,\  z e_z  \mapsto -1.
	\end{split}
\end{equation}
}
\end{notation}

\begin{void}
	{\rm \textbf{Special subgroups.} Consider the subgroups $G \leq \ws^\eta_m$ (where $m$ is a non-negative integer) and $H \leq \ws^\eta_n$, both containing the canonical central element $z$ (note that we do not distinguish between the $z$ in $G$ and that in $H$ as these elements will be identified anyway). If $g\in \ws^\eta_m\backslash \wa_m$ we write $|g|=1$ and if $g\in \wa_m$ we write $|g|=0$. Similarly for $h\in \ws^\eta_n$. We set $G \times_z H$ to be the subgroup of $\ws^\eta_{m+n}$ generated by $G$ and $H$, where we view $H \leq \ws^\eta_n$ as a subgroup of $\ws^\eta_{m+n}$ via $z \mapsto z$, $t_i \mapsto t_{m+i}$.
		
	As in \cite{KL25}, if $X\subseteq [n]$, then $S_X$ will signify the subgroup of $S_n$ consisting of all elements that fix $[n]\backslash X$ pointwise and $\ws^\eta_X$ will denote $\pi_n^{-1}(S_X) \leq \ws^\eta_n$. Moreover, for disjoint $X_1, \cdots, X_i \subseteq [n]$, let $\ws^\eta_{X_1,\cdots,X_i}$ be the subgroup of $\ws^\eta_n$ generated by $\ws^\eta_{X_1},\cdots,\ws^\eta_{X_i}$. In particular,
	\begin{align*}
		\ws^\eta_{X_1,\cdots,X_i} \cong \ws^\eta_{X_1} \times_z \cdots \times_z \ws^\eta_{X_i}
	\end{align*}
	and
	\begin{align}\label{EZCommAlg}
		R\ws^\eta_{X_1,\cdots,X_i}e_z \cong R\ws^\eta_{X_1}e_z \otimes_R \cdots \otimes_R R\ws^\eta_{X_i}e_z
	\end{align}
	as superalgebras. We also define the corresponding subgroups of $\wa_n$:
	$$\wa_X:=\ws^\eta_X \cap \wa_n\quad\text{and}\quad \wa_{X_1,\cdots,X_i}:=\ws^\eta_{X_1,\cdots,X_i} \cap \wa_n.
	$$
	In the important special case where $n=n_1+\cdots+n_i$ and 
	$Y_l=[n_1+\cdots+n_l]\setminus[n_1+\cdots+n_{l-1}]$ for $l=1,\cdots,r$, we have the {\it standard Young subgroups}
	\begin{equation*}\label{EYoungSubgr}
		\ws^\eta_{n_1,\cdots,n_i}:=\ws^\eta_{Y_1,\cdots,Y_i} \cong \ws^\eta_{n_1} \times_z \cdots \times_z \ws^\eta_{n_i}\leq \ws^\eta_n\quad\text{and}\quad
		\wa_{n_1,\cdots,n_i}=\ws^\eta_{n_1,\cdots,n_i}\cap\wa_{n}\leq \wa_n.
	\end{equation*}
	If some of the $n_l$'s are equal to each other, we use the exponential notation, as for example for the standard Young subgroup $\ws^\eta_{r+ip,p^{d-i}}=\ws^\eta_{r+ip,p,\cdots,p}$. 
	}	
\end{void}

Following Schur \cite{Sch}, to every strict $\lambda\in \mathscr{P}_0(n)$, one associated canonically a spin character $\chi_\lambda$ of $\ws^\eta_n$ such that $\chi_\lambda$ is irreducible  if $n-h(\lambda)$ is even and $\chi_\lambda=\chi_\lambda^++\chi_\lambda^-$ for irreducible characters $\chi_\lambda^\pm$ if $n-h(\lambda)$ is odd. Moreover,
$$\{\chi_\lambda\mid \lambda\in\mathscr{P}_0(n),~n-h(\lambda)~{\rm is~even}\}\cup \{\chi_\lambda^\pm\mid \lambda\in\mathscr{P}_0(n),~n-h(\lambda)~{\rm is~odd}\} $$
is a complete irredundant set of irreducible spin characters of $\ws^\eta_n$.

\begin{theorem}[{\cite{Hum}}, see also {\cite[Theorem 5.2.4]{KL22}}]
	With one exception, $\chi_\lambda^\pm$ and $\chi_\mu^\pm$ are in the same $p$-block if and only if ${\rm core}(\lambda)={\rm core}(\mu)$. The exception is $\lambda=\mu$ is a $\bar{p}$-core and $n-h(\lambda)$ is odd, in which case $\chi_\lambda^+$ and $\chi_\lambda^-$ are in different $p$-blocks.
\end{theorem}

The blocks of $R\ws^\eta_n$ are described by the following theorem: 

\begin{theorem}[{\cite[Theorem 1.1]{Hum}}, see also {\cite[Theorem 4.3]{KL25}}]\label{theo: blocks of double covers of Sn}
	\begin{enumerate}[{\rm (i)}]
		\item The spin blocks of $R\ws^\eta_n$ are labelled by pairs $(\rho,d)$, where $\rho$ is a $\bar{p}$-core and $d$ is a non-negative integer with $n=|\rho|+dp$. If $d=0$ and $n-h(\rho)$ is odd then $(\rho,0)$ labels two spin blocks of $R\ws^\eta_n$ with corresponding block idempotents $e_{\rho,0}^+,e_{\rho,0}^- \in R \ws^\eta_n e_z$. In all cases we denote the block idempotent (or sum of two block idempotents) corresponding to $(\rho,d)$ by $e_{\rho,d}\in R \ws^\eta_n e_z$.
		\item The character $\chi_\lambda^\pm$ lies in the (sum of) block(s) $\O \ws^\eta_n e_{\rho,d}$, where $\rho$ is the $\bar{p}$-core of $\lambda$. If $d=0$ and $\rho$ is odd, we choose the labeling so that $\chi_\rho^{\pm}$ lies in the block $\O \ws^\eta_n e_{\rho,0}^{\pm}$.
	\end{enumerate}
\end{theorem}

Defect groups and Brauer correspondents of blocks of $R \ws^\eta_n$ are described as follows:

\begin{theorem}[{\cite[Theorem A and Corollary 26]{Ca}}, see also {\cite[Theorem 4.3]{KL25}}]\label{theo:defect group and Brauer correspondent}
Keep the notation of Theorem \ref{theo: blocks of double covers of Sn}. Fix a $\bar{p}$-core $\rho$ and $d$ a non-negative integer with $n=|\rho|+dp$, and 
set $U:=[|\rho|]$ and $V:=[n]\setminus U$.
\begin{enumerate}
	\item[{\rm(i)}] If $d=0$, then each block of $R\ws^\eta_n e_{\rho,d}$ has trivial defect group. If $d>0$, then any Sylow $p$-subgroup $D$ of $\ws^\eta_{V} \leq \ws^\eta_n$ is a defect group of $R
	\ws^\eta_n e_{\rho,d}$.
	\item[{\rm(ii)}] If $d=0$, then each block of $R\ws^\eta_n e_{\rho,d}$ is its own Brauer correspondent. If $d>0$, then, setting $D$ to be as in (i), the Brauer correspondent of $R
	\ws^\eta_n e_{\rho,d}$ in $N_{\ws^\eta_n}(D)$ is $R N_{\ws^\eta_n}(D) e_{\rho,0}$. Here we view $e_{\rho,0}$ as an element of $R N_{\ws^\eta_n}(D)$ via $e_{\rho,0}\in R\ws^\eta_U \hookrightarrow R N_{\ws^\eta_n}(D)$. If $\rho = \varnothing$, we interpret $e_{\varnothing,0}$ as $e_z$.
\end{enumerate}
\end{theorem}

\begin{remark}[{\cite[Remark 4.4]{KL25}}]\label{remark:block idempotents}
{\rm The $R \ws^\eta_n e_{\rho,d}$'s are precisely the superblocks of $R \ws^\eta_n e_z$. That is, the finest decomposition of $R\ws^\eta_n e_z$ into a direct sum of two-sided ideals that are invariant under the involution $\sigma_{R \ws^\eta_n}$ (see (\ref{equation:involution map}) for this notation). This follows from the distribution of characters given in Theorem \ref{theo: blocks of double covers of Sn} (ii). In particular, we always have $e_{\rho,d} \in R \wa_n e_z$.}
\end{remark}

The spin blocks of $R\wa_n$ are described in \cite[Proposition 3.16]{K96}:

\begin{theorem}\label{thm:An_blocks}
For each $n>1$, a $\bar{p}$-core $\rho$ and a non-negative integer $d$ with $n=|\rho|+dp$, $R \wa_n e_{\rho,d}$ is a single block of $R \wa_n$ unless $d=0$ and $\rho$ is even. In this latter case $R\wa_n e_{\rho,d}$ is a direct sum of two blocks of $R \wa_n$. If $n=1$, then $R \wa_1 e_{(1),0} = R \wa_1 e_{(1),0} = R e_z$ is a single block of $R \wa_1$.
	
	If $d=0$, once again, the defect group is trivial. In all other cases the defect group of $R \wa_n e_{\rho,d}$ is the same as that of $R \ws^\eta_n e_{\rho,d}$.
\end{theorem}

We set
$$B^{\eta,\rho,d}_R:=R\ws^\eta_ne_{\rho,d}~~{\rm and}~~(B^{\eta,\rho,d}_{R})_{\0}:=R\wa_ne_{\rho,d}.$$
It is easy to see that $R \wa_n e_{\rho,d}$ is indeed the purely even subalgebra of the superalgebra $R\ws^\eta_ne_{\rho,d}$. We will often assume $d>0$ to ensure that blocks and superblocks coincide. However, the corresponding statements for trivial defect blocks are completely elementary. 

\begin{proposition}\label{prop:descent of block idempotents}
Keep the notation of Theorem \ref{theo: blocks of double covers of Sn}. The central idempotent $e_{\rho,d}$ of $R \ws^\eta_n$ is contained in the subring $R_p\ws^\eta_n$, hence also contained in $R_p\wa_n$ by Remark \ref{remark:block idempotents}. 
	\end{proposition}

\begin{proof}
It suffices to consider $R=\O$, since the case $R=k$ follows by reduction modulo $J(\O)$. Set
$$\Lambda_{\rho,d}:=\{\lambda\in\mathscr{P}_0(n)\mid {\rm core}(\lambda)=\rho\}.$$
For each $\lambda\in\Lambda_{\rho,d}$, let $a_\lambda$ be the common degree of the irreducible constituents of $\chi_\lambda$. Thus $a_\lambda=\chi_\lambda(1)$ if $\chi_\lambda$ is irreducible, and $a_\lambda=\chi_\lambda^+(1)=\chi_\lambda^-(1)$ otherwise. By Theorem~\ref{theo: blocks of double covers of Sn},
$$\Theta_{\rho,d}:=\sum_{\lambda\in\Lambda_{\rho,d}}a_\lambda\chi_\lambda =\sum_{\psi\in{\rm Irr}(K\ws^\eta_ne_{\rho,d})}\psi(1)\psi.$$
The formula for primitive central idempotents therefore gives
$$e_{\rho,d}=\frac{1}{|\ws^\eta_n|}\sum_{g\in\ws^\eta_n}\Theta_{\rho,d}(g^{-1})g.$$
By Remark~\ref{remark:block idempotents}, we have $e_{\rho,d}\in\O\wa_n$, so the coefficients of $e_{\rho,d}$ outside $\wa_n$ are zero. For every $\lambda\in\Lambda_{\rho,d}$, the restriction ${\rm Res}_{\wa_n}^{\ws^\eta_n}(\chi_\lambda)$ is integer-valued by the spin-character value formulas recalled in \cite[Section 2, page 355]{Bessenrodt18}. Hence $\Theta_{\rho,d}(h)\in\mathbb Z$ for every $h\in\wa_n$, and the displayed formula shows that $e_{\rho,d}\in\mathbb Q\wa_n$. Consequently, $e_{\rho,d}\in\mathbb Z_p\wa_n$, as required.
\end{proof}

We set
$$B^{\eta,\rho,d}_{R_p}:=R_p\ws^\eta_ne_{\rho,d}~~{\rm and}~~(B^{\eta,\rho,d}_{R_p})_{\0}:=R_p\wa_ne_{\rho,d}.$$
As before, it is easy to see that $R_p \wa_n e_{\rho,d}$ is indeed the purely even subalgebra of the superalgebra $R_p\ws^\eta_ne_{\rho,d}$.

Fix a $\bar p$-core $\rho$ and a non-negative integer $d$. Let $r:=|\rho|$ and $n:=r+dp$. As in Theorem~\ref{theo: blocks of double covers of Sn}, we set $U:=[r]$ and $V:=[n]\setminus U$. For $i\in\{1,\cdots,d\}$ we also set $
V_i:=[r+ip]\setminus[r+(i-1)p].
$ 
In other words, we have:
\begin{equation}\label{equation:graph}
	\underbrace{1,\cdots,r}_U\,\overbrace{\underbrace{r+1,\cdots,r+p}_{V_1}\,
		\,\underbrace{r+p+1,\dots,r+2p}_{V_2}\,
		\cdots\,\underbrace{r+p(d-1)+1,\cdots,r+dp}_{V_d}}^V
\end{equation}

Let $D_i$ is a Sylow $p$-subgroup of $\ws^\eta_{V_i}$. With this notation, the following lemma follows easily from Theorem \ref{theo: blocks of double covers of Sn}:

\begin{lemma}[{\cite[Lemma 4.6]{KL25}}]\label{lemma:defect group abelian} 
Let $D$ be a defect group of the block algebra $R\ws^\eta_ne_{\rho,d}$. 
	Then $D$ is abelian if and only if $d<p$. In this case we can choose $D =D_1\times \dots \times D_d$. 
\end{lemma}

\begin{void}\label{void:T_w}
{\rm Keep the notation in (\ref{equation:graph}). Let $\iota:S_d\to S_V$ be the natural permutation action of $S_d$ on the set $\{V_1,\cdots,V_d\}$: more precisely, for any $1\leq i\leq d$, $1\leq j\leq p$ and $w\in S_d$
$$\iota(w)\cdot (r+(i-1)p+j)=r+(w(i)-1)p+j.$$
For each $w\in S_d$, choose a lift $T_w$ of $\iota(w)$ to $\ws^\eta_V$, and set $T_i:=T_{s_i}$ for all $1\leq i\leq d-1$. 

Note that every Sylow $p$-subgroup of $S_{V_i}$ lifts uniquely and isomorphically to a Sylow $p$-subgroup of $\ws^\eta_{V_i}$. Hence we can choose the $D_i$'s such that $T_w D_i T_w^{-1} = D_{w(i)}$, for all $1\leq i\leq d$ and $w \in S_d$. Indeed, one can first construct Sylow $p$-subgroups of each of the $S_{V_i}$'s that get permuted by the $\iota(w)$'s and then lift to $\ws^\eta_{V_i}$. By the construction of the $T_w$'s we also have $T_w \ws^\eta_{V_i} T_w^{-1} = \ws^\eta_{V_{w(i)}}$, for all $1\leq i\leq d$ and $w \in S_d$.

Note that $R \ws^\eta_{V_1}e_{z} \cong R \ws^\eta_p e_{z}$ via $
t_{r+m} e_z
\mapsto t_m e_z $. We now fix isomorphisms between the $R\ws^\eta_{V_i}e_{z}$'s. For each $1 < i\leq d$, we identify $R \ws^\eta_{V_i}e_{z}$ with $R \ws^\eta_{V_1}e_{z} \cong R \ws^\eta_p e_{z}$ via the isomorphism
\begin{align}\label{Sp_ident}
	R \ws^\eta_{V_1}e_{z} \to R\ws^\eta_{V_i}e_{z},\quad
	a \mapsto (-1)^{i|a|}T_{(1,i)}aT_{(1,i)}^{-1}.
\end{align}
(This does not depend on the choice of $T_{(1,i)}$, because $T_{(1,i)}$ is uniquely determined up to multiplication by $z$.) Through these isomorphisms and using (\ref{EZCommAlg}), we can identify the superalgebra
\begin{align}\label{algn:V1_Vd_ident}
	R \ws^\eta_{V_1, \cdots, V_d} e_z \cong R \ws^\eta_{V_1}e_z \otimes_R \cdots \otimes_R R \ws^\eta_{V_d}e_z
\end{align}
with $(R \ws^\eta_p e_z)^{\otimes d}$. Clearly (\ref{EZCommAlg}), (\ref{Sp_ident}) and (\ref{algn:V1_Vd_ident}) still hold with $R$ replaced by $R_p$ - because all coefficients involved in these isomorphisms are contained in $R_p$.
}
\end{void}

\begin{lemma}[a slight refinement of {\cite[Lemma 4.7]{KL25}}]\label{lemma:kappa}
	Keep the notation in \ref{void:T_w}. Up to modifying the choice of $T_i$ for $1\leq i\leq d-1$, the set $\{T_i\mid1\leq i\leq d-1\}$ satisfies the following equalities:
	\begin{enumerate}[{\rm (i)}]
		\item 
		$( T_i e_z)(a_1\otimes \cdots \otimes a_d)( T_i^{-1} e_z) = (-1)^{\sum_{j\neq i,i+1}|a_j|}({}^{s_i}(a_1\otimes \cdots \otimes a_d)),$
		for all $1\leq i\leq d-1$ and $a_1 \otimes \dots \otimes a_d\in (R \ws^\eta_p e_z)^{\otimes d}$.
		\item
		$( T_i e_z)^2 = \eta(-1)^{\frac{(p-1)}{2}}e_z,$ 
		for all $1\leq i\leq d-1$.
		\item
		$(T_i e_z)( T_l e_z) = -( T_l e_z)( T_i e_z),$ 
		for all $1\leq i,l\leq d-1$ with $|i-l|>1$.
		\item
		$( T_i e_zT_{i+1} e_z)^3 =e_z= ( T_{i+1} e_z T_i e_z)^3,$
		for all $1\leq i\leq d-2$.
	\end{enumerate}
	The result also holds for $R_p$ instead of $R$. 
\end{lemma}

\begin{proof} Equality (i) does not depend on the choice of $T_i$'s. Indeed, if $a_1 \otimes \cdots \otimes a_d \in (R \ws^\eta_p e_z)^{\otimes d}$, we have
\begin{align*}
	&( zT_i e_z)(a_1\otimes \cdots \otimes a_d)(z^{-1} T_i^{-1} e_z) = T_i(a_1\otimes \cdots \otimes a_d)T_i^{-1} \\
	= &  T_{(1,i)}T_{(1,i+1)}T_{(1,i)} (a_1\otimes \cdots \otimes a_d) T_{(1,i)}^{-1}T_{(1,i+1)}^{-1}T_{(1,i)}^{-1}.
\end{align*}
Equality (i) now follows from a direct calculation using (\ref{Sp_ident}).

Equality (iii) also does not depend on the $T_i$'s. Noting that $ze_z=-e_z$, equality (iii) follows immediately from (\ref{EZCommAlg}).

Recall that for each $1 \leq i\leq d-1$, $T_i$ is a lift of $\iota(s_i)$, where 
$$\iota(s_i)=(r+p(i-1)+1, r+pi+1)(r+p(i-1)+2, r+pi+2)\cdots (r+pi, r+pi+p).$$   So $T_i^2 = z^{p(p-1)/2}$ if $\ws^\eta_n=\ws_n^+$, and $T_i^2 = z\cdot z^{p(p-1)/2}$ if $\ws^\eta_n=\ws_n^-$. Noting that $ze_z=-e_z$, we have $(T_ie_z)^2 = (-1)^{p(p-1)/2}e_z$ if $\ws^\eta_n=\ws_n^+$, and $(T_ie_z)^2 = (-1)\cdot (-1)^{p(p-1)/2}e_z$ if $\ws^\eta_n=\ws_n^-$. In other words, we have $(T_ie_z)^2 = \eta(-1)^{\frac{(p-1)}{2}}e_z$. This proves equality (ii). Again, equality (ii) does not depend on the choice of $T_i$.

Finally, since $\pi_n(T_iT_{i+1})=\iota(s_is_{i+1})$ has order $3$, we have $(T_iT_{i+1})^3\in\{1,z\}$ for each $1\leq i\leq d-2$.
We arrange $(T_iT_{i+1})^3=1$ successively: at step $i$, keep $T_1,\cdots,T_i$ fixed and replace $T_{i+1}$ with $zT_{i+1}$ if $(T_iT_{i+1})^3=z$. These replacements leave identities (i)--(iii) and the relations obtained at previous steps unchanged.
Moreover,
$$(T_{i+1}T_i)^3=T_i^{-1}(T_iT_{i+1})^3T_i=1.$$
Multiplying these cubic relations by $e_z$ proves (iv).

%Finally, since the $\pi_n(T_i)$'s satisfy the braid relations, property (iv) holds up to a $z$, for each $1\leq i\leq d-2$. We may therefore change the choice of some of $T_i$'s to ensure that equality (iv) holds.  

Since all coefficients involved in this lemma are contained in $R_p$, the same argument shows that the result also holds for $R_p$ instead of $R$.      \end{proof}

\begin{void}\label{void:special subgroups and their blocks}
{\rm Fix a $\bar{p}$-core $\rho$ and a non-negative integer $d$. We adopt the notation in (\ref{equation:graph}), in particular, $r:=|\rho|$ and $n:=r+dp$. Let 
$$G:= \ws^\eta_n,~~L:=\ws^\eta_{r,p^d},~~N^{\rho,d}:= N_G(\ws^\eta_{V_1,\cdots,V_d}).$$
}
\end{void}
Note that $N^{\rho,d}$ must permute  $V_1,\cdots,V_d$. In particular, $N^{\rho,d}  \leq \ws^\eta_n$ is the subgroup generated by $L=\ws^\eta_{U,V_1,\cdots,V_d}$ and the set $\{T_w\mid w\in S_d\}$ where the notation $T_w$ is from \ref{void:T_w}. Recall the notation in (\ref{equation:super-group-algebra-in-doublecover}) that we will also use the subgroups $G_{\0}=\wa_n$, $N_{\0}^{\rho,d}=N^{\rho,d}\cap \wa_n$, etc.  We make the following identifications using (\ref{EZCommAlg}):
\begin{align}
\begin{split}\label{algn:group_tensor_isom}
%	R H e_z & \cong R \ws^\eta_{U \cup V_1\cup\cdots\cup V_{d-1}}e_z \otimes_R \ws^\eta_{V_p}e_z  \\
	R L e_z & \cong R \ws^\eta_U e_z \otimes_R R \ws^\eta_{V_1}e_z \otimes_R \cdots \otimes_R R \ws^\eta_{V_d}e_z.
\end{split}
\end{align}
Set $b:=e_{\rho,d}\in RGe_z$ and $$f=:e_{\rho,0}\otimes e_{\varnothing,1}^{(1)}\otimes\cdots\otimes e_{\varnothing,1}^{(d)}\in R Le_z\subseteq R N^{\rho,d}e_z.$$

By Proposition \ref{prop:descent of block idempotents}, we have $f\in R_pL\subseteq R_pN^{\rho,d}$. Since $N^{\rho,d}$ is generated by $L$ and $\{T_w\mid w\in S_d\}$, the isomorphism 
$$R Le_z\cong R\ws^\eta_re_z\otimes_R (R\ws^\eta_pe_z)^{\otimes d}$$
from (\ref{algn:group_tensor_isom}) extends to an isomorphism
\begin{align*}
	R N^{\rho,d}  e_z  \cong R \ws^\eta_r e_z \otimes_R (R\ws^\eta_p e_z \wr_s \T^\eta_{d,R}),\ 
	T_i e_z  \mapsto (1 \otimes \xi_i),
\end{align*}
where $T_i$'s are chosen to satisfy Lemma \ref{lemma:kappa}. In particular, we have the isomorphisms of superalgebras
\begin{align}
	\begin{split}\label{algn:ONf_ident}
		R L f & \cong R \ws^\eta_r e_{\rho,0} \otimes_R (R\ws^\eta_p e_{\varnothing,1})^{\otimes d} = B^{\eta,\rho,0}_R\otimes_R (B^{\eta,\varnothing,1}_R)^{\otimes d} \\
		R N^{\rho,d} f & \cong R \ws^\eta_r e_{\rho,0} \otimes_R (R\ws^\eta_p e_{\varnothing,1} \wr_s \T^\eta_{d,R})= B^{\eta,\rho,0}_R\otimes_R (B^{\eta,\varnothing,1}_R\wr_s \T^\eta_{d,R}).
	\end{split}
\end{align}

\begin{remark}\label{remark:ONf_ident-over-R0}
{\rm Since all coefficients involved in the isomorphisms (\ref{algn:ONf_ident}) are contained in the subring $R_p$ of $R$. The isomorphisms (\ref{algn:ONf_ident}) of $R$-superalgebras restricts to isomorphisms of $R_p$-superalgebras
\begin{align*}
	\begin{split}\label{algn:ONf_ident-over-R0+}
		R_p L f & \cong R_p \ws^\eta_r e_{\rho,0} \otimes_{R_p} (R_p\ws^\eta_p e_{\varnothing,1})^{\otimes d} = B^{\eta,\rho,0}_{R_p}\otimes_{R_p} (B^{\eta,\varnothing,1}_{R_p})^{\otimes d} \\
		R_p N^{\rho,d} f & \cong R_p \ws^\eta_r e_{\rho,0} \otimes_{R_p} (R_p\ws^\eta_p e_{\varnothing,1} \wr_s \T^\eta_{d,R_p}) = B^{\eta,\rho,0}_{R_p}\otimes_{R_p} (B^{\eta,\varnothing,1}_{R_p}\wr_s \T^\eta_{d,R_p}) .
	\end{split}
\end{align*}
}
\end{remark}

\begin{lemma}[{\cite[Lemmas 4.9 and 4.10]{KL25}}]\label{lemma:f is a block idempotent}
Let $0<d<p$.	As in \ref{void:T_w}, for $i=1,\cdots,d$, let $D_i$ be a Sylow $p$-subgroup of $\ws^\eta_{V_i}$ such that $T_w D_i T_w^{-1} = D_{w(i)}$, for all $i=1,\cdots,d$ and $w \in S_d$, and let
	\begin{equation*}\label{ED}
		D  := D_1 \times \cdots \times D_d \leq \ws^\eta_V.
	\end{equation*}
\begin{enumerate}[{\rm (i)}]
	\item Each of $RGb$, $RLf$ and $RN^{\rho,d}f$ is a block with defect group $D$.
	\item $RN^{\rho,d}f$ (resp. $RN_{\0}^{\rho,d}f$) is the Brauer correspondent of $RGb$ (resp. $RG_{\0}b$) in $N^{\rho,d}$ (resp. $N_{\0}^{\rho,d}$). 
\end{enumerate}
\end{lemma}

\begin{void}[{\cite[Definition 5.1]{KL25}}]\label{void:Rouquier core}
{\rm \textbf{Rouquier cores and RoCK blocks.} Let $d\in\mathbb{N}$.	A {\it $d$-Rouquier ${\bar p}$-core} $\rho$ is a $\bar{p}$-core such that ${\rm Ab}_\rho$ has the following properties:
	\begin{enumerate}[{\rm (i)}]
		\item The first runner has at least $d$ beads.
		\item The $(i+1)$-th runner has at least $d-1$ more beads than the $i$-th runner, for $1\leq i\leq (p-3)/2$.
	\end{enumerate}
	If $\rho$ is a $d$-Rouquier $\bar{p}$-core, we refer to $B^{\eta,\rho,d}_R$ as a {\it RoCK block} of {\it weight $d$} over $R$. 
}
\end{void}

For more information on RoCK blocks we refer to \cite[\S5]{KL25}.

\section{On Brauer correspondents, and proofs of Theorems \ref{theo:intro Brauer corr over Fp} and \ref{theorem:introduction Brauer corr over Fp splendid}}\label{section:On Brauer corr}

As in Section \ref{section:On spin blocks of ws_n and wa_n}, we continue to assume that $K$ (and hence $\O$) contains a primitive $(2n!)$-th root of unity in this section. We also assume that $-1$ and $2$ have square roots in $K$ (and hence $\O$); this condition is automatic if $n\geq 4$. Let $(R,R_p)\in\{(\O,\Z_p),(k,\FF_p)\}$. In this section we continue to adopt the notation in \ref{void:special subgroups and their blocks}; we always assume that $d>0$. Let $D$ be a Sylow $p$-subgroup of $\ws^\eta_{V}$. Then by Theorem \ref{theo:defect group and Brauer correspondent}, $RN_{\ws^\eta_n}(D) e_{\rho,0}$ is the Brauer correspondent of $RGb=R\ws^\eta_ne_{\rho,d}$, and by \cite[Theorem 5.2.10]{KL22}, $RN_{\wa_n}(D) e_{\rho,0}$ is the Brauer correspondent of $RG_{\0}b=R\wa_ne_{\rho,d}$. Note that we interpret $e_{\varnothing,0}$ as $e_z$. We set
\begin{align*}
	\begin{split}
b^{\eta,\rho,d}_{R}&:=RN_{\ws^\eta_n}(D) e_{\rho,0},~~(b^{\eta,\rho,d}_{R})_{\0}:=RN_{\wa_n}(D) e_{\rho,0},\\
b^{\eta,\rho,d}_{R_p}&:=R_pN_{\ws^\eta_n}(D) e_{\rho,0},~~(b^{\eta,\rho,d}_{R_p})_{\0}:=R_pN_{\wa_n}(D) e_{\rho,0}.
\end{split}
\end{align*}
Note that $(b^{\eta,\rho,d}_R)_{\0}$ (resp. $(b^{\eta,\rho,d}_{R_p})_{\0}$) is indeed the even part of the superalgebra $b^{\eta,\rho,d}_R$ (resp. $b^{\eta,\rho,d}_{R_p}$).

\begin{theorem}[{\cite[Proposition 5.2.13]{KL22}}]\label{prop: Brauer corr over k}
	Suppose that $0<d<p$. Let $D_i$ be a Sylow $p$-subgroup of $\ws^\eta_{V_i}$ chosen as in \ref{void:T_w}, and let $D=D_1\times\cdots\times D_d$. Then we have an isomorphism of $k$-superalgebras
	\begin{equation}\label{equation:Brauer corr over k}
		b^{\eta,\rho,d}_k\cong k\ws^\eta_{U}e_{\rho,0}\otimes_k (b^{\eta,\varnothing,1}_k\wr_s \T_{d,k})\cong  B^{\eta,\rho,0}_k\otimes_k (b^{\eta,\varnothing,1}_k\wr_s \T_{d,k}).
	\end{equation}	
\end{theorem}

\begin{theorem}\label{theo: Brauer corr over Fp}
	Suppose that $0<d<p$. Let $D_i$ be a Sylow $p$-subgroup of $\ws^\eta_{V_i}$ chosen as in \ref{void:T_w}, and let $D=D_1\times\cdots\times D_d$. Then we have an isomorphism of $\FF_p$-superalgebras
	\begin{equation}\label{equation:Brauer corr over Fp}
		b^{\eta,\rho,d}_{\FF_p}\cong \FF_p\ws^\eta_{U}e_{\rho,0}\otimes_{\FF_p} (b^{\eta,\varnothing,1}_{\FF_p}\wr_s \T^\eta_{d,\FF_p})\cong  B^{\eta,\rho,0}_{\FF_p}\otimes_{\FF_p} (b^{\eta,\varnothing,1}_{\FF_p}\wr_s \T^\eta_{d,\FF_p}).
	\end{equation}	
\end{theorem}

\begin{proof}
By Proposition \ref{prop:descent of block idempotents}, we have $e_{\rho,0}\in\FF_p\wa_U$. The nontrivial orbits of $\pi_n(D)$ on $[n]$ are $V_1,\cdots,V_d$, and its fixed-point set is $U$. Recall that $\pi_n$ denotes the canonical surjective group homomorphism $\ws^\eta_n\to S_n$; see \ref{void:double covers}. Since the $T_w$'s permute the $D_i$'s as in \ref{void:T_w}, we have the disjoint coset decomposition
\[
N_{\ws^\eta_n}(D)=\bigsqcup_{w\in S_d}\bigl(N_{\ws^\eta_n}(D)\cap L\bigr)T_w.
\]
Moreover, $N_{\ws^\eta_n}(D)\cap L$ is generated by $\ws^\eta_U$ and $N_{\ws^\eta_{V_i}}(D_i)$ for $1\leq i\leq d$.
Since conjugation by $T_{(1,i)}$ sends $D_1$ to $D_i$, the identifications (\ref{EZCommAlg}) and (\ref{Sp_ident}) induce an isomorphism of $\FF_p$-superalgebras
$$\FF_p\bigl(N_{\ws^\eta_n}(D)\cap L\bigr)e_{\rho,0}
\cong \FF_p\ws^\eta_Ue_{\rho,0}\otimes_{\FF_p}
\bigl(b^{\eta,\varnothing,1}_{\FF_p}\bigr)^{\otimes d}.$$

Choose the lifts $T_i$ to satisfy the relations in
Lemma \ref{lemma:kappa}.
These relations give the defining relations of
$b^{\eta,\varnothing,1}_{\FF_p}\wr_s\T^\eta_{d,\FF_p}$.
Also, the elements supported on $U$ and $V$ satisfy the super tensor product relations by (\ref{EZCommAlg}).
Hence the inverse of the above isomorphism extends to a homomorphism of $\FF_p$-superalgebras
\[
\FF_p\ws^\eta_Ue_{\rho,0}\otimes_{\FF_p}
\bigl(b^{\eta,\varnothing,1}_{\FF_p}\wr_s\T^\eta_{d,\FF_p}\bigr)
\longrightarrow b^{\eta,\rho,d}_{\FF_p},
\]
sending $e_{\rho,0}\otimes\xi_i$ to $T_i e_{\rho,0}$.
This homomorphism is surjective because $N_{\ws^\eta_n}(D)$ is generated by $N_{\ws^\eta_n}(D)\cap L$ and the $T_i$'s.
By Proposition \ref{proposition:epsilon-AwrT-graded} (iii) and the coset decomposition above, both sides have dimension
\[
d!\,\dim_{\FF_p}\bigl(\FF_p\ws^\eta_Ue_{\rho,0}\bigr)
\bigl(\dim_{\FF_p}b^{\eta,\varnothing,1}_{\FF_p}\bigr)^d.
\]
Thus the homomorphism is an isomorphism. The second isomorphism in (\ref{equation:Brauer corr over Fp}) follows from the definition of $B^{\eta,\rho,0}_{\FF_p}$.
\end{proof}

\begin{void}
	{\rm Keep the notation of Theorem \ref{theo: Brauer corr over Fp}. We have the following isomorphisms all of which are obtained by restricting either the isomorphism (\ref{equation:Brauer corr over Fp}) or the isomorphisms in Remark \ref{remark:ONf_ident-over-R0} to appropriate subalgebras:
		\begin{align}
			\begin{split}\label{align:restriction of isomorphisms}
				\FF_p(N_G(D)\cap \ws^\eta_V)e_{\varnothing,0}^{\otimes d} &\cong b^{\eta,\varnothing,1}_{\FF_p}\wr_s \T^\eta_{d,\FF_p}\\
				\FF_p(N^{\rho,d}\cap \ws^\eta_V)e_{\varnothing,1}^{\otimes d}&\cong B^{\eta,\varnothing,1}_{\FF_p}\wr_s \T^\eta_{d,\FF_p}
			\end{split}
		\end{align}
	}
\end{void}

We realise \cite[Lemma 5.2.19]{KL22} over $\FF_p$: 

\begin{lemma}\label{lemma:descent of splendid Rickard equivalence for Ap} Let $P$ be a Sylow $p$-subgroup of $\wa_p$. There is a splendid Rickard complex of $(\FF_p\wa_pe_{\varnothing,1},\FF_pN_{\wa_p}(P)e_{\varnothing,0})$-bimodules that extends to a complex of $(\FF_p\ws^\eta_pe_{\varnothing,1},\FF_pN_{\ws^\eta_p}(P)e_{\varnothing,0})_{C_2}$-modules.
\end{lemma}

\begin{proof} 
First suppose that $p=3$. Then $\wa_3\cong C_6$, and $P$ is normal in $\ws^\eta_3$. Hence $N_{\wa_3}(P)=\wa_3$ and $N_{\ws^\eta_3}(P)=\ws^\eta_3$. Moreover, $e_{\varnothing,1}=e_{\varnothing,0}=e_z$. The bimodule $\FF_3\wa_3e_z$, concentrated in degree zero, induces a splendid Rickard equivalence with itself. Its action extends to
$(\FF_3\ws^\eta_3e_z,\FF_3\ws^\eta_3e_z)_{C_2}$ by letting $se_z\otimes s^{-1}e_z$ act as $x\mapsto sxs^{-1}$, for $s\in\ws^\eta_3\setminus\wa_3$.
We may therefore assume that $p\geq5$. In \cite[Theorem 1.10]{Kessar_Linckelmann}, Kessar and Linckelmann construct a splendid Rickard complex for blocks with a cyclic defect group over a non-splitting field, hence a splendid Rickard complex  $X$ of $(\FF_p\wa_pe_{\varnothing,1},\FF_pN_{\wa_p}(P)e_{\varnothing,0})$-bimodules. We review the construction of $X$ from the proof of \cite[Theorem 1.10]{Kessar_Linckelmann}. Note that $|P|=p$ and $O_p(\wa_p)=1$. The complex $X$ is of the form 
\begin{equation}\label{equation:splendid Rickard complex for Ap}
	X:=\cdots \to 0\to N\xrightarrow{\varphi} M\to 0\to \cdots
\end{equation}   
where $M$ is in degree zero. Moreover $M$ is up to isomorphism the unique nonprojective indecomposable summand of the $(\FF_p\wa_pe_{\varnothing,1},\FF_pN_{\wa_p}(P)e_{\varnothing,0})$-bimodule $e_{\varnothing,1}\FF_p\wa_pe_{\varnothing,0}$.  Let $Z\to M$ be a projective cover of the $(\FF_p\wa_pe_{\varnothing,1},\FF_pN_{\wa_p}(P)e_{\varnothing,0})$-bimodule $M$. Then $N$ is a direct summand of $Z$ and $\varphi$ is the restriction of the map $Z\to M$ to $N$; moreover, if we write $Z=N\oplus N'$, then $N$ and $N'$ have no nonzero isomorphic direct summand.           

As in the proof of \cite[Lemma 5.2.19]{KL22}, we transport our setup to the group algebra setting. Set
\begin{align*}
	\mathcal{A} = \langle\wa_p \times N_{\wa_p}(P),(s,s)\rangle\leq \ws^\eta_p \times N_{\ws^\eta_p}(P),
\end{align*}
for some fixed $s\in N_{\ws^\eta_p}(P)\backslash N_{\wa_p}(P)$. We have an isomorphism of $\FF_p$-superalgebras
\begin{align*}
	\FF_p \mathcal{A} (e_{\varnothing,1}\otimes e_{\varnothing,0})&\cong (\FF_p\ws^\eta_pe_{\varnothing,1},\FF_pN_{\ws^\eta_p}(P)e_{\varnothing,0})_{C_2},\\ 
	(a_1, a_2)(e_{\varnothing,1}\otimes e_{\varnothing,0})&\mapsto a_1 e_{\varnothing,1}\otimes a_2^{-1}e_{\varnothing,0}, \\
	(s,s)(e_{\varnothing,1}\otimes e_{\varnothing,0})&\mapsto s e_{\varnothing,1}\otimes s^{-1} e_{\varnothing,0},
\end{align*}
for all $a_1\in \wa_p$ and $a_2\in N_{\wa_p}(P)$. We identify $\FF_p\wa_pe_{\varnothing,1}\otimes_{\FF_p} (\FF_pN_{\wa_p}(P)e_{\varnothing,0})^{\rm sop}$ with the subalgebra $\FF_p(\wa_p \times N_{\wa_p}(P))(e_{\varnothing,1}\otimes e_{\varnothing,0})$ of $\FF_p \mathcal{A} (e_{\varnothing,1}\otimes e_{\varnothing,0})$.

In this setting, we need to show that the complex $X$ from (\ref{equation:splendid Rickard complex for Ap}) extends to a complex of $\FF_p\mathcal{A}$-modules. Clearly the $\FF_p(\wa_p \times N_{\wa_p}(P))$-module $e_{\varnothing,1}\FF_p \wa_pe_{\varnothing,0}$ extends to an $\FF_p\mathcal{A}$-module via $(s,s)\cdot x=sxs^{-1}$ for all $x\in e_{\varnothing,1}\FF_p \wa_pe_{\varnothing,0}$. If we write $e_{\varnothing,1}\FF_p \wa_pe_{\varnothing,0}=M\oplus M'$, then (by the construction of $X$) $M$ is the nonprojective part and $M'$ is the projective part of $e_{\varnothing,1}\FF_p \wa_pe_{\varnothing,0}$. Conjugation by $(s,s)$ preserves projectivity. Thus, by the Krull--Schmidt theorem, the isomorphism classes of the
nonprojective part $M$ and the projective part $M'$ are both fixed by $(s,s)$; that is, $(s,s)M\cong M$ and $(s,s)M'\cong M'$ as $\FF_p(\wa_p \times N_{\wa_p}(P))$-modules. By \cite[Theorem 3.11]{Da84}, $M$ and $M'$ extend to $\FF_p\mathcal{A}(e_{\varnothing,1}\otimes e_{\varnothing,0})$-modules. Let $\tilde{M}$ be such an extension of $M$ to $\FF_p\mathcal{A}(e_{\varnothing,1}\otimes e_{\varnothing,0})$. Let $\tilde{Z}$ be a projective cover of the $\FF_p\mathcal{A}$-module $\tilde{M}$. Since $2$ is invertible in $\FF_p$, we have $$\tilde{Z}/J(\FF_p\mathcal{A})\tilde{Z}\cong\tilde{M}/J(\FF_p\mathcal{A})\tilde{M}\cong\tilde{M}/J(\FF_p(\wa_p \times N_{\wa_p}(P)))\tilde{M}\cong M/{\rm rad}(M)$$ as $\FF_p(\wa_p \times N_{\wa_p}(P))$-modules; see e.g. \cite[Theorem 1.11.10]{Lin18a} for the second isomorphism. So ${\rm Res}_{\wa_p \times N_{\wa_p}(P)}^{\mathcal{A}}(\tilde{Z})$ is a projective cover of $M$, and hence it is isomorphic to $Z$. 

Next we show that $(s,s)N\cong N$ and $(s,s)N'\cong N'$ as $\FF_p(\wa_p \times N_{\wa_p}(P))$-modules. Choose a finite extension $k_1$ of $\FF_p$ which is a splitting field for $\wa_p\times N_{\wa_p}(P)$. As recalled in the proof of \cite[Lemma 5.2.19]{KL22}, the Brauer tree of $k_1\wa_pe_{\varnothing,1}$ has an exceptional vertex of multiplicity $2$.
Since $(s,s)M\cong M$, conjugation by $(s,s)$ permutes the set of isomorphism classes of indecomposable summands of $k_1\otimes_{\FF_p}Z$. The summands $k_1\otimes_{\FF_p}N$ and $k_1\otimes_{\FF_p}N'$ of $k_1\otimes_{\FF_p}Z$ correspond to the sets $I_1$ and $I_0$ in \cite[page 206]{Kessar_Linckelmann}, respectively. The sets $I_0$ and $I_1$ there are invariant under the action of $s$, because conjugation by $s$ fixes the exceptional vertex and preserves distances in the Brauer tree. Hence
$$k_1\otimes_{\FF_p}(s,s)N \cong k_1\otimes_{\FF_p}N~~~{\rm and}~~~k_1\otimes_{\FF_p}(s,s)N' \cong k_1\otimes_{\FF_p}N'.$$
Restricting scalars to $\FF_p$ and applying the Krull--Schmidt theorem implies $(s,s)N\cong N$ and $(s,s)N'\cong N'$ as $\FF_p(\wa_p \times N_{\wa_p}(P))$-modules, as claimed.

Since $(s,s)N\cong N$ and $(s,s)N'\cong N'$ as $\FF_p(\wa_p \times N_{\wa_p}(P))$-modules, by \cite[Theorem 3.11]{Da84}, there is an $\FF_p\mathcal{A}$-direct summand $\tilde{N}$ of $\tilde{Z}$ such that $${\rm Res}_{\wa_p \times N_{\wa_p}(P)}^{\mathcal{A}}\tilde{N}\cong N.$$ Now let $\tilde{\varphi}$ be the restriction of the map $\tilde{Z}\to \tilde{M}$ to $\tilde{N}$. Then the complex
$$\tilde{X}:=\cdots \to 0\to \tilde{N}\xrightarrow{\tilde{\varphi}} \tilde{M}\to 0\to \cdots$$
is an extension of $X$ to $\FF_p\mathcal{A}(e_{\varnothing,1}\otimes e_{\varnothing,0})$.
\end{proof}

In \cite[Proposition 5.2.21]{KL22}, Kleshchev and Livesey proved the following result for Brauer correspondents of blocks of $\ws^\eta_n$ and $\wa_n$ with abelian defect groups.

\begin{theorem}[{\cite[Proposition 5.2.33]{KL22}}]\label{prop:Brauer corr over k splendid} Suppose that $0<d<p$. Let $D_i$ be a Sylow $p$-subgroup of $\ws^\eta_{V_i}$ chosen as in \ref{void:T_w}, and let $D=D_1\times\cdots\times D_d$. Then $kN_{\ws^\eta_n}(D)e_{\rho,0}$, i.e. $b^{\eta,\rho,d}_k$, is splendidly Rickard equivalent to $kNf$, and $kN_{\wa_n}(D)e_{\rho,0}$, i.e. $(b^{\eta,\rho,d}_k)_{\0}$, is splendidly Rickard equivalent to $kN_{\0}f$.
\end{theorem}

Now we realise Theorem \ref{prop:Brauer corr over k splendid} over $\FF_p$:

\begin{theorem}\label{theorem:Brauer corr over Fp splendid}
 Suppose that $0<d<p$. Let $D_i$ be a Sylow $p$-subgroup of $\ws^\eta_{V_i}$ chosen as in \ref{void:T_w}, and let $D=D_1\times\cdots\times D_d$. Then $\FF_pN_G(D)e_{\rho,0}$, i.e. $b^{\eta,\rho,d}_{\FF_p}$, is splendidly Rickard equivalent to $\FF_pN^{\rho,d}f$, and $\FF_pN_{G_{\0}}(D)e_{\rho,0}$, i.e. $(b^{\eta,\rho,d}_{\FF_p})_{\0}$, is splendidly Rickard equivalent to $\FF_pN_{\0}^{\rho,d}f$.
\end{theorem}

\begin{proof} By Remark \ref{remark:ONf_ident-over-R0}, we have 
$$\FF_p N^{\rho,d} f  \cong B^{\eta,\rho,0}_{\FF_p} \otimes_{\FF_p} (B^{\eta,\varnothing,1}_{\FF_p} \wr_s \T^\eta_{d,\FF_p}).$$ 
Let $P$ be a Sylow $p$-subgroup of $\ws^\eta_p$. The block algebra $\FF_p\wa_pe_{\varnothing,1}$ has a cyclic defect group $P$ and Brauer correspondent $\FF_pN_{\wa_p}(P)e_{\varnothing,0}$. By Lemma \ref{lemma:descent of splendid Rickard equivalence for Ap}, there is a splendid Rickard complex $X$ of $((B^{\eta,\varnothing,1}_{\FF_p})_{\0},(b^{\eta,\varnothing,1}_{\FF_p})_{\0})$-bimodule that extends to a complex $\tilde{X}$ of $(B^{\eta,\varnothing,1}_{\FF_p},b^{\eta,\varnothing,1}_{\FF_p})_{C_2}$-modules. Let $H$ be the subgroup of $G$ generated by $\ws^\eta_U$ and all $\wa_{V_i}$'s. Now (by Propositions \ref{prop:boxtimes Rickard equivalences}, \ref{prop:boxtimes splendid} and Remark \ref{remark: usual algebras as superalgebras}) $X^{\otimes d}$ induces a splendid Rickard equivalence between
$$\FF_p(H\cap \ws^\eta_V)e_{\varnothing,1}^{\otimes d}\cong \FF_p\wa_{V_1}e_{\varnothing,1}\otimes_{\FF_p}\cdots\otimes_{\FF_p} \FF_p\wa_{V_d}e_{\varnothing,1}\cong (B^{\eta,\varnothing,1}_{\0})^{\otimes d}$$
and 
$$\FF_p(N_G(D)\cap H\cap \ws^\eta_{V})e_{\varnothing,0}\cong \FF_pN_{\wa_{V_1}}(D_1)e_{\varnothing,0}\otimes_{\FF_p}\cdots\otimes_{\FF_p} \FF_pN_{\wa_{V_d}}(D_d)e_{\varnothing,0}\cong (b^{\eta,\varnothing,1}_{\0})^{\otimes d}.$$
(Note that by Remark \ref{remark: usual algebras as superalgebras}, Theorem \ref{prop:boxtimes splendid} applies to usual algebras and bimodules: in this case the box tensor products in Theorem \ref{prop:boxtimes splendid} are the usual tensor products.)

We now adopt the identifications (\ref{align:restriction of isomorphisms}). By Proposition \ref{prop:epsilon-ext_der_wreath}, $X^{\otimes d}$ extends to a complex $\tilde{Y}$ of modules for the group algebra of
$$\langle (H\cap \ws^\eta_V)\times (N_G(D)\cap H\cap \ws^\eta_{V}), (u_i,u_i), (T_w,T_w) \rangle \leq (N^{\rho,d}\cap \ws^\eta_V)\times(N_G(D)\cap \ws^\eta_V),$$
where $i$ ranges over $\{1,\cdots,d\}$ with $u_i\in N_{\ws^\eta_{V_i}}(D_i)\backslash N_{\wa_{V_i}}(D_i)$, and $w$ ranges over $S_d$ with $T_w$ chosen as in Lemma \ref{lemma:kappa}. We then induce $\tilde{Y}$ up to
$$(N^{\rho,d}\cap \ws^\eta_V)\times(N_G(D)\cap \ws^\eta_V)$$ 
to obtain a complex $Y$. Since $d<p$, we have $p\nmid |C_2\wr S_d|$. It follows from Theorem \ref{theorem:Marcus CA Theorem 4.8} that $Y$ induces a $(C_2\wr S_d)$-graded Rickard equivalence between 
$$\FF_p(N^{\rho,d}\cap \ws^\eta_V)e_{\varnothing,1}^{\otimes d}\cong B^{\eta,\varnothing,1}_{\FF_p} \wr_s \T^\eta_{d,\FF_p}$$
and 
$$\FF_p(N_G(D)\cap \ws^\eta_V)e_{\varnothing,0}^{\otimes d}\cong b^{\eta,\varnothing,1}_{\FF_p}\wr_s \T^\eta_{d,\FF_p}.$$
By Proposition \ref{prop:splendid modules ind and res} (iii), $Y$ is a splendid complex.           

Next consider the group homomorphism 
$$C_2\wr S_d\to \Z/2\Z,~~g_i\mapsto \1,~~s_j\mapsto \1$$
for $1\leq i\leq d$ and $1\leq j\leq d-1$, where $g_i$ is the generator of the $i$-th $C_2$ and $s_j$ is the appropriate elementary transposition in $S_d$. Let $\mathcal{K}$ be the kernel of this homomorphism. Since $B^{\eta,\varnothing,1}_{\FF_p} \wr_s \T^\eta_{d,\FF_p}$ is a $C_2\wr S_d$-graded algebra, we use the notation $(B^{\eta,\varnothing,1}_{\FF_p} \wr_s \T^\eta_{d,\FF_p})_g$ to denote the $g$-component of $B^{\eta,\varnothing,1}_{\FF_p} \wr_s \T^\eta_{d,\FF_p}$ for any $g\in C_2\wr S_d$. Then we have
$$(B^{\eta,\varnothing,1}_{\FF_p} \wr_s \T^\eta_{d,\FF_p})_{\0}=\bigoplus_{g\in \mathcal{K}}(B^{\eta,\varnothing,1}_{\FF_p} \wr_s \T^\eta_{d,\FF_p})_g.$$
Since $Y$ is a $C_2\wr S_d$-graded bimodule, $Y$ can be regarded as a bisupermodule ($\Z/2\Z$-graded bimodule) via
$$Y_{\0}=\bigoplus_{g\in \mathcal{K}} Y_g~~{\rm and}~~Y_{\1}=\bigoplus_{g\in (C_2\wr S_d)\backslash \mathcal{K}}Y_g.$$
By the definition of a $(C_2\wr S_d)$-graded Rickard equivalence, we see that $Y$ (considered as a bisupermodule) induces a $\Z/2\Z$-graded Rickard equivalence between  
$$\FF_p(N^{\rho,d}\cap \ws^\eta_V)e_{\varnothing,1}^{\otimes d}\cong B^{\eta,\varnothing,1}_{\FF_p} \wr_s \T^\eta_{d,\FF_p}$$
and 
$$\FF_p(N_G(D)\cap \ws^\eta_V)e_{\varnothing,0}^{\otimes d}\cong b^{\eta,\varnothing,1}_{\FF_p}\wr_s \T^\eta_{d,\FF_p}.$$

It is straightforward to check that the $(B^{\eta,\rho,0}_{\FF_p},B^{\eta,\rho,0}_{\FF_p})$-bisupermodule $B^{\eta,\rho,0}_{\FF_p}$ induces a Morita self-superequivalence of $B^{\eta,\rho,0}_{\FF_p}$; this bimodule is splendid because it is isomorphic to a direct summand of $\FF_p\ws^\eta_r\cong {\rm Ind}_{\Delta \ws^\eta_r}^{\ws^\eta_r\times \ws^\eta_r}(\FF_p)$. Hence by Proposition \ref{prop:boxtimes Rickard equivalences}, the complex $B^{\eta,\rho,0}_{\FF_p}\boxtimes_{\FF_p} Y$ of $\Z/2\Z$-graded $(B^{\eta,\rho,0}_{\FF_p} \otimes_{\FF_p} \FF_p(N^{\rho,d}\cap \ws^\eta_V)e_{\varnothing,1}^{\otimes d})$-$(B^{\eta,\rho,0}_{\FF_p}\otimes_{\FF_p} \FF_p(N_G(D)\cap \ws^\eta_V)e_{\varnothing,0}^{\otimes d})$-bimodule induces a $\Z/2\Z$-graded Rickard equivalence between
$$\FF_p N^{\rho,d} f\cong B^{\eta,\rho,0}_{\FF_p} \otimes_{\FF_p} \FF_p(N^{\rho,d}\cap \ws^\eta_V)e_{\varnothing,1}^{\otimes d}$$
and
$$\FF_pN_G(D)e_{\rho,0}\cong B^{\eta,\rho,0}_{\FF_p} \otimes_{\FF_p} \FF_p(N_G(D)\cap \ws^\eta_V)e_{\varnothing,0}^{\otimes d}.$$
By Proposition \ref{prop:boxtimes splendid}, $B^{\eta,\rho,0}_{\FF_p}\boxtimes_{\FF_p} Y$ is a splendid complex of $(\FF_p N^{\rho,d} f, \FF_pN_G(D)e_{\rho,0})$-bimodules.

Now by Theorem \ref{theorem:Marcus CA Theorem 4.7}, $(B^{\eta,\rho,0}_{\FF_p}\boxtimes_{\FF_p} Y)_{\0}$ induces a Rickard equivalence between
$\FF_p N^{\rho,d}_{\0}f$ and $\FF_pN_{G_{\0}}(D)e_{\rho,0}$.
Since $(B^{\eta,\rho,0}_{\FF_p}\boxtimes_{\FF_p} Y)_{\0}$ is isomorphic to a direct summand of the complex
$${\rm Res}_{N^{\rho,d}_{\0}\times N_{G_{\0}}(D)}^{N^{\rho,d}\times N_G(D)}(B^{\eta,\rho,0}_{\FF_p}\boxtimes_{\FF_p} Y),$$
by Proposition \ref{prop:splendid modules ind and res} (ii), $(B^{\eta,\rho,0}_{\FF_p}\boxtimes_{\FF_p} Y)_{\0}$ is splendid.    \end{proof}

\section{Morita equivalences with endopermutation source over arbitrary fields}\label{s: Morita equivalences with endopermutation source over arbitrary fields}

In this section, we extend some known results on Morita equivalences with endopermutation source over splitting $p$-modular systems to arbitrary $p$-modular systems. These results will be used in the next section. In this section $p$ is not necessarily being odd.

\begin{definition}[{\cite[Definition 8.5.1]{Lin18b}}]\label{defi:fusion system of block}
{\rm Let $R\in \{\O,k\}$, $G$ a finite group, $b$ a block of $RG$ and $(P,e)$ a maximal $RGb$-Brauer pair. (See \cite[\S 6.3]{Lin18b} for details on Brauer pairs.) For any subgroup $Q$ of $P$ denote by $e_Q$ the unique block of $kC_G(Q)$ such that $(Q,e_Q)\leq (P,e)$. We define a category $\F=\F_{(P,e)}(RGb)$ as follows: the objects of $\F$ are the subgroups of $P$, and for any two subgroups $Q$ and $S$ of $P$, the morphism set $\Hom_\F(Q,S)$ is the set of all group homomorphism $\varphi:Q\to S$ for which there exists an element $x\in G$ satisfying $\varphi(u)=xux^{-1}$ for all $u\in Q$ and such that ${}^x(Q,e_Q)\leq (S,e_S)$.
}	
\end{definition}

\begin{remark}\label{remark:extension of fusion}
{\rm Keep the notation of Definition \ref{defi:fusion system of block}. It is well known that the category $\F=\F_{(P,e)}(RGb)$ satisfies the extension axiom \cite[Definition 8.1.6 (II)]{Lin18b}; see the proof of \cite[Theorem 8.5.2]{Lin18b}. Only for the Sylow axiom \cite[Definition 8.1.6 (I)]{Lin18b}, one needs to assume that $k$ is large enough. Assume that $P$ is abelian. Then by the extension axiom \cite[Definition 8.1.6 (II)]{Lin18b}, for any subgroup $Q,S$ of $P$ and any morphism $\varphi:Q\to S$ in $\F$, there exists a morphism $\psi\in\Aut_\F(P)$ such that $\psi|_Q=\varphi$.
}	
\end{remark}

\begin{definition}[cf. {\cite[Definitions 9.9.1, 9.9.2]{Lin18b}}]
{\rm We refer to \cite[Definition 7.3.1]{Lin18b} for the definition of an endopermutation module and to \cite[Definition 7.11.1]{Lin18b} for the definition of an endosplit $p$-permutation resolution. Keep the notation of Definition \ref{defi:fusion system of block}. Let $Q$ be a subgroup of $P$ and $V$ an endopermutation $RQ$-module (resp. an endosplit $p$-permutation resolution of some $RQ$-module). We say that $V$ is {\it $\F$-stable} if for any subgroup $S$ of $Q$ and any morphism $\varphi:S\to Q$ in $\F$, the sets of isomorphism classes of indecomposable direct summands with vertex $S$ of the $RS$-modules (resp. complexes of $RS$-modules) $\Res_S^Q(V)$ and ${}_\varphi V$ are equal (including the possibility that both sets may be empty).}
\end{definition}

\begin{theorem}[A slight generalization of {\cite[Theorem 8.7.1 (i)]{Lin18b}} to arbitrary fields]\label{theo:8.7.1 (i)}
Let $R\in \{\O,k\}$, $G$ a finite group and $b$ a block of $RG$ with a defect group $P$. Let $i\in (RGb)^P$ be an almost source idempotent of $b$; see \cite[Definition 6.4.3]{Lin18b}. Let $e$ be the block of $RC_G(P)$ determined by $i$ and let $\F=\F_{(P,e)}(RGb)$. Let $Q$ and $S$ be subgroups of $P$. Then every indecomposable direct summand of $iRGi$ as an $RQ$-$RS$-bimodule is isomorphic to $RQ\otimes_{R T} {}_\varphi RS$ for some subgroup $T$ of $Q$ and some morphism $\varphi:T\to S$ in $\F$.
\end{theorem}	

\begin{proof}
One checks that the proof of \cite[Theorem 8.7.1 (i)]{Lin18b} does not require $k$ to be large enough.
\end{proof}

\begin{proposition}[A slight generalization of {\cite[9.11.4 (i),(ii),(iii) and 9.11.5 (i)]{Lin18b}} to arbitrary fields]\label{prop:9.11.4+9.11.5}
Let $R\in \{\O,k\}$. Let $G$, $H$ be finite groups, $b$ a block of $RG$ with a defect group $P$ and $c$ a block of $RH$ with a defect group $Q$. Let $i\in (RGb)^P$ and $j\in (RHc)^Q$ be source idempotents. Let $(Q,f)$ be the maximal $kHc$-Brauer pair determined by $j$ and let $\F=\F_{(Q,f)}(RHc)$. Assume that there is a group isomorphism $\varphi:P\xrightarrow{\cong} Q$. Set $\Delta\varphi:=\{(u,\varphi(u))\mid u\in P\}$. Let $V$ be an indecomposable endopermutation $R \Delta\varphi$-module with vertex $\Delta\varphi$. We also view $V$ as an $RQ$-module via the isomorphism $Q\cong \Delta \varphi$, $v\mapsto (\varphi^{-1}(v),v)$. Suppose that $V$ has an $\F$-stable endosplit $p$-permutation resolution $Y_V$. Set 
$$U=RGi\otimes_{RP}\Ind_{\Delta \varphi}^{P\times Q}(V)\otimes_{RQ} jRH,$$
$$Y=RGi\otimes_{RP}\Ind_{\Delta \varphi}^{P\times Q}(Y_V)\otimes_{RQ} jRH.$$
The following hold:
\begin{enumerate}[{\rm (i)}]
	\item The complex $Y$ is a splendid complex of $RGb$-$RHc$-bimodules, and $Y$ has homology  concentrated in degree zero, isomorphic to $U$.
	\item The complex $RGb$-$RGb$-bimodule $Y\otimes_{RHc} Y^*$ is split.
	\item Let $M$ be a direct summand of $U$. Then there is a direct summand $X$ of $Y$ such that $X$ has homology concentrated in degree zero, isomorphic to $M$.
	\item Keep the Notation of (iii). If $M$ is indecomposable and induces a Morita equivalence between $RGb$ and $RHc$, then $X$ is a splendid Rickard complex of $RGb$-$RHc$-bimodules.
\end{enumerate}
\end{proposition}	

\begin{proof}
Without loss of generality we may identify $P$ and $Q$ via $\varphi$. In other words, we may assume that $P=Q$ and $\varphi={\rm id}_P$. Then the proposition follows from the proof of \cite[9.11.4 (i),(ii),(iii) and 9.11.5 (i)]{Lin18b}: Note that since \cite[Theorem 8.7.1 (i)]{Lin18b} does not require $k$ to be large enough (see Theorem \ref{theo:8.7.1 (i)}), the proof there for these four statements is independent of the size of $k$. We also note that for these four statements, the proof there does not use the assumption that $i$ and $j$ determine the same fusion category. We remind the reader that even if in the proof of \cite[9.11.5 (i)]{Lin18b}, it only shows that $X\otimes_{RHc}X^*$ is homotopy equivalent to $RGb$. By \cite[Theorem 4.14.16]{Lin18a} we obtain that $X^*\otimes_{RGb}X$ is homotopy equivalent to $RHc$.
\end{proof}

\section{The Morita superequivalence bimodule, and proofs of Theorems \ref{theo:main1} and \ref{theo:main2}}\label{section:The Morita superequivalence bimodule}
\begin{void}\label{void:s5-notation}
{\rm Throughout this section we set $d$ to be an integer with $1\leq d<p$ and $\rho$ a $d$-Rouquier $\bar{p}$-core. As in Section \ref{section:On spin blocks of ws_n and wa_n}, we continue to assume that $K$ (and hence $\O$) contains a primitive $(2n!)$-th root of unity in this section. We also assume that $-1$ and $2$ have square roots in $K$ (and hence $\O$); this condition is automatic if $n\geq 4$. Let $(R,R_p)\in\{(\O,\Z_p),(k,\FF_p)\}$. We may assume that $R$ is finitely generated as $R_p$-module. We will adopt the notation in \ref{void:special subgroups and their blocks}, in particular, $r=|\rho|$, $n=r+dp$, 
\begin{align*}
	&G=\ws^\eta_n,\quad L= \ws^\eta_{U,V_1,\dots,V_d},
	\quad N^{\rho,d}= N_G(\ws^\eta_{V_1,\cdots,V_d}),\quad   
	\\
	&b=e_{\rho,d}\in RG,\quad  f= e_{\rho,0} \otimes e_{\varnothing,1}^{(1)} \otimes \dots \otimes e_{\varnothing,1}^{(d)} \in R L e_z,
\end{align*}
and the defect group 
$$D=D_1\times\dots\times D_d$$ is chosen as in Lemma \ref{lemma:f is a block idempotent}. We identify $RLf$ with $B^{\eta,\rho,0}_{R} \otimes_R (B^{\eta,\varnothing,1}_R)^{\otimes d}$ and $RN^{\rho,d}f$ with $B^{\eta,\rho,0}_R \otimes_R (B^{\eta,\varnothing,1}_R \wr_s \T^\eta_{d,R})$ via the isomorphisms (\ref{algn:ONf_ident}). Let $\ell:=(p-1)/2$.
}
\end{void}

\begin{void}\label{void:The bimodule which induces a Morita superequivalences}
{\rm \textbf{The bisupermodule which induces a Morita superequivalence.} We set $M^R$ to be the $(B^{\eta,\varnothing,1}_R,B^{\eta,\varnothing,1}_R)$-bisupermodule 
$$M^R:=\Omega_{B^{\eta,\varnothing,1}_R\otimes_R (B^{\eta,\varnothing,1}_R)^{\rm sop}}^\ell(B^{\eta,\varnothing,1}_R).$$
Then we have the $((B^{\eta,\varnothing,1}_R)^{\otimes d},(B^{\eta,\varnothing,1}_R)^{\otimes d})$-bisupermodule $(M^R)^{\boxtimes d}$. For $1\leq i\leq d$, we identify $(B^{\eta,\varnothing,1}_R)^{(i)}:=R\ws^\eta_{V_i}e_{\varnothing,1}^{(i)}$ with $B^{\eta,\varnothing,1}_R=R\ws^\eta_pe_{\varnothing,1}$ via (\ref{Sp_ident}), as in \ref{void:T_w}. Denote by $(M^R)^{(i)}$ the $i$-th factor in $(M^R)^{\boxtimes d}$. That is, $(M^R)^{(i)}$ is a $((B^{\eta,\varnothing,1}_R)^{(i)},(B^{\eta,\varnothing,1}_R)^{(i)})$-bisupermodule that we identify with $M^R$. Denote by $(M^R)_{S_d}^{\boxtimes d}$ the $(B^{\eta,\varnothing,1}_R\wr_s \T^\eta_{d,R},B^{\eta,\varnothing,1}_R\wr_s \T^\eta_{d,R})_{S_d}$-supermodule defined as in (\ref{algn:epsilon-ext_M^n}).

Recalling the identification $RLf \cong B^{\eta,\rho,0}_R \otimes_R (B^{\eta,\varnothing,1}_R)^{\otimes d}$, define the $(RLf,RLf)$-bisupermodule
$$M_L^R:=B^{\eta,\rho,0}_R\boxtimes_R (M^R)^{\boxtimes d}.$$
We have also identified $RN^{\rho,d}f$ with $B^{\eta,\rho,0}_R \otimes_R (B^{\eta,\varnothing,1}_R \wr_s \T^\eta_{d,R})$. Recalling (\ref{algn:epsilon-sup_wreath_bimod}), define the $(RN^{\rho,d}f,RN^{\rho,d}f)$-bisupermodule
\begin{align*}
	M_{N^{\rho,d}}^R  := B^{\eta,\rho,0}_R \boxtimes_R (M^R \wr_s \T^\eta_{d,R})
	=B^{\eta,\rho,0}_R \boxtimes_R \Big(\Ind_{(B^{\eta,\varnothing,1}_R \wr_s \T^\eta_{d,R},B^{\eta,\varnothing,1}_R \wr_s \T^\eta_{d,R})_{S_d}}^{(B^{\eta,\varnothing,1}_R \wr_s \T^\eta_{d,R}) \otimes_R (B^{\eta,\varnothing,1}_R \wr_s \T^\eta_{d,R})^{\rm sop}}(M^R)^{\boxtimes d}_{S_d}\Big).
\end{align*}

By \cite[Lemma 7.2 (ii)]{KL25} and Proposition \ref{proposition:epsilon-MwrT} (i), the $(RN^{\rho,d}f,RN^{\rho,d}f)$-bimodule $M_{N^{\rho,d}}^R$ is module-indecomposable with vertex $\Delta D:=\{(u,u)\mid u\in D\}$, and by \cite[Lemma 4.8 (ii)]{KL25}, we have $N_{G\times N^{\rho,d}}(\Delta D)\leq N^{\rho,d}\times N^{\rho,d}$. Hence we can now define $X^R$ to be the super Green correspondent of $M_{N^{\rho,d}}^R$ in $G\times N^{\rho,d}$; see Theorem \ref{theo:defi-of-supergreen} for the definition of the super Green correspondence. 
By \cite[Lemma 7.5]{KL25}, $X^R$ is an $(RGb,RN^{\rho,d}f)$-bisupermodule. 
}
\end{void}

\begin{theorem}[{\cite[Theorem 8.13, Corollary 8.14 and Lemma 9.2]{KL25}}]
The bisupermodules $X^R$ and $(X^R)^*$ induce a Morita superequivalence between $RGb$ and $RN^{\rho,d}f$. The bimodules $X^R_{\0}$ and $(X^R_{\0})^*$ induce a Morita equivalence between $RG_{\0}b$ and $RN_{\0}^{\rho,d}f$. Both the bimodules $X^R$ and $X^R_{\0}$ have endopermutation sources.
\end{theorem}	
	
\begin{theorem}\label{theorem: Proof of Theorem main1}
The blocks $R_p Gb$ and $R_pNf$ are Morita superequivalent. More precisely, there is an $(R_pGb,R_pN^{\rho,d}f)$-bisupermodule $X^{R_p}$, together with its dual $(X^{R_p})^*$, inducing a Morita superequivalence between $R_p Gb$ and $R_pN^{\rho,d}f$ and fulfilling $X^R\cong R\otimes_{R_p}X^{R_p}$. The bimodule $X^{R_p}$ has the endopermutation $R_p\Delta D$-source
$U^{R_p}:=\Omega_{R_p \Delta D_1}^\ell(R_p)\otimes_{R_p} \cdots\otimes_{R_p} \Omega_{R_p \Delta D_d}^\ell(R_p)$.
\end{theorem}

\begin{proof} For the first statement, by Proposition \ref{proposition:Descent of Morita superequivalences}, it suffices to show that there exists an $(R_pGb,R_pN^{\rho,d}f)$-bisupermodule $X^{R_p}$ such that $X^R\cong R\otimes_{R_p} X^{R_p}$. By Remark \ref{remark:ONf_ident-over-R0} we have an isomorphism of $R_p$-superalgebras
$$R_pN^{\rho,d} f  \cong B^{\eta,\rho,0}_{R_p} \otimes_{R_p} (B^{\eta,\varnothing,1}_{R_p} \wr_s \T^\eta_{d,R_p}).$$ We let $M^{R_p}$ be the $(B^{\eta,\varnothing,1}_{R_p},B^{\eta,\varnothing,1}_{R_p})$-bimodule 
$$M^{R_p}:=\Omega_{B^{\eta,\varnothing,1}_{R_p}\otimes_{R_p} (B^{\eta,\varnothing,1}_{R_p})^{\rm sop}}^\ell(B^{\eta,\varnothing,1}_{R_p}).$$
Then by Proposition \ref{proposition:Descent of Heller translates} we have $M^R\cong R\otimes_{R_p} M^{R_p}$. Define the $(R_pN^{\rho,d}f,R_pN^{\rho,d}f)$-bisupermodule
\begin{align*}
\begin{split}
	M_{N^{\rho,d}}^{R_p}  &:= B^{\eta,\rho,0}_{R_p} \boxtimes_{R_p} (M^{R_p} \wr_s \T^\eta_{d,R_p})\\
	& =B^{\eta,\rho,0}_{R_p} \boxtimes_{R_p} \Big(\Ind_{(B^{\eta,\varnothing,1}_{R_p} \wr_s \T^\eta_{d,R_p},B^{\eta,\varnothing,1}_{R_p} \wr_s \T^\eta_{d,R_p})_{S_d}}^{(B^{\eta,\varnothing,1}_{R_p} \wr_s \T^\eta_{d,R_p}) \otimes_{R_p} (B^{\eta,\varnothing,1}_{R_p} \wr_s \T^\eta_{d,R_p})^{\rm sop}}(M^{R_p})^{\boxtimes d}_{S_d}\Big).
	\end{split}
	\end{align*}
Then by Proposition \ref{proposition:epsilon-MwrT}, we have $M_{N^{\rho,d}}^R\cong R\otimes_{R_p} M_{N^{\rho,d}}^{R_p}$ as $(R N^{\rho,d}f,R N^{\rho,d}f)$-bisupermodules. Now by Proposition \ref{prop:descent super Green cor}, there is an $(R_pGb,R_pN^{\rho,d}f)$-bisupermodule $X^{R_p}$ such that $X^R\cong R\otimes_{R_p} X^{R_p}$, as desired. 

Denote by $U^{R_p}$ the $R_p\Delta D$-module $\Omega_{R_p \Delta D_1}^\ell(R_p)\otimes_{R_p} \cdots\otimes_{R_p} \Omega_{R_p \Delta D_d}^\ell(R_p)$. By Proposition \ref{proposition:Descent of Heller translates}, we have
$$R\otimes_{R_p} U^{R_p}\cong\Omega_{R \Delta D_1}^\ell(R)\otimes_{R} \cdots\otimes_{R} \Omega_{R \Delta D_d}^\ell(R)$$
as $R\Delta D$-modules. Then by \cite[Lemma 9.2]{KL25}, $R\otimes_{R_p} U^{R_p}$ is a source of $X^R$. Hence by \cite[Lemma 5.2]{Kessar_Linckelmann}, $X^{R_p}$ has $U^{R_p}$ as a source. \end{proof}

% By \cite[Lemma 3.10]{KL25}, Theorem \ref{theorem: Proof of Theorem main1} has the following corollary:

\begin{corollary}\label{corollary:Morita equ for wa_n}
With the notation of Theorem \ref{theorem: Proof of Theorem main1}, the $(R_pG_{\0}b,R_pN_{\0}^{\rho,d}f)$-bimodule $X_{\0}^{R_p}$ induces a Morita equivalence between $R_pG_{\0}b$ and $R_pN_{\0}^{\rho,d}f$ and has an endopermutation source.
\end{corollary}

\begin{proof} By Theorem \ref{theorem: Proof of Theorem main1} and \cite[Lemma 3.10]{KL25}, $X_{\0}^{R_p}$ induces a Morita equivalence between $R_pG_{\0}b$ and $R_pN_{\0}^{\rho,d}f$. Since ${\rm Res}_{G_{\0}\times N_{\0}^{\rho,d}}^{G\times N^{\rho,d}}(X^{R_p})=X_{\0}^{R_p}\oplus X^{R_p}_{\1}$, $X^{R_p}_{\0}$ is isomorphic to a direct summand of
$${\rm Res}_{G_{\0}\times N_{\0}^{\rho,d}}^{G\times N^{\rho,d}}{\rm Ind}_{\Delta D}^{G\times N^{\rho,d}}(U^{R_p}).$$
By the Mackey formula,  $X_{\0}^{R_p}$ has an endopermutation source.
      \end{proof}

\begin{remark}
{\rm In \cite[Remark 8.15]{KL25}, Kleshchev and Livesey asked that if $K$ only contains a primitive $n!$-th root of unity (and not necessarily a primitive $(2n!)$-th root of unity), can we still construct a Morita equivalence between $\O G_{\0}$ and $\O N^{\rho,d}_{\0}f$? Corollary \ref{corollary:Morita equ for wa_n} answers this question with the highest generality.}
\end{remark}

\begin{theorem}\label{theorem: Rock splendid equ over Zp}
With the notation of Theorem \ref{theorem: Proof of Theorem main1}, $\Z_pGb$ is splendidly Rickard equivalent to $\Z_p N^{\rho,d}f$.
\end{theorem}

\begin{proof} 
	
By Theorem \ref{theorem: Proof of Theorem main1}, the bimodule $X^{\Z_p}$ induces a Morita equivalence between $\Z_p Gb$ and $\Z_pN^{\rho,d}f$ and has the endopermutation $\Z_p\Delta D$-source $U^{\Z_p}=\Omega_{\Z_p \Delta D_1}^\ell(\Z_p)\otimes_{\Z_p} \cdots\otimes_{\Z_p} \Omega_{\Z_p \Delta D_d}^\ell(\Z_p)$. 
We can regard $U^{\Z_p}$ as an endopermutation $\Z_p D$-module via the isomorphism $D\cong\Delta D$. Note that for any $1\leq m\leq d$, since $D_m$ is a cyclic group of order $p$, we have
$$\Omega_{\Z_p  D_m}^\ell(\Z_p)\cong \left\{ \begin{gathered}
	\Z_p,~{\rm if}~\ell~{\rm is~even}, \hfill \\
\Omega_{\Z_p D_m}(\Z_p),~{\rm if}~\ell~{\rm is~odd}. \hfill \\
\end{gathered} \right. $$
By \cite[7.11.7 and 7.11.6 (ii)]{Lin18b}, the $\Z_p D$-source $U^{\Z_p}$ has an endosplit $p$-permutation resolution
$$\Y^{\Z_p}:= \left\{ \begin{gathered}
	\Z_p,~{\rm if}~\ell~{\rm is~even}, \hfill \\
	\Y_1^{\Z_p}\otimes_{\Z_p}\cdots\otimes_{\Z_p}\Y_d^{\Z_p},~{\rm if}~\ell~{\rm is~odd}. \hfill \\
\end{gathered} \right. $$
where
$$\Y_m^{\Z_p}:=\cdots\to \Z_p D_m\to \Z_p\to 0\to\cdots~~~(1\leq m\leq d)$$
with $\Z_p D_m$ in degree $0$ and the non-zero boundary map $\Z_p D_m\to \Z_p$ being given by $u\mapsto 1$ for all $u\in D_m$.  In the proof of \cite[Lemma 9.3]{KL25}, it is shown that the complex $\Y^\O:=\O\otimes_{\Z_p}\Y^{\Z_p}$ is $N_G(D)$-stable. One easily note that the argument there does not depend on the coefficient ring, and hence $\Y^{\Z_p}$ is $N_G(D)$-stable as well. By \cite[Proposition 3.3]{H26}, the $(\Z_pGb,\Z_pN^{\rho,d}f)$-bimodule $X^{\Z_p}$ is isomorphic to a direct summand of
$$\Z_p Gi\otimes_{\Z_p D}{\rm Ind}_{\Delta D}^{D\times D}(U^{\Z_p})\otimes_{\Z_pD} j\Z_p N^{\rho,d}$$
for some source idempotents $i\in (\Z_p Gb)^D$ and $j\in (\Z_pN^{\rho,d}f)^D$. Let $(D,f_D)$ be the maximal $\Z_pN^{\rho,d}f$-Brauer pair determined by $j$ and let $\F:=\F_{(D,f_D)}(\Z_pN^{\rho,d}f)$. Since $N_{N^{\rho,d}}(D)\leq N_G(D)$, $\Y^{\Z_p}$ is $N_{N^{\rho,d}}(D)$-stable. Since $D$ is abelian, by Remark \ref{remark:extension of fusion} we see that $\Y^{\Z_p}$ is $\F$-stable. Now by Proposition \ref{prop:9.11.4+9.11.5} (iv), $\Z_pGb$ and $\Z_p N^{\rho,d}f$ are splendidly Rickard equivalent.
\end{proof}

\begin{theorem}\label{theorem: RoCK splendid over Zp for alternating double cover}
The block $\Z_pG_{\0}b$ is splendidly Rickard equivalent to $\Z_p N^{\rho,d}_{\0}f$.
\end{theorem}

\begin{proof} Since $p\nmid |G\times N^{\rho,d}:G_{\0}\times N_{\0}^{\rho,d}|$, any $p$-subgroup of $G\times N^{\rho,d}$ is contained in $G_{\0}\times N^{\rho,d}_{\0}$. Since ${\rm Res}_{G_{\0}\times N_{\0}^{\rho,d}}^{G\times N^{\rho,d}}(X^{\Z_p})=X_{\0}^{\Z_p}\oplus X^{\Z_p}_{\1}$, $X^{\Z_p}_{\0}$ is isomorphic to a direct summand of 
$${\rm Res}_{G_{\0}\times N_{\0}^{\rho,d}}^{G\times N^{\rho,d}}{\rm Ind}_{\Delta D}^{G\times N^{\rho,d}}(U^{\Z_p})\cong \bigoplus_{t\in G\times N^{\rho,d}/G_{\0}\times N_{\0}^{\rho,d}}{\rm Ind}_{{}^t (\Delta D)}^{G_{\0}\times N^{\rho,d}_{\0}}({}^tU^{\Z_p}),$$   
where the isomorphism is by the Mackey formula. Again since $p\nmid |G\times N^{\rho,d}:G_{\0}\times N_{\0}^{\rho,d}|$, any vertex of the bimodule $X_{\0}^{\Z_p}$ has the same order with any vertex of the bimodule $X^{\Z_p}$.  It follows that $X_{\0}^{\Z_p}$ has an endopermutation $\Z_p ({}^t(\Delta D))$-source ${}^tU^{\Z_p}$ for some $t\in  G\times N^{\rho,d}/G_{\0}\times N_{\0}^{\rho,d}$. Write $t=(t_1,t_2)$ with $t_1\in G$ and $t_2\in N^{\rho,d}$; then ${}^t(\Delta D)=\{({}^{t_1}u,{}^{t_2}u)\mid u\in D\}$. We can regard $U^{\Z_p}$ as a $\Z_pD$-module via the isomorphism $D\cong \Delta D$ and view ${}^tU^{\Z_p}$ as a $\Z_p {}^{t_2}D$-module via the isomorphism ${}^{t_2}D\cong {}^t(\Delta D)$. Since $\Y^{\Z_p}$ is an endosplit $p$-permutation resolution of the $\Z_pD$-module $U^{\Z_p}$, ${}^{t_2}\Y^{\Z_p}$ is an endosplit $p$-permutation resolution of the $\Z_p{}^{t_2}D$-module ${}^tU^{\Z_p}$. Since $\Y^{\Z_p}$ is $N_{N^{\rho,d}}(D)$-stable, we see that ${}^{t}\Y^{\Z_p}$ is $N_{N^{\rho,d}_{\0}}({}^{t_2}D)$-stable. Write $P:={}^{t_1}D$, $Q:={}^{t_2}D$ and $\varphi:P\cong Q$ the group isomorphism induced by conjugation action of $t_2t_1^{-1}$. By \cite[Proposition 3.2]{H26}, $P$ is a defect group of $\Z_p G_{\0}b$ and $Q$ is a defect group of $\Z_p N^{\rho,d}_{\0}f$. By \cite[Proposition 3.3]{H26} 
$X_{\0}^{\Z_p}$ is isomorphic to a direct summand of
$$\Z_p G_{\0}i_0\otimes_{\Z_p P}{\rm Ind}_{\Delta \varphi}^{P\times Q}({}^{t}U^{\Z_p})\otimes_{\Z_p Q} j_0\Z_p N^{\rho,d}_{\0}$$
for some source idempotents $i_0\in(\Z_p G_{\0}b)^P$ and $j_0\in(\Z_p N^{\rho,d}_{\0}f)^Q$. The rest is analogous with the proof of Theorem \ref{theorem: Rock splendid equ over Zp}:  Let $(Q,f_Q)$ be the maximal $\Z_pN^{\rho,d}_{\0} f$-Brauer pair determined by $j_0$ and let $\F_{\0}:=\F_{(Q,f_Q)}(\Z_pN^{\rho,d}_{\0}f)$. Since $Q$ is abelian, by Remark \ref{remark:extension of fusion} we see that ${}^t\Y^{\Z_p}$ is $\F_0$-stable. Now by Proposition \ref{prop:9.11.4+9.11.5} (iv), $\Z_pG_{\0}b$ and $\Z_p N^{\rho,d}_{\0}f$ are splendidly Rickard equivalent.   \end{proof}

\begin{remark}
{\rm When $p\equiv 1~{\rm mod}~4$ (i.e. $\ell$ is even), the $\Z_p D$-source $U^{\Z_p}$ of $X^{\Z_p}$ is the trivial module. Hence in this case $X^{\Z_p}$ induces a splendid Morita equivalence between $\Z_p Gb$ and $\Z_p N^{\rho,d}f$, and $X^{\Z_p}_{\0}$ induces a splendid Morita equivalence between $\Z_p G_{\0}b$ and $\Z_p N^{\rho,d}_{\0}f$. %So by Theorem \ref{theorem: Proof of Theorem main1} we immediately have the following:
}
\end{remark}

%\begin{theorem}
%If $p\equiv 1~{\rm mod}~4$, then $X^{\Z_p}$ induces a splendid Morita equivalence between $\Z_p Gb$ and $\Z_p N^{\rho,d}d$, and $X^{\Z_p}_{\0}$ induces a splendid Morita equivalence between $\Z_p G_{\0}b$ and $\Z_p N^{\rho,d}_{\0}f$.
%\end{theorem}	

We can now construct splendid Rickard equivalences between $\Z_p Gb$ (resp. $\Z_pG_{\0}b$) and its Brauer correspondent in $
\Z_pN_G(D)$ (resp. $\Z_pN_{G_{\0}}(D)$):

\begin{theorem}\label{theo: proof of main2}
The block $R_pGb$ (resp. $R_pG_{\0}b$) is splendidly Rickard equivalent to its Brauer correspondent in $R_pN_G(D)$ (resp. $R_pN_{G_{\0}}(D)$).
	\end{theorem}
\begin{proof}	
 By the lifting theorem of splendid Rickard equivalences (Theorem \ref{theorem: lifting theorem of splendid Rickard complexes}), it suffices to consider $R_p=\FF_p$. By Theorem \ref{theorem:Brauer corr over Fp splendid}, there is a splendid Rickard equivalence between the Brauer correspondent of $\FF_pGb$ (resp. $\FF_pG_{\0}b$) in $\FF_pN_G(D)$ (resp. $\FF_pN_{G_0}(D)$) and  $\FF_pN^{\rho,d}f$ (resp. $\FF_pN^{\rho,d}_{\0}f$). Now the statements follows from Theorems \ref{theorem: Rock splendid equ over Zp}, \ref{theorem: RoCK splendid over Zp for alternating double cover} and Proposition \ref{prop:tensor products and duality}. 
\end{proof}

\bigskip\noindent\textbf{Acknowledgements.}\quad We thank Prof. Pham Huu Tiep for suggesting us to consider spin blocks of symmetric and alternating groups in SUSTech, May 2024, and thank Prof. Alexander Kleshchev for helpful explanation. We are grateful to the referee for careful reading and many constructive suggestions.

\bigskip\noindent\textbf{Funding information.} Y. Du is supported by National Key R \& D Program of China (2020YFE0204200) and NSFC (12601033, 12526574, 12431001, 12350710787). X. Huang is supported by NSFC (12501024, 12471016), Fundamental Research Funds for the Central Universities (CCNU24XJ028), China Postdoctoral Science Foundation (GZC20252006, 2025T001HB) and China Scholarship Council (202506770066).
%CCNU25JC025, CCNUJCPT031


\begin{thebibliography}{1}
\expandafter\ifx\csname url\endcsname\relax
  \def\url#1{\texttt{#1}}\fi
\expandafter\ifx\csname urlprefix\endcsname\relax\def\urlprefix{URL }\fi
\expandafter\ifx\csname href\endcsname\relax
  \def\href#1#2{#2} \def\path#1{#1}\fi



\bibitem{Bessenrodt18}
C. Bessenrodt, Critical classes, Kronecker products of spin characters, and the Saxl conjecture, Algebr. Comb. {\bf 1} (2018) 353--369.

\bibitem{Boltje}
R. Boltje, The Brou\'{e} invariant of a $p$-permutation equivalence, Ann. Represent. Theory {\bf 1} (2024) 517--537.

\bibitem{BP}
R. Boltje, P. Perepelitsky, $p$-Permutation equivalences between blocks of group algebras, J. Algebra {\bf 664} (2025) 815--887.

%\bibitem{BKL}
%R. Boltje, R. Kessar, M. Linckelmann, On Picard groups of blocks of finite groups, J. Algebra {\bf 558} (2020) 70--101.

\bibitem{BY}
R. Boltje, D. Y{\i}lmaz, Galois descent of equivalences between blocks of $p$-nilpotent groups, Proc. Amer. Math. Soc. {\bf150} (2022) 559--573.

%\bibitem{Broue85}
%M. Brou\'{e}, On Scott modules and $p$-permutation modules: an approach through the Brauer morphism, Proc. Amer. Math. Soc. {\bf 93} (1985) 401--408.

\bibitem{BM}
J. V. Breland, S. K. Miller, Brauer pairs for splendid Rickard equivalences, J. Algebra {\bf 691} (2026) 694--729.

\bibitem{BK}
J. Brundan and A. Kleshchev, Odd Grassmannian bimodules and derived equivalences for spin symmetric groups, Selecta Math. (N.S.) {\bf 32} (2026) Paper No. 9.

\bibitem{Ca}
M. Cabanes, Local structure of the $p$-blocks of $\tilde{S}_n$, Math. Z. {\bf 198} (1988) 519--543.

\bibitem{CK}
J. Chuang, R. Kessar, Symmetric groups, wreath products, Morita equivalences, and Brou\'{e}'s abelian defect group conjecture, Bull. London Math. Soc. {\bf 34} (2002) 174--184.

\bibitem{CR}
J. Chuang, R. Rouquier, Derived equivalences for symmetric groups and $\mathfrak{sl}_2$-categorification, Ann. Math. {\bf 167} (2008) 245--298.

\bibitem{Da84}
E. C. Dade, Extending group modules in a relatively prime case, Math. Z. {\bf 186} (1984) 81--98.


\bibitem{DH}
Y. Du, X. Huang, J. Zhang, The blockwise Navarro Alperin Weight Conjecture for double covers of symmetric and alternating groups, arXiv:2509.13673. 

\bibitem{ELV}
M. Ebert, A. D. Lauda, L. Vera, Derived superequivalences for spin symmetric groups and odd $\mathfrak{sl}_2$-categorifications, arXiv:2203.14153.

\bibitem{FH}
P. Fong, M.E. Harris, On perfect isometries and isotypies in alternating groups, Trans, Amer. Math. Soc. {\bf 349} (1997) 3469--3516.

%\bibitem{Hi-Ho}
%T. Hirai, A. Hora, Spin representations of twisted central products of double covering finite groups and the case
%of permutation groups, J. Math. Soc. Japan {\bf 66} (2014) 1191--1226.

\bibitem{HL}
M. E. Harris, M. Linckelmann, Splendid derived equivalences for blocks of finite $p$-solvable groups, J. London Math. Soc. (2) {\bf 62} (2000) 85--96.

\bibitem{H23}
X. Huang, Descent of equivalences for blocks with Klein four defect groups, J. Algebra {\bf 614} (2023) 898--905. 

\bibitem{H24}
X. Huang, Descent of splendid Rickard equivalences in alternating groups, J. Algebra {\bf 658} (2024) 277--293. 

%\bibitem{H24c}
%X. Huang, Stable equivalences with endopermutation source, slash functors, and funtorial equivalences, J. Algebra {\bf 650} (2024) 335--353.

\bibitem{H24b}
X. Huang, Virtual Morita equivalences and Brauer character bijections, Arch. Math. {\bf 123} (2024) 117--121.

\bibitem{H26}
X. Huang, On stable equivalences of Morita type with twisted diagonal vertices, Accepted by Represent. Theory (2026).

\bibitem{HLZ23}
X. Huang, P. Li, J. Zhang, The strengthened Brou\'{e} abelian defect group conjecture for ${\rm SL}(2,p^n)$ and ${\rm GL}(2,p^n)$, J. Algebra {\bf 633} (2023) 114--137.

\bibitem{HLZ24}
X. Huang, P. Li, J. Zhang, Descent of splendid Rickard equivalences in ${\rm GL}_n(q)$, J. Pure Appl. Algebra {\bf 228} (2024) 107711.

\bibitem{Hum}
J. F. Humphreys, Blocks of projective representations of symmetric groups, J. London Math. Soc. (2) {\bf 33} (1986) 441--452. 

\bibitem{K96}
R. Kessar, Blocks and source algebras for the double covers of the symmetric and alternating groups, J. Algebra {\bf 186} (1996) 872--933.

\bibitem{Kessar_Linckelmann}
R. Kessar, M. Linckelmann, Descent of equivalences and character bijections, Geometric and topological aspects of the representation theory of finite groups, Springer Proc. Math. Stat., {\bf 242} (2018) 181--212.

\bibitem{KS}
R. Kessar, M. Schaps, Crossover Morita equivalences for blocks of the covering groups of the symmetric and alternating groups, J. Group Theory {\bf 9} (2006) 715--730.


\bibitem{KL22}
A. Kleshchev, M. Livesey, RoCK blocks for double covers of the symmetric groups and quiver Hecke superalgebras, Mem. Amer. Math. Soc. {\bf 309} (2025) no. 1564.

\bibitem{KL25} 
A. Kleshchev, M. Livesey, RoCK blocks for double covers of symmetric groups over a complete discrete valuation ring, Ann. Represent. Theory {\bf2} (2025) 85--171.



\bibitem{Lin18a}
M. Linckelmann, The Block Theory of Finite Group Algebras I, London Math. Soc. Student Texts, vol. {\bf91}, Cambridge University Press, 2018.

\bibitem{Lin18b}
M. Linckelmann, The Block Theory of Finite Group Algebras II, London Math. Soc. Student Texts, vol. {\bf92}, Cambridge University Press, 2018.

%\bibitem{MarJA}
%A. Mǎrcuş, On equivalences between blocks of group algebras: reduction to simple components, J. Algebra {\bf 184} (1996) 372--396.

\bibitem{MarCA}
A. Mǎrcuş, Equivalences induced by graded bimodules, Comm. Algebra {\bf 26} (1998)  713--731.

\bibitem{M1}
S. K. Miller, Galois descent of splendid Rickard equivalences for blocks of $p$-nilpotent groups, Proc. Amer. Math. Soc. {\bf 153} (2025) 1893--1902.

\bibitem{M2}
S. K. Miller, On endosplit $p$-permutation resolutions and Brou\'{e}'s conjecture for $p$-solvable groups, Represent. Theory {\bf 30} (2026) 151--174.

\bibitem{Morris}
A. O. Morris, The spin representations of the symmetric group, Can. J. Math. {\bf 17} (1965) 543--549.

\bibitem{NT}
H. Nagao, Y. Tsushima, Representations of Finite Groups, Academic Press, 1989.

\bibitem{Rickard96}
J. Rickard, Splendid equivalences: derived categories and permutation modules, Proc. Lond. Math. Soc. {\bf 72} (1996) 331--358.

%\bibitem{Rou}
%R. Rouquier, The derived category of blocks with cyclic defect groups, in: Derived Equivalences
%for Group Rings (S. K\"{o}nig, A. Zimmermann), Lecture Notes in Math. {\bf1685},
%Springer Verlag, Berlin-Heidelberg, 1998, 199--220.


\bibitem{Sch}
J. Schur, \"{U}ber die darstellung der symmetrischen und der alternierenden Gruppe durch gebrochene lineare substitutionen, J. Reine. Angew. Math. {\bf 139} (1911) 155--250.




%\bibitem{Rickard}
%J. Rickard, Splendid equivalences: derived categories and permutation modules, Proc. Lond. Math. Soc. {\bf 72} (1996) 331--358.



%\bibitem{Thevenaz}
%J. Th\'{e}venaz, $G$-algebras and Modular Representation Theory, Oxford Science Publications, Clarendon, Oxford, 1995.

%\bibitem{Th07}
%J. Th\'{e}venaz, Endo-permutation modules, a guided tour,  in: M. Geck, D. Testerman, J. Th\'{e}venaz (Eds.), Group Representation Theory, EPFL Press, Lausanne, 2007, pp.115--147. 


\end{thebibliography}
\end{document}